\documentclass[11pt,a4paper]{article}
\IfFileExists{lmodern.sty}{\usepackage[T1]{fontenc}\usepackage{lmodern}}{}
\usepackage{microtype}
\usepackage[margin=25mm]{geometry}
\usepackage{amsmath,amssymb,amsthm,mathtools,bm}
\usepackage{enumitem,booktabs,graphicx}
\usepackage[font=small,labelfont=bf,labelsep=period]{caption}
\usepackage[round,authoryear]{natbib}
\usepackage[dvipsnames]{xcolor}
\usepackage[colorlinks=true,linkcolor=MidnightBlue,citecolor=MidnightBlue,urlcolor=MidnightBlue]{hyperref}
\usepackage[capitalise,nameinlink]{cleveref}
\usepackage{setspace}
\usepackage{placeins}
\hypersetup{pdftitle={Finite-horizon Optimal Consumption and Investment with Costly Income Adjustment},
            pdfauthor={Gugyum Ha, Junkee Jeon}}
\numberwithin{equation}{section}
\newtheorem{theorem}{Theorem}[section]
\newtheorem{proposition}[theorem]{Proposition}
\newtheorem{lemma}[theorem]{Lemma}
\newtheorem{corollary}[theorem]{Corollary}
\theoremstyle{definition}
\newtheorem{assumption}[theorem]{Assumption}

\theoremstyle{remark}
\newtheorem{remark}[theorem]{Remark}
\crefname{assumption}{Assumption}{Assumptions}

\newcommand{\E}{\mathbb E}
\newcommand{\R}{\mathbb R}
\newcommand{\X}{\mathcal X}
\newcommand{\PP}{\mathcal P}
\newcommand{\LL}{\mathcal L}
\newcommand{\ind}{\mathbf 1}
\newcommand{\dd}{\,\mathrm d}
\newcommand{\ee}{\mathrm e}
\newcommand{\Zr}{\mathcal Z}
\newcommand{\Nr}{\mathcal N}
\newcommand{\Ar}{\mathcal A}
\newcommand{\adm}{\mathfrak A}
\newcommand{\zb}{\bar z}
\DeclareMathOperator{\sgn}{sgn}
\makeatletter
\def\@fnsymbol#1{\ensuremath{\ifcase#1\or *\or \dagger\or \ddagger\or \mathsection\or \mathparagraph\or \|\or **\or \dagger\dagger\or \ddagger\ddagger\else\@ctrerr\fi}}
\makeatother

\title{\bf \LARGE Finite-horizon Optimal Consumption and Investment with Costly Income Adjustment}
\author{Gugyum Ha\thanks{Department of Mathematics, Sogang University, Seoul 04107, Republic of Korea. Email: \href{mailto:ggha@sogang.ac.kr}{ggha@sogang.ac.kr}.}
\and Junkee Jeon\thanks{Department of Applied Mathematics, Kyung Hee University, Yongin-si 17104, Republic of Korea. Email: \href{mailto:junkeejeon@khu.ac.kr}{junkeejeon@khu.ac.kr}.}}
\date{October 2026}

\begin{document}
\maketitle

\begin{abstract}
\noindent
We study optimal consumption and investment for an investor who receives labor income until a fixed retirement date, can raise the income rate irreversibly at a constant price per unit of income flow, and cannot borrow against future income. A pathwise argument shows that no income is purchased during a final period before retirement whose length depends only on the interest rate and the price. Before this period the dual problem has two free boundaries: a lower one generated by the borrowing constraint and an upper one, at which income is purchased, that escapes to infinity as the final period begins. The analysis of this problem is our main mathematical contribution. During the final period it is a single obstacle problem, which we solve, together with the regularity of its free boundary, by established methods. Before that period a partial hodograph transformation straightens both free boundaries and turns the problem into a quasilinear parabolic equation with oblique boundary conditions on a fixed interval. Comparison principles, explicit barriers and regularity theory for quasilinear equations then give existence, uniqueness and smoothness. Transforming back and gluing the two constructions yields global regularity: both free boundaries are infinitely differentiable and the value function is a classical solution at every positive wealth level. A dual budget process verifies optimality. We also derive coupled Volterra integral equations for the free boundaries, prove exponential convergence to the perpetual problem under a scalar condition, and find numerically a purchase threshold that is U-shaped in the time to retirement.
\end{abstract}



\section{Introduction}\label{sec:intro}

Optimal consumption and investment in the presence of labor income is a basic problem of continuous-time finance. When the investor cannot borrow against future income, wealth must stay nonnegative, and the problem becomes a stochastic control problem with a state constraint whose solution is described by a free boundary at zero wealth \citep{HePages1993,ElKarouiJeanblanc1998,DybvigLiu2010}. In these models the income stream is exogenous. In practice an investor can often raise her future income by spending wealth, for example on degrees, licenses and training, or by giving up current pay to move into a better paid position. Such an investment is irreversible and pays off only until retirement, so its value falls as retirement approaches and the horizon is finite in an essential way. It must also be financed out of liquid wealth, so it competes with the precautionary role of savings under the borrowing constraint. This paper studies a continuous-time model in which both features are present and derives the optimal consumption, portfolio and income adjustment policies together with a detailed description of the associated free boundaries.

We consider an investor with constant relative risk aversion who trades a riskless bond and a stock, consumes, and receives an income flow that she may raise at any time before a fixed retirement date, irreversibly, at a fixed price per unit of income flow. Wealth must remain nonnegative. Income stops at retirement, the income stock cannot be resold, and after retirement the investor solves the Merton problem with the same preferences. The value function therefore has the Merton value as terminal condition. Income adjustment is a singular control of monotone follower type, and the purchase decision depends on calendar time through the time remaining until retirement.

The first result is elementary and exact. An additional unit of income is worth the present value of the income it generates until retirement, and this value falls to zero as retirement approaches. Once it is below the price, a pathwise comparison shows that any strategy that buys income is weakly dominated by the strategy that keeps the same consumption and portfolio and never buys, and the comparison respects the borrowing constraint path by path. Hence nothing is purchased during a final window before retirement, the period in which the present value of the remaining income is below the price. The length of this window depends only on the interest rate and the price, not on risk aversion, impatience or the equity premium. We show that the window is exact: before it, the purchase region is nonempty. Without the borrowing constraint the problem would not be well posed before the window, because an investor who can borrow against future income could buy income at its price and borrow more than the price against it (\cref{rem:natural}).

Before the window the structure is richer. After the standard dual transformation, the value function is described by the inverse of the consumption policy, which gives the wealth--income ratio at which a given consumption--income ratio is optimal. The borrowing constraint produces a lower free boundary in the consumption--income ratio, below which wealth is zero, and the option to buy income produces an upper free boundary, above which income is bought. As calendar time approaches the start of the window, the upper boundary escapes to infinity. In the time to retirement, which is the natural time variable of the dual problem, the upper boundary therefore enters from infinity at the length of the window: at this time the problem changes from one free boundary to two, and the new boundary has no finite initial position.

The main mathematical contribution of the paper is the analysis of this free boundary problem, which differs from the approaches commonly used in mathematical finance. These approaches either postulate smooth free boundaries and verify them afterwards, or obtain the regularity of a single free boundary by methods specific to obstacle problems, as for the American put \citep{ChenChadam2007}. Neither applies once two boundaries are present. Where the upper boundary enters, the problem has no natural front fixing formulation, and a direct construction between the two moving boundaries would require regularity of the upper boundary that has not yet been established. We proceed instead in three steps. During the window only the zero-wealth boundary is present, and the problem is a parabolic obstacle problem with a single free boundary. We solve it by penalization and prove that its free boundary is smooth with the higher order boundary Harnack inequality of \citet{Kukuljan2022}, which is the established route. Before the window we straighten the two free boundaries simultaneously by a partial hodograph transformation in the sense of \citet{KinderlehrerNirenberg1977}. The transformed problem is posed on a fixed interval, with the two free boundaries on its two sides, and it is a quasilinear parabolic equation with oblique boundary conditions. General results for such equations therefore become available: a comparison principle, which makes barrier arguments possible, and H\"older gradient and Schauder estimates. Finally, we invert the transformation and glue the resulting solution to the solution of the obstacle problem along the time at which the window starts. The gluing is smooth to every order, and it gives the regularity of the value function and of both free boundaries on the whole horizon. The hodograph transformation thus serves as the construction itself rather than as a device for improving the regularity of a free boundary already known to exist, and no regularity of the free boundaries is assumed at any stage.

The transformation acts on the Euler derivative of the inverse consumption map, a function that plays the role of effective human wealth per unit of income. It solves a linear parabolic equation with a constant source, vanishes at the lower boundary, equals the price at the upper boundary, and is strictly increasing in between. Inverting it in the spatial variable gives a problem on an interval whose length is the price. The free boundaries become the boundary values of the new unknown, the free boundary conditions become an oblique nonlinear condition at one end of the interval and a Neumann condition at the other, and the escape of the upper boundary to infinity becomes an infinite initial value at one corner. On this interval we construct the solution with three explicit barriers. The decisive one solves the differentiated equation exactly, because its diffusion term and its quadratic gradient term cancel. It gives uniform parabolicity away from the corner and a gradient bound that does not deteriorate as the truncation of the initial datum is removed. An energy estimate, based on an exact algebraic identity for the nonlinear boundary flux, gives convergence of the truncations and uniqueness. As a result, both free boundaries are infinitely differentiable, the upper boundary tends to infinity at the start of the window, and the value function is a classical solution at every strictly positive wealth level. A comparison with the frozen-income problem, in which income cannot be adjusted, shows that the upper boundary blows up only mildly: its logarithm grows at most like the square root of the logarithm of the inverse of the time left before the window starts.

We also show that the zero-wealth boundary vanishes linearly at retirement, with slope equal to the Merton consumption--wealth ratio. At zero wealth the value function is not twice differentiable. Its curvature blows up like the inverse square root of wealth, the marginal propensity to consume out of wealth is infinite, and the risky position grows like the square root of wealth, with explicit coefficients. Optimality of the reflected feedback controls is proved with a positive dual budget process built from two shadow prices, one for wealth and one for income. The argument uses only nonnegativity and the supermartingale property. It covers relative risk aversion below and above one and needs no integrability condition on competing strategies.

The optimal strategy keeps the ratio of liquid wealth to income between zero and a threshold and buys income whenever the ratio reaches the threshold. The threshold is infinite at the start of the window and falls steeply as the time to retirement lengthens. For the baseline parameters it then rises slowly, so that it is U-shaped in the time to retirement, while the threshold measured in years of consumption falls monotonically. Simulated income paths are concave and front loaded, the shape that \citet{BenPorath1967} derived for optimal human capital accumulation: income grows early, flattens and stays constant during the final window, and wealth accumulation for retirement starts in earnest only after purchases stop. In our numerical experiments, more risk averse investors buy income at lower wealth, since income is riskless in the model, whereas a larger market price of risk raises the threshold, since liquid wealth gives access to the equity premium. The optimal risky position is the Merton fraction of liquid wealth plus an effective human wealth. The investor thus treats part of her future income as a bond, but this effective human wealth never exceeds the purchase price of the income and is well below its market value long before retirement. Consumption is also continuous at the retirement date even though income drops to zero: the investor enters retirement with exactly the wealth that the Merton rule requires for her current consumption.

We also describe how the solution depends on the time to retirement. A pathwise argument shows that the value function increases with it. In the purchase region the value is a power function of liquid wealth plus the purchase value of income, with a coefficient that depends only on the time to retirement. The monotonicity of the value therefore forces this coefficient to move in opposite directions according as relative risk aversion is above or below one. A longer time to retirement raises consumption in the purchase region if the income effect dominates and lowers it if the substitution effect dominates. When relative risk aversion exceeds one, a maximum principle for the time derivative of the inverse consumption map shows that consumption increases with the time to retirement at every wealth level and that the zero-wealth boundary is monotone. Otherwise it need not be. Far from retirement the problem approaches the perpetual income problem, in which there is no retirement date. The stationary hodograph equation reduces to a first order ordinary differential equation with a terminal condition, whose solution determines the long-run behavior. A stationary solution exists if and only if a scalar inequality on the initial value of this solution holds, and then the time dependent free boundaries converge to it exponentially fast, because the time derivative of the hodograph solves a linear equation with a dissipative boundary condition. When the inequality fails, which requires relative risk aversion below one, the consumption--income ratio at zero wealth is not bounded away from zero. The same differential equation shows that in the long run the purchase option lowers consumption at zero wealth whenever the stationary problem is solvable.

Finally, we derive integral equations for the two free boundaries. In logarithmic variables the Euler derivative solves, in the sense of distributions on the whole line, a linear parabolic equation whose right-hand side consists of the indicator functions of the regions and of two single layers carried by the free boundaries. The densities of the layers are fixed by the second order conditions at the boundaries. Duhamel's formula then expresses the Euler derivative through the two boundaries alone, and its prescribed values on the boundaries become two coupled nonlinear Volterra equations with explicit Gaussian kernels. This is an analogue of the early exercise premium representation of the American put \citep{Kim1990,CarrJarrowMyneni1992}. In addition, the coefficient of the purchase region has a closed form as an integral over the upper boundary. The numerical results in the paper are obtained by solving these equations.

\paragraph{Related literature.}
Optimal consumption and portfolio choice goes back to \citet{Merton1969,Merton1971}. The martingale and duality methods we use were developed by \citet{KaratzasLehoczkyShreve1987}, \citet{CoxHuang1989} and, for constrained problems, \citet{CvitanicKaratzas1992}. Labor income with a borrowing constraint was studied by \citet{HePages1993} and \citet{ElKarouiJeanblanc1998}, who related the constraint to an optimal stopping problem. See also \citet{Koo1998} for uninsurable income risk, \citet{DybvigLiu2010} for retirement with constrained borrowing, and \citet{KimJeon2024} for income adjustment on an infinite horizon. Endogenous labor supply and retirement choice appear in \citet{BodieMertonSamuelson1992}, \citet{BodieDetempleOtrubaWalter2004}, \citet{ChoiShim2006}, \citet{FarhiPanageas2007} and \citet{YangKoo2018}, and \citet{JeonPark2020} combine retirement with consumption ratcheting. In those papers the choice concerns when to stop working. Here the retirement date is fixed and the choice concerns how much to invest in income before it. Life cycle models with labor income, such as \citet{Viceira2001} and \citet{CoccoGomesMaenhout2005}, emphasize the role of human capital as an implicit bond holding in portfolio choice, and our effective human wealth makes this role explicit in the presence of a liquidity constraint. The model can also be read as a stylized continuous-time version of the human capital investment problem of \citet{BenPorath1967} with a financial market and a hard retirement date, a feature whose economic rationale is discussed by \citet{Lazear1979}. The concavity of the consumption function and the precautionary role of liquidity constraints relate our findings to \citet{Deaton1991} and \citet{CarrollKimball1996}.

Mathematically, income purchase is a monotone follower or singular control problem in the tradition of \citet{KaratzasShreve1984}. Irreversible investment problems with a similar monotone structure have been analyzed, for example, by \citet{RiedelSu2011} and \citet{DeAngelisFedericoFerrari2017}. Irreversible controls on consumption rather than income appear in the ratcheting and drawdown literature \citep{Dybvig1995,ElieTouzi2008,JeonKooShin2018,DengLiPhamYu2022,JeonKwak2024}. For free boundaries on a finite horizon we rely on the theory of parabolic obstacle problems \citep{Friedman1975,PeskirShiryaev2006}, on the boundary Harnack principle of \citet{Kukuljan2022}, and on the theory of quasilinear oblique derivative problems \citep{Lieberman1996}. The American put is the best known financial example of a free boundary on a finite horizon whose regularity has been analyzed in detail \citep{ChenChadam2007}. There a single free boundary is present from the start, whereas in our problem a second boundary appears from infinity at a positive time. For the integral equation of the American put, uniqueness was proved by \citet{Peskir2005}. Reflected diffusions in time dependent intervals are constructed as in \citet{BurdzyKangRamanan2009}.

\paragraph{Organization.}
\Cref{sec:model} sets up the model and states the main results together with their financial consequences. \Cref{sec:construction} contains the analysis of the free boundary problem: the dual reduction and the no-purchase window, the obstacle problem during the window, the hodograph construction before it, the gluing of the two constructions and the reconstruction of the value function, and the integral equations. \Cref{sec:verification} constructs the optimal strategy, proves optimality and the main results, and studies the dependence on the time to retirement and the long-run limit. \Cref{sec:numerics} presents numerical results and discusses their implications, and \cref{sec:conclusion} concludes.

\section{Model and main results}\label{sec:model}

\subsection{Market, income and controls}
Let $(\Omega,\mathcal F,\mathbb F,\mathbb P)$ be a filtered probability space satisfying the usual conditions and supporting a one-dimensional Brownian motion $B$, and $\mathbb F=(\mathcal F_s)_{s\ge0}$ is the augmented filtration of $B$. The market consists of a bond with constant interest rate $r>0$ and a stock with price $\mathrm dS_s=S_s(\mu\dd s+\sigma\dd B_s)$, $\sigma>0$. We write $\theta=(\mu-r)/\sigma$ for the market price of risk.

The investor works until a mandatory retirement date $T>0$. At time $t\in[0,T)$ she has liquid wealth $x\ge0$ and receives labor income at the rate $\xi>0$. Before $T$ she may raise her income rate irreversibly, and an increase $\mathrm d\Xi_s$ of the income rate costs $p\,\mathrm d\Xi_s$ units of wealth, with $p>0$ fixed. Income ceases at $T$, the accumulated income stock cannot be resold, and there is no pension. A strategy on $[t,T]$ is a triple $(c,\pi,\Xi)$, where the consumption rate $c\ge0$ and the amount $\pi$ invested in the stock are progressively measurable with $\int_t^T(c_s+\pi_s^2)\dd s<\infty$ almost surely, and $\Xi$ is adapted, right-continuous and nondecreasing with $\Xi_{t-}=\xi$. An initial jump $\Xi_t-\xi>0$ is allowed. Wealth evolves according to
\begin{equation}\label{eq:wealth}
 \mathrm dX_s=(rX_s+\sigma\theta\pi_s-c_s+\Xi_s)\dd s+\sigma\pi_s\dd B_s-p\dd\Xi_s,
 \qquad X_{t-}=x,\qquad t\le s\le T.
\end{equation}
The investor cannot borrow against future income: the strategy is admissible if $X_s\ge0$ for all $s\in[t,T]$. In particular every purchase satisfies $p\,\Delta\Xi_s\le X_{s-}$. We denote the set of admissible strategies by $\adm(t,x,\xi)$.

\subsection{Preferences and the retirement payoff}
Preferences are time-additive with discount rate $\beta>0$ and utility $u(c)=c^{1-\gamma}/(1-\gamma)$, $\gamma>0$, $\gamma\ne1$. After retirement the investor has no income and faces the infinite-horizon problem of \citet{Merton1969,Merton1971} with the same preferences. Put
\begin{equation}\label{eq:K}
 K=r+\frac{\beta-r}{\gamma}+\frac{\gamma-1}{2\gamma^2}\theta^2 .
\end{equation}

\begin{assumption}\label{ass:standing}
$\theta\ne0$ and $K>0$.
\end{assumption}

Under \cref{ass:standing} the post-retirement value function and its optimal feedbacks are
\begin{equation}\label{eq:retirement}
 R(x)=\frac{K^{-\gamma}x^{1-\gamma}}{1-\gamma},\qquad c^R=KX,\qquad \pi^R=\frac{\theta}{\gamma\sigma}X,
\end{equation}
with $R(0)=0$ if $\gamma<1$ and $R(0)=-\infty$ if $\gamma>1$. The condition $\theta\ne0$ excludes the degenerate case in which the stock is never held and the dual equation below loses its diffusion term. The pre-retirement problem is
\begin{equation}\label{eq:value}
 V(t,x,\xi)=\sup_{(c,\pi,\Xi)\in\adm(t,x,\xi)}\E\Bigl[\int_t^T\ee^{-\beta(s-t)}u(c_s)\dd s+\ee^{-\beta(T-t)}R(X_T)\Bigr].
\end{equation}
Since $u$ and $R$ have the sign of $1-\gamma$, the expectation is well defined in $[-\infty,\infty]$. The value is finite (\cref{thm:optimal}). For a lower bound it suffices to note that keeping income constant, holding no stock and consuming $\xi/2$ gives a finite objective, because terminal wealth is then at least $\xi\varphi(T-t)/2>0$, where
\begin{equation}\label{eq:annuity}
 \varphi(\tau)=\frac{1-\ee^{-r\tau}}{r}
\end{equation}
is the present value of a unit income flow received for $\tau$ years.

\begin{remark}\label{rem:natural}
The nonnegativity constraint is essential. Suppose instead that wealth were only required to stay above the natural limit $-\Xi_s\varphi(T-s)$. If $\varphi(T-t)>p$, buying $n$ units of income flow and financing them by borrowing raises total wealth $x+\xi\varphi(T-t)$ by $n(\varphi(T-t)-p)$, and the constraint is respected for every $n$. Letting $n\to\infty$ drives the value to $+\infty$ if $\gamma<1$, and to the unattained supremum $0$ if $\gamma>1$. With nonnegative wealth, income can be bought only out of accumulated savings, and the trade-off between liquidity and a higher income becomes the object of the analysis.
\end{remark}

\subsection{Main results}\label{sec:results}

Throughout, $\tau=T-t$ denotes the time to retirement and $z=x/\xi$ the ratio of wealth to income. We use the constants
\begin{equation}\label{eq:constants}
 a=\frac{\theta^2}{2\gamma^2},\qquad b=r-K,\qquad b_0=b-a,\qquad
 \tau_p=\begin{cases}-\dfrac1r\log(1-rp),& rp<1,\\[4pt] +\infty,& rp\ge1,\end{cases}
\end{equation}
so that $\varphi(\tau_p)=p$ when $rp<1$. Scaling the controls shows that
\begin{equation}\label{eq:homogeneity}
 V(t,x,\xi)=\xi^{1-\gamma}v(T-t,x/\xi)
\end{equation}
for a function $v$ of $(\tau,z)\in(0,T]\times[0,\infty)$. In regions where $v_z>0$ and $v_{zz}<0$, the dynamic programming equation for $v$ reads
\begin{equation}\label{eq:primalHJB}
 \max\Bigl\{-v_\tau-\frac{\theta^2}{2}\frac{v_z^2}{v_{zz}}+(rz+1)v_z-\beta v+\frac{\gamma}{1-\gamma}v_z^{(\gamma-1)/\gamma},\ (1-\gamma)v-(z+p)v_z\Bigr\}=0,
\end{equation}
with $v(0,z)=K^{-\gamma}z^{1-\gamma}/(1-\gamma)$. The second term is $\xi^{\gamma}(V_\xi-pV_x)$ and records that one unit of income costs $p$ units of wealth. The boundary $z=0$ carries a state constraint rather than a boundary condition.

\begin{theorem}\label{thm:window}
If $\varphi(T-t)\le p$, then for every $(c,\pi,\Xi)\in\adm(t,x,\xi)$ the strategy $(c,\pi,\xi)$ with constant income is admissible and has terminal wealth at least as large. Consequently, during the last $\tau_p$ years of work, $V$ coincides with the value $V^0$ of the frozen-income problem, in which income cannot be adjusted. If moreover $\varphi(T-t)<p$, a strategy that buys income with positive probability and has a finite objective is strictly worse than $(c,\pi,\xi)$. In particular, no optimal strategy buys income on $(T-\tau_p,T]$. If $rp\ge1$, these statements hold for every $t$.
\end{theorem}

The next theorem shows that the window is exact and describes the structure for longer horizons. We call the set of wealth--income ratios at which the investor buys income immediately the \emph{purchase region}.

\begin{theorem}\label{thm:main}
Let \cref{ass:standing} hold and $rp<1$. There are functions $L\in C^\infty((0,\infty))$, $H\in C^\infty((\tau_p,\infty))$, $\zb\in C^\infty((\tau_p,\infty))$ and $d\in C((0,\infty))\cap C^\infty((\tau_p,\infty))$ with the following properties.
\begin{enumerate}[label=\textup{(\roman*)},leftmargin=*,itemsep=2pt]
\item $0<L<1$. On $(0,\tau_p]$ the function $L$ is nondecreasing and
$L(\tau)\le\min\{K\varphi(\tau),L^0_\infty\}$, where
$L^0_\infty\in(0,1)$ is the free boundary of the stationary
frozen-income obstacle problem \eqref{eq:stationarysingleVI}, given
explicitly in \eqref{eq:Linfty}.
Moreover $L(\tau)/(K\tau)\to1$ as $\tau\downarrow0$.
\item For $\tau>\tau_p$, $L(\tau)<H(\tau)<\infty$ and $0<\zb(\tau)<\infty$. As $\tau\downarrow\tau_p$, $H(\tau)\to\infty$ and $\zb(\tau)\to\infty$, but slowly: there are $h_*>0$ and $C>0$ such that
\begin{equation}\label{eq:blowup}
 \log H(\tau)\le|b_0|\tau+2\sqrt{a\tau\log\frac{C}{\tau-\tau_p}},\qquad \zb(\tau)\le C\,H(\tau),\qquad \tau_p<\tau\le\tau_p+h_* .
\end{equation}
\item The function $v$ in \eqref{eq:homogeneity} belongs to $C^{1,2}((0,T]\times(0,\infty))\cap C((0,T]\times[0,\infty))$, satisfies $v_z>0$ and $v_{zz}<0$, and solves \eqref{eq:primalHJB} classically for $z>0$. The purchase region is $[\zb(\tau),\infty)$ for $\tau>\tau_p$ and is empty for $\tau\le\tau_p$. On the purchase region $v(\tau,z)=d(\tau)^\gamma(z+p)^{1-\gamma}/(1-\gamma)$.
\item At zero wealth, $v_z(\tau,0+)=L(\tau)^{-\gamma}$ and
\begin{equation}\label{eq:curvature}
 v_{zz}(\tau,z)=-\gamma L^{-\gamma}\sqrt{\frac{a}{2(1-L)}}\;z^{-1/2}+O(1)\qquad(z\downarrow0),
\end{equation}
locally uniformly in $\tau>0$, with $L=L(\tau)$. In particular $v(\tau,\cdot)$ is not $C^2$ up to $z=0$.
\end{enumerate}
The boundaries do not depend on $T$, and the horizon only restricts the range $0<\tau\le T$ that is used.
\end{theorem}

When $rp\ge1$, parts (i), (iii) and (iv) hold with $H\equiv\infty$ and an empty purchase region.

The optimal strategy is expressed through the consumption--income ratio $\Lambda_s=c^*_s/\Xi^*_s$, which is kept in the interval $[L(T-s),H(T-s)]$ by reflection. The function $W(\tau,\lambda)\in[0,p]$ that appears in the portfolio rule is introduced in \eqref{eq:XQF}.

\begin{theorem}\label{thm:optimal}
Let \cref{ass:standing} hold. For every $(t,x,\xi)$ the value $V(t,x,\xi)$ is finite and is attained by the following strategy. If $T-t>\tau_p$ and $x/\xi>\zb(T-t)$, the investor first buys income until $X_t/\Xi_t=\zb(T-t)$. Thereafter $\Xi^*$ increases only at times $s$ with $X^*_s/\Xi^*_s=\zb(T-s)$, by the minimal amount that keeps $X^*/\Xi^*\le\zb(T-\cdot)$, and it is constant on $[T-\tau_p,T]$. Consumption and investment are
\begin{equation}\label{eq:policies}
 c^*_s=\Xi^*_s\,\Lambda_s,\qquad
 \pi^*_s=\frac{\theta}{\gamma\sigma}\bigl(X^*_s+W(T-s,\Lambda_s)\,\Xi^*_s\bigr),
\end{equation}
where $\Lambda$ is the process \eqref{eq:logreflection} reflected at $L(T-s)$ and $H(T-s)$. Terminal wealth is $X^*_T=c^*_T/K>0$. The result holds for $0<\gamma<1$ and for $\gamma>1$.
\end{theorem}

The following properties of the optimal policies drive the financial interpretation in \cref{sec:numerics}. We write $\lambda(\tau,z)$ for the optimal consumption--income ratio at the state $(\tau,z)$, so that $c^*=\xi\lambda(\tau,x/\xi)$, and $Q=V_\xi/V_x$ for the marginal rate of substitution between income and wealth.

\begin{proposition}\label{prop:properties}
Let $\tau>0$.
\begin{enumerate}[label=\textup{(\alph*)},leftmargin=*,itemsep=2pt]
\item \emph{Consumption.} $\lambda(\tau,\cdot)$ is continuous, strictly increasing and concave on $[0,\infty)$, $C^1$ on $(0,\infty)$ and strictly concave on $(0,\zb(\tau))$, with $\lambda(\tau,0)=L(\tau)$ ($\zb=\infty$ if $\tau\le\tau_p$). The marginal propensity to consume out of wealth, $\partial c^*/\partial x=\lambda_z$, satisfies
\[
 \lambda_z(\tau,z)=L\sqrt{\frac{a}{2(1-L)z}}\,\bigl(1+O(\sqrt z)\bigr)\quad(z\downarrow0),
 \qquad \lambda_z(\tau,z)=\frac1{d(\tau)}\quad(z\ge\zb(\tau)).
\]
\item \emph{Portfolio.} $\pi^*/\xi=\frac{\theta}{\gamma\sigma}(z+W)$, where $W=0$ at $z=0$, $W=p$ for $z\ge\zb(\tau)$, and $0<W<\min\{p,\varphi(\tau)\}$ for $0<z<\zb(\tau)$. Near zero wealth,
\[
 \frac{\pi^*}{\xi}=\frac{2\sgn\theta}{\sigma}\sqrt{(1-L)z}+O(z).
\]
\item \emph{Shadow price of income.} $W\le Q\le p$ and $Q<\varphi(\tau)$ everywhere, with $Q=p$ exactly for $z\ge\zb(\tau)$.
\item \emph{Retirement.} $c^*_{T-}=KX^*_T=c^R_T$: consumption does not jump at the retirement date.
\end{enumerate}
\end{proposition}

\paragraph{Horizon effects and the long run.}
The remaining results describe how the policies depend on the time to retirement. The first is a monotonicity property of the value function, which constrains how consumption responds.

\begin{theorem}\label{thm:horizon}
Let \cref{ass:standing} hold and $rp<1$.
\begin{enumerate}[label=\textup{(\alph*)},leftmargin=*,itemsep=2pt]
\item For every $z\ge0$, the function $\tau\mapsto v(\tau,z)$ is nondecreasing on $(0,\infty)$.
\item If $\gamma>1$, $d$ is nonincreasing on $(\tau_p,\infty)$ and $d\le1/K$, and if $\gamma<1$, $d$ is nondecreasing and $d\ge1/K$. Hence, at every wealth--income ratio in the purchase region, the consumption--income ratio $\lambda=(z+p)/d(\tau)$ is nondecreasing in the time to retirement when $\gamma>1$ and nonincreasing when $\gamma<1$.
\item If $\gamma>1$, then $\X_\tau\le0$ on $(0,\infty)\times(0,\infty)$, where $\X(\tau,\lambda)$ is the wealth--income ratio at which the consumption--income ratio $\lambda$ is optimal \textup{(}see \eqref{eq:XQF}\textup{)}. Consequently $L$ is nondecreasing on $(0,\infty)$, and $\lambda(\tau,z)$ is nondecreasing in $\tau$ for every $z\ge0$.
\end{enumerate}
\end{theorem}

For $\gamma<1$ the zero-wealth boundary need not be monotone, and \cref{sec:numerics} gives an example in which $L$ first rises and then falls. The dichotomy in (b) is the familiar tension between income and substitution effects: a longer time to retirement raises lifetime resources, which pushes consumption up, and improves the terms on which savings can be converted into income, which pushes it down. The coefficient $\gamma$, the inverse elasticity of intertemporal substitution, decides which effect wins.

The next result compares the problem with the one in which income is frozen. For $\tau>\tau_p$ let $L^0(\tau)$ and $W^0(\tau,\cdot)$ be the zero-wealth boundary and the Euler derivative of the frozen-income problem (\cref{thm:single,lem:singleW}), and let $\Lambda^0_p(\tau)$ be the unique solution of $W^0(\tau,\Lambda^0_p(\tau))=p$, which exists because $W^0(\tau,\cdot)$ increases from $0$ to $\varphi(\tau)>p$.

\begin{theorem}\label{thm:compare}
Let \cref{ass:standing} hold, $rp<1$, and let
\[
 h_1=\sup\Bigl\{h>0:\ \lambda W^0_\lambda\bigl(\tau,\Lambda^0_p(\tau)\bigr)\le\frac{1-rp}{a\gamma}\ \text{for all }\tau\in(\tau_p,\tau_p+h]\Bigr\}.
\]
Then $h_1>0$, and for $\tau\in(\tau_p,\tau_p+h_1]$
\begin{equation}\label{eq:compare}
 L(\tau)\le L^0(\tau),\qquad H(\tau)\le\Lambda^0_p(\tau),\qquad W(\tau,\lambda)\ge W^0(\tau,\lambda)\quad(L(\tau)\le\lambda\le H(\tau)).
\end{equation}
Moreover $\log\Lambda^0_p(\tau)\le|b_0|\tau+2\sqrt{a\tau\log\bigl(\varphi(\tau)\ee^{r\tau}/(\tau-\tau_p)\bigr)}$ for every $\tau>\tau_p$.
\end{theorem}

Thus, at least during the first years in which purchases are possible, the option to buy income makes a constrained investor consume less at zero wealth: savings have acquired an additional use. Whether the inequality $L\le L^0$ persists at longer horizons depends on the parameters. In the numerical examples it holds for $\gamma=2$ and fails for $\gamma=5$ at $\tau=40$, where the income effect of the option dominates.

Finally, far from retirement the problem approaches the perpetual income problem, in which there is no retirement date. Let $u_*$ be the solution of the terminal value problem
\begin{equation}\label{eq:ustar}
 u'(w)=-\frac1a\Bigl(b_0+\frac{1-rw}{u(w)}\Bigr),\quad 0\le w\le p,\qquad u(p)=\frac{1-rp}{a\gamma},
\end{equation}
it exists and is positive on $[0,p]$ (\cref{lem:stationary}). Put $\varsigma=-m_->0$, where $m_-$ is the negative root in \eqref{eq:Linfty}.

\begin{theorem}\label{thm:longrun}
Let \cref{ass:standing} hold and $rp<1$.
\begin{enumerate}[label=\textup{(\alph*)},leftmargin=*,itemsep=2pt]
\item The stationary free boundary problem \eqref{eq:stationaryproblem} has a solution if and only if $u_*(0)<1/a$. The solution is then unique, with
\[
 L^\infty=1-a\,u_*(0),\qquad H^\infty=L^\infty\exp\Bigl(\int_0^p\frac{\mathrm dw}{u_*(w)}\Bigr).
\]
\item If $u_*(0)<1/a$, then, exponentially fast as $\tau\to\infty$, $L(\tau)\to L^\infty$, $H(\tau)\to H^\infty$, $d(\tau)\to d^\infty=\int_0^p\ee^{-\rho_s(w)}\dd w$ and $\zb(\tau)\to\zb^\infty=d^\infty H^\infty-p$, where $\rho_s$ is the stationary solution, and $\X$, $W$ and $Q$ converge locally uniformly.
\item If $\gamma>1$, then $u_*(0)<1/a$.
\item If $u_*(0)\ge1/a$, which can happen only when $\gamma<1$, then $\liminf_{\tau\to\infty}L(\tau)=0$.
\item Whenever $u_*(0)<1/a$, $L^\infty<L^0_\infty$: in the long run the purchase option lowers the consumption--income ratio at zero wealth.
\end{enumerate}
\end{theorem}

In the regime of part (d) the stationary problem has no solution. An investor with a high elasticity of intertemporal substitution and a very long time to retirement then consumes, at zero wealth, a fraction of her income that is not bounded away from zero, because every unit saved can be converted into income on attractive terms. Numerically, both boundaries decay exponentially (\cref{sec:numerics}). Part (e) is the long-run counterpart of \cref{thm:compare}: by \cref{lem:fixedlimit} the frozen-income boundary $L^0(\tau)$ increases to $L^0_\infty$, so for large $\tau$ the purchase option lowers consumption at zero wealth whenever $u_*(0)<1/a$, and along a sequence $\tau\to\infty$ otherwise, although the inequality $L\le L^0$ can fail at intermediate horizons.

\begin{remark}\label{rem:log}
For $\gamma=1$ we have $u=\log$, $K=\beta$, and scaling gives $V(t,x,\xi)=\beta^{-1}\log\xi+v(T-t,x/\xi)$ instead of \eqref{eq:homogeneity}, and in \eqref{eq:primalHJB} the consumption term becomes $-\log v_z-1$ and the gradient term becomes $1/\beta-(z+p)v_z$. The constructions of \cref{sec:single,sec:hodo,sec:reconstruction} involve $\gamma$ only through $a$, $b$, $K$, $c_H=a\gamma/(1-rp)$ and the kernel in \eqref{eq:QFrecon}, and they remain valid for $\gamma=1$. In that case \eqref{eq:dODE} gives $d\equiv1/K$, and the reconstruction gives $Q=\lambda/K-\X$, which is exactly $V_\xi/V_x=1/(\beta v_z)-z$ with $v_z=\lambda^{-1}$. The verification in \cref{sec:verification} goes through with $\widetilde u(y)=-\log y-1$ and the logarithmic Merton value. Only the closed form identity for the optimal objective, which contains the factor $1/(1-\gamma)$, has to be replaced by the expectation itself. The identity $v_z=\lambda^{-1}$ also survives: the budget inequality \eqref{eq:budget}, written with the candidate's shadow prices \eqref{eq:prices} at $(t,x,\xi)$, holds for every other initial wealth $x'$, so the dual bound gives $V(t,x',\xi)\le V(t,x,\xi)+Y_t(x'-x)$ with equality at $x'=x$. Thus $Y_t=(\xi\lambda)^{-1}$ is a supergradient of $V(t,\cdot,\xi)$ that depends continuously on $x$, and $V_x=Y_t$. Consequently \cref{thm:window,thm:optimal}, \cref{prop:properties}, parts (i), (ii) and (iv) of \cref{thm:main}, \cref{thm:horizon}(a),(c), \cref{thm:compare}, and \cref{thm:longrun}(a),(b),(e) hold for $\gamma=1$, and so does part (c): the proof of \cref{lem:stationary}(iv) covers $\gamma\ge1$.
\end{remark}

\paragraph{Discussion of \cref{prop:properties}.}
Part (a) is a finite-horizon version of the concavity of the consumption function studied by \citet{CarrollKimball1996}. Here it is driven by the borrowing constraint rather than by income risk, and the marginal propensity to consume is infinite at zero wealth. Part (b) shows how human capital enters the portfolio. For an investor with constant income who could borrow against it, the risky position would be $\frac\theta{\gamma\sigma}$ times total wealth $X+\Xi\varphi(\tau)$. Under the borrowing constraint the investor acts as if her human wealth were $W\Xi$, and $W$ never exceeds the purchase price $p$. The shadow price $Q$ in (c) lies below the market value $\varphi(\tau)$ of a unit of income, and the difference $\varphi(\tau)-Q$ is an illiquidity discount. Income is bought when this private value reaches the price. Part (d) says that the retirement date, although it removes all income, does not create a consumption drop.

\section{Construction of the free boundaries}\label{sec:construction}

This section constructs the two free boundaries and the dual value function. After the dual reduction and the proof of the no-purchase window (\cref{sec:reduction}), the construction proceeds in three steps, split at the time $\tau_p$ at which the purchase boundary appears. Below $\tau_p$ the problem is a parabolic obstacle problem with a single free boundary, and we solve it and prove the smoothness of its free boundary by established methods: penalization, comparison, and the higher order boundary Harnack inequality (\cref{sec:single}). Above $\tau_p$ both free boundaries are present, and the upper one enters from infinity. We place them on the two sides of a fixed interval by a partial hodograph transformation. The transformed problem is a quasilinear parabolic equation with oblique boundary conditions, to which the comparison principle and the regularity theory for such equations apply, and explicit barriers adapted to the singular corner complete the construction (\cref{sec:hodo}). Finally, we transform back and glue the two constructions along $\tau=\tau_p$ (\cref{cor:bothboundaries}), which gives the regularity of the value function and of both free boundaries on the whole horizon (\cref{sec:reconstruction}). The section ends with integral equations for the free boundaries (\cref{sec:integral}).

\subsection{Dual reduction and the no-purchase window}\label{sec:reduction}

\subsubsection{Linearization}
For $\tau>0$ let $\widehat v(\tau,y)=\sup_{z\ge0}\{v(\tau,z)-yz\}$ be the convex dual of $v(\tau,\cdot)$, and introduce the consumption--income variable $\lambda$ through $y=\lambda^{-\gamma}$. Since $u'(c)=c^{-\gamma}$ and $V_x=\xi^{-\gamma}v_z$, the optimal consumption at a state where $v_z=y$ is $c=\xi\lambda$. We set
\begin{equation}\label{eq:Fdef}
 F(\tau,\lambda)=\lambda^\gamma\,\widehat v(\tau,\lambda^{-\gamma}),
\end{equation}
and define three further functions of $(\tau,\lambda)$:
\begin{equation}\label{eq:XQF}
 \X=\frac{\lambda F_\lambda-\gamma F}{\gamma},\qquad
 Q=F-\lambda F_\lambda,\qquad
 W=\lambda\X_\lambda-\X .
\end{equation}
At an interior maximizer in the definition of $\widehat v$ we have $z=-\widehat v_y$, and a direct computation gives $z=\X(\tau,\lambda)$: the function $\X$ is the inverse of the consumption policy, that is, the wealth--income ratio at which $\lambda$ is optimal. Differentiation of \eqref{eq:XQF} gives
\begin{equation}\label{eq:Qmeaning}
 Q=\frac{V_\xi}{V_x},\qquad F=\frac{Q+\gamma\X}{1-\gamma},\qquad
 \lambda Q_\lambda=-\lambda^2F_{\lambda\lambda}=\gamma(Q-W),\qquad
 W_\lambda=\lambda\X_{\lambda\lambda}.
\end{equation}
In particular the gradient constraint $V_\xi\le pV_x$ becomes $Q\le p$, and the state constraint becomes $\X\ge0$.

Where no constraint is active, the dual function satisfies the linear equation
\[
 \widehat v_\tau=\frac{\theta^2}{2}y^2\widehat v_{yy}+(\beta-r)y\widehat v_y-\beta\widehat v+\frac{\gamma}{1-\gamma}y^{(\gamma-1)/\gamma}+y .
\]
With $D=\lambda\partial_\lambda$ we have $y\partial_y=-D/\gamma$, and conjugating the operator by $\lambda^{-\gamma}$ yields
\begin{equation}\label{eq:P}
 \PP=\partial_\tau-\LL,\qquad \LL=a\lambda^2\partial_{\lambda\lambda}+b\lambda\partial_\lambda-r,
\end{equation}
with $a$ and $b$ as in \eqref{eq:constants}. The Euler operators $D-\gamma$, $1-D$ and $D-1$ commute with $\PP$. Applying them to the equation for $F$ gives
\begin{equation}\label{eq:linearPDEs}
 \PP F=\frac{\gamma}{1-\gamma}\lambda+1,\qquad
 \PP\X=\lambda-1,\qquad
 \PP Q=1,\qquad
 \PP W=1 .
\end{equation}
The retirement payoff gives the initial data
\begin{equation}\label{eq:initial}
 F(0,\lambda)=\frac{\gamma\lambda}{(1-\gamma)K},\qquad
 \X(0,\lambda)=\frac\lambda K,\qquad Q(0,\lambda)=W(0,\lambda)=0 .
\end{equation}
In the logarithmic variable $\zeta=\log\lambda$,
\begin{equation}\label{eq:logoperator}
 \PP=\partial_\tau-a\partial_{\zeta\zeta}-b_0\partial_\zeta+r ,
\end{equation}
a parabolic operator with constant coefficients, uniformly parabolic because $a>0$. Subscripts $\zeta$ denote derivatives in this variable, so that $\X_\zeta=\lambda\X_\lambda$ and $W_\zeta=\lambda W_\lambda=\lambda^2\X_{\lambda\lambda}$.

In the variable $\lambda$ the three regions of the problem are
\[
 \Zr=\{\lambda\le L(\tau)\},\qquad \Nr=\{L(\tau)<\lambda<H(\tau)\},\qquad \Ar=\{\lambda\ge H(\tau)\},
\]
with $H=\infty$ when there is no purchase region. On $\Zr$ wealth is zero and $\X=0$, and on $\Ar$ income is bought and $Q=p$. The two free boundaries $L$ and $H$ are the curves of the consumption--income ratio at which the investor has zero wealth and at which she buys income. In the wealth variable they correspond to $z=0$ and $z=\zb(\tau)=\X(\tau,H(\tau))$.

The argument proceeds by constructing $\X$, $Q$ and $F$ directly, then checking that the reconstructed $v$ solves \eqref{eq:primalHJB} and finally proving optimality by a duality argument that does not assume any regularity of the true value function. The formal derivation above only motivates the construction.

\subsubsection{The final interval}
\begin{proof}[Proof of \cref{thm:window}]
Write $\Xi_s=\xi+I_s$, where $I$ is nondecreasing, right-continuous and includes a possible jump at $t$. Keep $c$ and $\pi$ and replace $\Xi$ by the constant $\xi$. The difference $\chi_s$ between the new and the old wealth solves $\mathrm d\chi_s=(r\chi_s-I_s)\dd s+p\dd I_s$ with $\chi_{t-}=0$. Integrating and using Fubini's theorem together with $\int_u^s\ee^{r(s-v)}\dd v=\ee^{r(s-u)}\varphi(s-u)$,
\begin{equation}\label{eq:surplus}
 \chi_s=\int_{[t,s]}\ee^{r(s-u)}\bigl[p-\varphi(s-u)\bigr]\dd I_u .
\end{equation}
Since $s-u\le T-t$ and $\varphi$ is increasing, every integrand is nonnegative when $\varphi(T-t)\le p$. The new wealth $X_s+\chi_s$ is therefore nonnegative, consumption is unchanged, and the terminal wealth is at least $X_T$. Because $R$ is nondecreasing, the new objective is at least the old one, and $V=V^0$.

If $\varphi(T-t)<p$, the integrand in \eqref{eq:surplus} at $s=T$ is strictly positive, so $\chi_T>0$ on the event $\{I_T>0\}$. Since $R$ is strictly increasing on $[0,\infty)$, $R(X_T+\chi_T)>R(X_T)$ on this event, and the new objective is strictly larger than a finite old one whenever $\mathbb P(I_T>0)>0$. If $rp\ge1$, then $\varphi(\tau)<1/r\le p$ for every $\tau$.
\end{proof}

The comparison is pathwise and respects the borrowing constraint at every date and does not require borrowing against the present value of the purchased income. On $(0,\tau_p]$ the problem is thus a frozen-income problem with a single free boundary, which we treat next.

\subsection{The frozen-income problem}\label{sec:single}

In this subsection income is frozen. We solve the resulting problem for every horizon $\tau>0$, because on $(0,\tau_p]$ it coincides with the original problem and because its solution at $\tau=\tau_p$ is the initial datum of the two-boundary regime. The free boundary of the frozen-income problem is denoted by $L$ in this subsection. When it is compared with the boundary of the original problem for $\tau>\tau_p$, we write $L^0$.

For later comparison, consider also the stationary counterpart of the
frozen-income obstacle problem,
\begin{equation}\label{eq:stationarysingleVI}
    \min\{\X_\infty,\,-\LL\X_\infty-\lambda+1\}=0,
    \qquad \lambda>0,
\end{equation}
in the class of functions for which
$\X_\infty-\lambda/K$ is bounded. We denote by $L^0_\infty$ the
boundary of its contact set. Below we construct its explicit solution
and verify that the contact set is $(0,L^0_\infty]$. Later we show that
$L^0(\tau)\to L^0_\infty$ as $\tau\to\infty$.

Let $m_-<0$ be the negative root of $am^2+b_0m-r=0$ and put
\begin{equation}\label{eq:Linfty}
 L^0_\infty=\frac Kr\,\frac{m_-}{m_--1}=1+\frac{am_-}{r}\in(0,1).
\end{equation}
The two expressions agree because $m_-$ solves the quadratic, and $L^0_\infty<1$ because $m_-<0$.

\begin{theorem}\label{thm:single}
For every $\bar\tau>0$ there is a unique solution of the obstacle problem
\begin{equation}\label{eq:singleVI}
 \min\{\X,\ \PP\X-\lambda+1\}=0\quad\text{in }(0,\bar\tau]\times(0,\infty),\qquad \X(0,\lambda)=\lambda/K,
\end{equation}
in the class of functions for which $\X-\lambda/K$ is bounded. In the variables $(\tau,\zeta)$ it belongs to $W^{1,2}_{\mathfrak q,\mathrm{loc}}$ for every $\mathfrak q<\infty$, and it is smooth on $\{\X>0\}$. Moreover:
\begin{enumerate}[label=\textup{(\alph*)},leftmargin=*,itemsep=2pt]
\item $\max\{0,\lambda/K-\varphi(\tau)\}\le\X(\tau,\lambda)\le\lambda/K$, $-\ee^{-r\tau}\le\X_\tau\le0$ and $\X_{\tau\tau}\ge0$.
\item $\X_\lambda>0$ and $\X_{\lambda\lambda}>0$ on $\{\X>0\}$, and $\{\X(\tau,\cdot)=0\}=(0,L(\tau)]$ for a nondecreasing function $L\in C^\infty((0,\bar\tau])$ with
\begin{equation}\label{eq:singleL}
 0<L(\tau)\le\min\{K\varphi(\tau),L^0_\infty\},\qquad \lim_{\tau\downarrow0}\frac{L(\tau)}{K\tau}=1.
\end{equation}
\item $\X$ is smooth up to the free boundary from the side $\lambda>L$, and there
\begin{equation}\label{eq:lowercurv}
 \X=\X_\lambda=\X_\tau=0,\qquad \X_{\lambda\lambda}(\tau,L+)=\frac{1-L}{aL^2}.
\end{equation}
\end{enumerate}
\end{theorem}

The solution has the stopping representation
\begin{equation}\label{eq:stoprep}
 \X(\tau,\lambda)=\sup_{\varrho\le\tau}\E\Bigl[\int_0^\varrho\ee^{-ru}(\widetilde\Lambda_u-1)\dd u+\ind_{\{\varrho=\tau\}}\ee^{-r\tau}\frac{\widetilde\Lambda_\tau}{K}\Bigr],
 \qquad \mathrm d\widetilde\Lambda_u=b\widetilde\Lambda_u\dd u+\sqrt{2a}\,\widetilde\Lambda_u\dd\widetilde B_u,\ \widetilde\Lambda_0=\lambda,
\end{equation}
where $\widetilde B$ is a Brownian motion on an auxiliary probability space and $\varrho$ ranges over its stopping times. It expresses $\X$ as the largest expected present value of consumption in excess of income that can be financed before wealth is exhausted, in the spirit of \citet{HePages1993} and \citet{ElKarouiJeanblanc1998}. Stopping at $\varrho<\tau$ corresponds to hitting zero wealth, and the payoff $\widetilde\Lambda_\tau/K$ at $\varrho=\tau$ is the wealth needed to continue consumption after retirement.

\subsubsection{Penalization, bounds and monotonicity in time}

\begin{proof}[Proof of existence, uniqueness and part \textup{(a)}]
We work in the variable $\zeta=\log\lambda$, in which $\PP$ has the constant coefficients of \eqref{eq:logoperator}, and we keep the same symbols for functions written in $(\tau,\zeta)$.

Choose a smooth convex nonincreasing function $\psi\ge0$ such that
$\psi=0$ on $[0,\infty)$, $\psi(-1)=1$ and $\psi'$ is bounded. For $\varepsilon>0$ set
$\psi_\varepsilon(s)=\psi((s+\varepsilon)/\varepsilon)$.
Then
\[
    \psi_\varepsilon=0 \quad\text{on }[-\varepsilon,\infty),
    \qquad
    \psi_\varepsilon(-2\varepsilon)=1,
    \qquad
    \psi_\varepsilon'\le0,
    \qquad
    \psi_\varepsilon''\ge0.
\]

For $\mathcal C_n=(0,\bar\tau]\times(-n,n)$ consider
\begin{equation}\label{eq:penalty-n}
\begin{cases}
    \PP\X_{\varepsilon,n}
    =\ee^\zeta-1+\psi_\varepsilon(\X_{\varepsilon,n})
        & \text{in }\mathcal C_n,\\
    \X_{\varepsilon,n}(0,\zeta)=\ee^\zeta/K
        & \text{for }|\zeta|\le n,\\
    \X_{\varepsilon,n}(\tau,\pm n)=\ee^{\pm n}/K
        & \text{for }0\le\tau\le\bar\tau.
\end{cases}
\end{equation}
The $W^{1,2}_{\mathfrak q}$ estimate of \citet[Theorem 7.32]{Lieberman1996}, combined with the Schauder fixed point theorem as in \citet{YangYiDai2006}, gives a unique strong solution
$\X_{\varepsilon,n}\in W^{1,2}_{\mathfrak q}(\mathcal C_n)\cap C(\overline{\mathcal C}_n)$ for
$\mathfrak q>3$. Interior parabolic regularity gives smoothness in $\mathcal C_n$.

Recall that $\varphi(\tau)=(1-\ee^{-r\tau})/r$, so that
$\varphi'+r\varphi=1$. Since
$\PP(\ee^\zeta/K)=\ee^\zeta$, we have
\[
\begin{aligned}
    \PP(-2\varepsilon)
        &=-2r\varepsilon
        \le \ee^\zeta
        =\ee^\zeta-1+\psi_\varepsilon(-2\varepsilon),\\
    \PP\left(
        \frac{\ee^\zeta}{K}-\varphi(\tau)\right)
        &=\ee^\zeta-1
        \le \ee^\zeta-1+
        \psi_\varepsilon\left(
            \frac{\ee^\zeta}{K}-\varphi(\tau)\right),\\
    \PP\left(\frac{\ee^\zeta}{K}\right)
        &=\ee^\zeta
        \ge \ee^\zeta-1+
        \psi_\varepsilon\left(\frac{\ee^\zeta}{K}\right)
        =\ee^\zeta-1.
\end{aligned}
\]
Thus $-2\varepsilon$ and
$\ee^\zeta/K-\varphi(\tau)$ are subsolutions of
\eqref{eq:penalty-n}, while $\ee^\zeta/K$ is a supersolution. These
functions are ordered with the initial and lateral data, so comparison
gives
\begin{equation}\label{eq:penalty-bounds-n}
    \max\left\{
        -2\varepsilon,\frac{\ee^\zeta}{K}-\varphi(\tau)
    \right\}
    \le\X_{\varepsilon,n}\le\frac{\ee^\zeta}{K},
    \qquad
    0\le\psi_\varepsilon(\X_{\varepsilon,n})\le1.
\end{equation}

Fix $\varepsilon>0$ and $m\ge1$. If $n>2m$, then $\mathcal C_m$ is contained in
$(0,\bar\tau]\times(-2m,2m)\subset\mathcal C_n$. By
\eqref{eq:penalty-bounds-n}, the right-hand side of
\eqref{eq:penalty-n} is bounded on
$(0,\bar\tau]\times(-2m,2m)$ by a constant independent of $n$.
The interior parabolic $W^{1,2}_{\mathfrak q}$ estimates \citep[p.~355]{LadyzhenskayaSolonnikovUraltseva1968} therefore give
\[
    \sup_{n>2m}
    \|\X_{\varepsilon,n}\|_{W^{1,2}_{\mathfrak q}(\mathcal C_m)}
    <\infty.
\]
Hence, by weak compactness and compact embedding, one may extract a
subsequence converging weakly in $W^{1,2}_{\mathfrak q}(\mathcal C_m)$ and uniformly on
$\overline{\mathcal C}_m$. Applying this successively for $m=1,2,\ldots$ and
using a diagonal argument, we obtain a single subsequence, still
denoted by $\X_{\varepsilon,n}$, and a function $\X_\varepsilon$
such that
\[
    \X_{\varepsilon,n}\rightharpoonup\X_\varepsilon
    \quad\text{in }W^{1,2}_{\mathfrak q,\mathrm{loc}},
    \qquad
    \X_{\varepsilon,n}\to\X_\varepsilon
    \quad\text{locally uniformly}
\]
on $[0,\bar\tau]\times\R$. In particular,
\[
    \max\left\{
        -2\varepsilon,\frac{\ee^\zeta}{K}-\varphi(\tau)
    \right\}
    \le\X_\varepsilon\le\frac{\ee^\zeta}{K},
    \qquad
    0\le\psi_\varepsilon(\X_\varepsilon)\le1.
\]
Moreover,
\begin{equation}\label{eq:penalty}
    \PP\X_\varepsilon
    =\ee^\zeta-1+\psi_\varepsilon(\X_\varepsilon)
    \qquad\text{in }(0,\bar\tau]\times\R.
\end{equation}
Since $\varphi(\tau)\to0$ as $\tau\downarrow0$, the preceding bounds
give
$\X_\varepsilon(0,\zeta)=\ee^\zeta/K$.

We next let $\varepsilon\downarrow0$. By the preceding bounds, the
right-hand side of \eqref{eq:penalty} is bounded uniformly in
$\varepsilon$ on compact subsets of
$(0,\bar\tau]\times\R$. The local parabolic $W^{1,2}_{\mathfrak q}$ estimates and
compact embedding therefore give, along a subsequence,
$\X_\varepsilon\rightharpoonup\X$ in
$W^{1,2}_{\mathfrak q,\mathrm{loc}}$ and
$\X_\varepsilon\to\X$ locally uniformly.

If $\mathcal U\Subset\{\X>0\}$ is compact, local uniform convergence gives
$\X_\varepsilon>0$ on $\mathcal U$ for all sufficiently small
$\varepsilon$. Hence
$\psi_\varepsilon(\X_\varepsilon)=0$ on $\mathcal U$, and therefore
$\PP\X=\ee^\zeta-1$ in $\{\X>0\}$. Moreover,
\eqref{eq:penalty} gives
$\PP\X_\varepsilon\ge \ee^\zeta-1$, so weak passage to
the limit yields $\PP\X\ge \ee^\zeta-1$ a.e. Consequently,
\[
    \min\{\X,\PP\X-\ee^\zeta+1\}=0
    \qquad\text{in }(0,\bar\tau]\times\R.
\]
The bounds on $\X_\varepsilon$ pass to the limit and give
$\max\{0,\ee^\zeta/K-\varphi(\tau)\}\le\X(\tau,\zeta)\le\ee^\zeta/K$. Since
$\varphi(\tau)\to0$ as $\tau\downarrow0$, they also give
$\X(0,\zeta)=\ee^\zeta/K$.

For uniqueness, let $\X_1$ and $\X_2$ be two solutions in the class
of the theorem. Then $\X_1-\X_2$ is bounded on
$[0,\bar\tau]\times\mathbb R$. Choose $c>0$ so large that
\[
    \Gamma(\tau,\zeta)=\ee^{c\tau}(1+\zeta^2)
\]
satisfies $\PP\Gamma>0$, and let $\eta>0$. On the open set where
$\X_1>\X_2+\eta\Gamma$ we have $\X_1>0$, hence
$\PP\X_1=\ee^\zeta-1\le\PP\X_2<\PP(\X_2+\eta\Gamma)$. If
$M=\|\X_1-\X_2\|_\infty$, then
$\eta\Gamma(\tau,\pm R)\ge M$ for all $0\le\tau\le\bar\tau$ once
$R$ is sufficiently large. Hence
$\X_1\le\X_2+\eta\Gamma$ on the lateral boundary of
$(0,\bar\tau]\times(-R,R)$, while the same inequality holds at
$\tau=0$ because $\X_1$ and $\X_2$ have the same initial value.
The comparison principle on the bounded cylinder gives
$\X_1\le\X_2+\eta\Gamma$. Letting first $R\to\infty$ and then
$\eta\downarrow0$ yields $\X_1\le\X_2$, and interchanging
$\X_1$ and $\X_2$ proves uniqueness.

To derive the estimates in part \textup{(a)}, we use the penalized
solution $\X_\varepsilon$ on the whole line obtained above. By
uniqueness, the whole family converges to $\X$. Since $\PP$ has
constant coefficients, the local estimates on cylinders of unit size
hold with constants that do not depend on the center. Together with
the bounds $-2\varepsilon\le\X_\varepsilon\le\ee^\zeta/K$, they show
that the derivatives of $\X_\varepsilon$ used below grow at most like
$\ee^{|\zeta|}$. The maximum principle is therefore available on the
whole line for these functions: it suffices to repeat the comparison
argument above with $\ee^{c\tau}\cosh(2\zeta)$ in place of $\Gamma$.

Differentiating \eqref{eq:penalty} in $\zeta$ and $\tau$ gives
\begin{equation}\label{eq:penalty-derivatives}
\begin{aligned}
    \left(
        \PP
        -\psi_\varepsilon'(\X_\varepsilon)
    \right)(\X_\varepsilon)_\zeta
    &=\ee^\zeta,\\
    \left(
        \PP
        -\psi_\varepsilon'(\X_\varepsilon)
    \right)(\X_\varepsilon)_\tau
    &=0.
\end{aligned}
\end{equation}
Since $(\X_\varepsilon)_\zeta(0,\zeta)=\ee^\zeta/K>0$ and
$-\psi_\varepsilon'(\X_\varepsilon)\ge0$, the maximum principle gives
$(\X_\varepsilon)_\zeta\ge0$.

To identify the initial value of $(\X_\varepsilon)_\tau$, note that
$\X_\varepsilon(0,\zeta)=\ee^\zeta/K>0$, so the penalty term vanishes at
$\tau=0$. Using
$\PP(\ee^\zeta/K)=\ee^\zeta$ in
\eqref{eq:penalty} gives
$(\X_\varepsilon)_\tau(0,\zeta)=-1$. Moreover,
\[
    \left(
        \PP
        -\psi_\varepsilon'(\X_\varepsilon)
    \right)(-\ee^{-r\tau})
    =
    \psi_\varepsilon'(\X_\varepsilon)\ee^{-r\tau}
    \le0.
\]
Comparison with $-\ee^{-r\tau}$ and $0$ therefore gives
$-\ee^{-r\tau}\le(\X_\varepsilon)_\tau\le0$.

Differentiating once more in $\tau$ gives
\[
    \left(
        \PP
        -\psi_\varepsilon'(\X_\varepsilon)
    \right)(\X_\varepsilon)_{\tau\tau}
    =
    \psi_\varepsilon''(\X_\varepsilon)
    \bigl((\X_\varepsilon)_\tau\bigr)^2
    \ge0.
\]
At $\tau=0$, $\psi_\varepsilon'(\X_\varepsilon)$ and
$\psi_\varepsilon''(\X_\varepsilon)$ vanish because
$\X_\varepsilon(0,\zeta)>0$. Since
$(\X_\varepsilon)_\tau(0,\zeta)=-1$ is independent of $\zeta$,
evaluating the equation for $(\X_\varepsilon)_\tau$ at $\tau=0$
gives
\[
    (\X_\varepsilon)_{\tau\tau}(0,\zeta)
    +r(\X_\varepsilon)_\tau(0,\zeta)=0,
\]
and hence
$(\X_\varepsilon)_{\tau\tau}(0,\zeta)=r$. Therefore, the maximum principle
gives $(\X_\varepsilon)_{\tau\tau}\ge0$. Thus
$-(\X_\varepsilon)_\tau$ is nonincreasing in $\tau$, and
\begin{equation}\label{eq:penalty-time-bound}
    -(\X_\varepsilon)_\tau(\tau,\zeta)
    \le
    \frac{\X_\varepsilon(0,\zeta)
          -\X_\varepsilon(\tau,\zeta)}{\tau}
    \le
    \frac{\ee^\zeta/K+2\varepsilon}{\tau}.
\end{equation}

Passing to the limit $\varepsilon\downarrow0$ gives
$\X_\zeta\ge0$, $-\ee^{-r\tau}\le\X_\tau\le0$,
$\X_{\tau\tau}\ge0$, and
$-\X_\tau(\tau,\zeta)\le \ee^\zeta/(K\tau)$.
Returning to $\lambda=\ee^\zeta$, for which
$\X_\zeta=\lambda\X_\lambda$, gives part \textup{(a)}.
\end{proof}

\subsubsection{Convexity in \texorpdfstring{$\lambda$}{lambda} and the free boundary}
The representation \eqref{eq:stoprep} follows from the obstacle problem by the usual verification argument with the generalized It\^o formula for $W^{1,2}_{\mathfrak q,\mathrm{loc}}$ functions. For a fixed stopping rule, the expression in \eqref{eq:stoprep} is affine in the initial value $\lambda$ with nonnegative slope, because $\widetilde\Lambda$ is proportional to $\lambda$. The supremum is therefore nondecreasing and convex in $\lambda$, and the contact set $\{\X(\tau,\cdot)=0\}$ is an interval of the form $(0,L(\tau)]$ (possibly empty at this stage). As $\X_\tau\le0$, contact sets grow in $\tau$ and $L$ is nondecreasing. On $\{\X>0\}$ the function $\lambda\X_\lambda\ge0$ solves $\PP(\lambda\X_\lambda)=\lambda>0$, and the strong maximum principle gives $\X_\lambda>0$.

\begin{proof}[Proof of the bounds in \eqref{eq:singleL}]
The stationary problem introduced above can be solved explicitly. In
the continuation region its equation is
$-\LL\X_\infty=\lambda-1$. Solving this equation under the growth
condition that $\X_\infty-\lambda/K$ is bounded, and imposing value
matching and smooth fit at the free boundary $L^0_\infty$, gives
\begin{equation}\label{eq:stationarysingle}
 \X_\infty(\lambda)=
 \begin{cases}
 0,&\lambda\le L^0_\infty,\\[2pt]
 \dfrac{\lambda}{K}-\dfrac1r
 -\dfrac{(L^0_\infty)^{1-m_-}}{Km_-}\lambda^{m_-},
 &\lambda>L^0_\infty .
 \end{cases}
\end{equation}
To verify that \eqref{eq:stationarysingle} solves
\eqref{eq:stationarysingleVI}, note first that it is $C^1$. On
$\lambda>L^0_\infty$, it is positive and satisfies
$-\LL\X_\infty=\lambda-1$, while on
$\lambda\le L^0_\infty$ it vanishes and
$1-\lambda\ge0$ since $L^0_\infty<1$ (see \eqref{eq:Linfty}). Hence
\eqref{eq:stationarysingle} solves \eqref{eq:stationarysingleVI},
with contact set exactly $(0,L^0_\infty]$. Moreover,
$0\le\X_\infty\le\lambda/K$, so that
$\X_\infty-\lambda/K$ is bounded as required.

Viewed as a time independent function, $\X_\infty$ satisfies the same
variational inequality as \eqref{eq:singleVI} and lies below the
prescribed initial datum,
$\X_\infty(\lambda)\le\lambda/K=\X(0,\lambda)$. The comparison
principle therefore gives
$\X(\tau,\lambda)\ge\X_\infty(\lambda)$, and hence
$L(\tau)\le L^0_\infty$. On the other hand, the lower bound
$\X(\tau,\lambda)\ge\lambda/K-\varphi(\tau)$, evaluated at
$\lambda=L(\tau)$ where $\X(\tau,L(\tau))=0$, yields
$L(\tau)\le K\varphi(\tau)$.

To obtain a lower bound on $L$, fix an integer $n\ge2$ and set
\[
    \ell_n(\tau)
    =\frac{K(n-1)}{nr}(1-\ee^{-r\tau}).
\]
Define $\Upsilon_n=0$ for $\lambda\le\ell_n$ and
\begin{equation}\label{eq:Sn}
 \Upsilon_n(\tau,\lambda)
 =\frac1K\Bigl[
   \lambda-\frac{n}{n-1}\ell_n
   +\frac{\ell_n^{\,n}}{(n-1)\lambda^{n-1}}
 \Bigr],
 \qquad \lambda>\ell_n .
\end{equation}
The function $\Upsilon_n$ is $C^1$, nonnegative and convex in
$\lambda$, and satisfies
$\Upsilon_n(0,\lambda)=\lambda/K$. Since
$\ell_n'=K(n-1)/n-r\ell_n$, a direct computation
(\cref{app:computations}) gives, for $\lambda>\ell_n$,
\[
 \PP\Upsilon_n-\lambda+1
 =
 \Bigl[1-\frac{(an+K)\ell_n}{K}\Bigr]
 \Bigl(\frac{\ell_n}{\lambda}\Bigr)^{n-1}.
\]
This is nonnegative whenever
$\ell_n\le K/(an+K)$. On $\lambda\le\ell_n$, the residual is
$1-\lambda\ge0$ for sufficiently small $\tau$. Thus
$\Upsilon_n$ is a supersolution of \eqref{eq:singleVI} for
sufficiently small $\tau$, and comparison gives
$\X\le\Upsilon_n$. Since $\Upsilon_n=0$ on
$\lambda\le\ell_n$, it follows that
\[
    L(\tau)\ge\ell_n(\tau)
\]
for sufficiently small $\tau$. As $L$ is nondecreasing, this also
shows that $L(\tau)>0$ for every $\tau>0$.

Finally, combining the preceding lower bound with
$L(\tau)\le K\varphi(\tau)$ gives
\[
 \frac{n-1}{n}
 \le
 \liminf_{\tau\downarrow0}\frac{L(\tau)}{K\tau}
 \le
 \limsup_{\tau\downarrow0}\frac{L(\tau)}{K\tau}
 \le1.
\]
Letting $n\to\infty$ yields
$L(\tau)/(K\tau)\to1$ as $\tau\downarrow0$.
\end{proof}

\subsubsection{Lipschitz continuity of the free boundary}

The next estimate controls the time derivative by the spatial
derivative in the positive region. It is obtained from the penalized
equations and does not require any a priori regularity of the free
boundary.

\begin{lemma}\label{lem:cone}
For every $\delta>0$ and $2\delta\le\tau\le\bar\tau$, on
$\{\X>0\}$,
\begin{equation}\label{eq:cone}
    0\le-\X_\tau
    \le\frac{2}{K\delta^2}\X_\zeta .
\end{equation}
Consequently $\log L$ is Lipschitz on $[2\delta,\bar\tau]$.
\end{lemma}

\begin{proof}
Set $C_\delta=2/(K\delta)$ and
$E_\varepsilon=4\varepsilon/r$. For $\tau\ge\delta$, define
\[
    \Theta_\varepsilon
    :=
    C_\delta(\tau-\delta)(\X_\varepsilon)_\zeta
    +(\tau-\delta)^2(\X_\varepsilon)_\tau
    +E_\varepsilon .
\]
By \eqref{eq:penalty-derivatives},
\[
\begin{aligned}
    &\left(
        \PP
        -\psi_\varepsilon'(\X_\varepsilon)
    \right)\Theta_\varepsilon\\
    &\quad=
    C_\delta(\X_\varepsilon)_\zeta
    +C_\delta(\tau-\delta)\ee^\zeta
    +2(\tau-\delta)(\X_\varepsilon)_\tau
    +(r-\psi_\varepsilon'(\X_\varepsilon))E_\varepsilon .
\end{aligned}
\]
For $\tau\ge\delta$, \eqref{eq:penalty-time-bound} gives
\[
\begin{aligned}
    -2(\tau-\delta)(\X_\varepsilon)_\tau
    &\le
    \frac{2(\tau-\delta)}{\tau}
    \left(\frac{\ee^\zeta}{K}+2\varepsilon\right)\\
    &\le
    \frac{2(\tau-\delta)}{K\delta}\ee^\zeta
    +4\varepsilon\\
    &=
    C_\delta(\tau-\delta)\ee^\zeta+rE_\varepsilon .
\end{aligned}
\]
Since $\psi_\varepsilon'\le0$, we therefore obtain
\[
    \left(
        \PP
        -\psi_\varepsilon'(\X_\varepsilon)
    \right)\Theta_\varepsilon
    \ge
    C_\delta(\X_\varepsilon)_\zeta
    \ge0.
\]
Moreover,
$\Theta_\varepsilon(\delta,\zeta)=E_\varepsilon\ge0$.

We first establish the weighted estimate on
$[\delta,\bar\tau]\times\mathbb R$.
Since $(\X_\varepsilon)_\zeta\ge0$ and
$-1\le(\X_\varepsilon)_\tau\le0$, we have
$\Theta_\varepsilon\ge-(\bar\tau-\delta)^2$.
Choose $A>0$ so large that
$\Gamma_A(\tau,\zeta)=\ee^{A(\tau-\delta)}(1+\zeta^2)$ satisfies
\[
    \left(
        \PP
        -\psi_\varepsilon'(\X_\varepsilon)
    \right)\Gamma_A\ge0
    \quad\text{on }[\delta,\bar\tau]\times\mathbb R.
\]
For each $\eta>0$,
$\Theta_\varepsilon+\eta\Gamma_A$ is nonnegative on $\tau=\delta$.
Since $\Theta_\varepsilon$ is bounded from below and
$\Gamma_A(\tau,\zeta)\to\infty$ as $|\zeta|\to\infty$ uniformly for
$\tau\in[\delta,\bar\tau]$, the same function is nonnegative on
$\zeta=\pm R$ for all sufficiently large $R$.
The maximum principle on
$[\delta,\bar\tau]\times[-R,R]$ gives
$\Theta_\varepsilon+\eta\Gamma_A\ge0$.
Letting $R\to\infty$ and then $\eta\downarrow0$ yields
$\Theta_\varepsilon\ge0$ on
$[\delta,\bar\tau]\times\mathbb R$.

We next pass to the limit.
Fix an open set $U\Subset\{\X>0\}\cap((\delta,\bar\tau]\times\mathbb R)$.
The local uniform convergence $\X_\varepsilon\to\X$ and the
positivity of $\X$ on $\overline U$ imply that the penalty vanishes
on $U$ for all sufficiently small $\varepsilon$.
The interior $W^{1,2}_p$ estimate then gives, after passing to a
subsequence,
$\X_\varepsilon\rightharpoonup\X$ weakly in
$W^{1,2}_{\mathfrak q}(U)$.
Passing to the limit in $\Theta_\varepsilon\ge0$ in the weak sense
gives
\[
    (\tau-\delta)^2(-\X_\tau)
    \le C_\delta(\tau-\delta)\X_\zeta
    \quad\text{a.e. in }U.
\]
Since $\X$ satisfies the linear equation in $\{\X>0\}$, interior
parabolic regularity implies that its first derivatives are
continuous. Hence the preceding inequality holds pointwise in
$\{\X>0\}$ for $\tau>\delta$.

Finally, for $\tau\ge2\delta$, we have
$\tau-\delta\ge\delta$. Recalling
$C_\delta=2/(K\delta)$, we obtain
\[
    0\le-\X_\tau
    \le\frac{C_\delta}{\tau-\delta}\X_\zeta
    \le\frac{2}{K\delta^2}\X_\zeta .
\]
This proves \eqref{eq:cone}.

We next turn to the free boundary. Since
$\X_\zeta\ge0$ and, in $\{\X>0\}$,
\[
    \PP\X_\zeta=\ee^\zeta>0,
\]
the strong maximum principle gives
$\X_\zeta>0$ in $\{\X>0\}$. Thus, for every $\kappa>0$, the implicit function theorem gives a
unique smooth function $\zeta^\kappa$ satisfying
\[
    \X(\tau,\zeta^\kappa(\tau))=\kappa .
\]
Differentiating this identity gives
\[
    \frac{\mathrm d}{\mathrm d\tau}
        \zeta^\kappa(\tau)
    =
    -\frac{\X_\tau}{\X_\zeta}
      \bigg|_{\zeta=\zeta^\kappa(\tau)}.
\]
By \eqref{eq:cone},
\[
    0\le
    \frac{\mathrm d}{\mathrm d\tau}
        \zeta^\kappa(\tau)
    \le
    \frac{2}{K\delta^2},
    \qquad
    2\delta\le\tau\le\bar\tau.
\]
As $\kappa\downarrow0$,
$\zeta^\kappa(\tau)\downarrow\log L(\tau)$ for every $\tau$.
The uniform Lipschitz bound therefore passes to the limit, and
$\log L$ is Lipschitz on $[2\delta,\bar\tau]$.
\end{proof}

\subsubsection{Smoothness of the free boundary}
\begin{proof}[Proof of the regularity in \textup{(b)} and \textup{(c)}]
Write $\zeta_L=\log L$ and
$\mathcal O=\{(\tau,\zeta):\zeta>\zeta_L(\tau)\}$.

Fix a compact interval
$I=[\tau_1,\tau_2]\Subset(0,\bar\tau]$ and choose $\delta>0$ with
$2\delta<\tau_1$. All arguments below are local in a neighborhood of
$\{(\tau,\zeta_L(\tau)):\tau\in I\}$. Since $I$ is arbitrary,
the resulting regularity holds on $(0,\bar\tau]$.

Since $\X\in W^{1,2}_{\mathfrak q,\mathrm{loc}}$ for all $\mathfrak q<\infty$,
the parabolic Sobolev embedding, applied with $\mathfrak q>3$, implies that
$\X$ and $\X_\zeta$ are continuous. Both vanish on the contact set,
hence on $\partial\mathcal O$.
By \cref{lem:cone}, $\zeta_L$ is Lipschitz
on a neighborhood of $I$.
Its oscillation over a time interval of length $\delta_1$ is
$O(\delta_1)=o(\sqrt{\delta_1})$, and $\mathcal O$ is a parabolic
$C^1$ domain.

The function $\Phi=\X_\zeta$ is positive in $\mathcal O$, vanishes
on $\partial\mathcal O$ and solves
$\PP\Phi=\ee^\zeta$.
The function $\Psi=-\X_\tau\ge0$ solves $\PP\Psi=0$ in
$\mathcal O$.

Moreover, \eqref{eq:cone} gives
$0\le\Psi\le C_I\Phi$ with $C_I=2/(K\delta^2)$ near the portion of the
free boundary corresponding to $I$. Since $\Phi$ is continuous and vanishes on
$\partial\mathcal O$, it follows that $\Psi$ extends continuously to
this portion of the boundary with zero boundary value.

After multiplication by $\ee^{r\tau}$ both equations are in
non-divergence form without zero order term, with constant principal
coefficient $a$ and constant drift.

More precisely, with
$\widehat\Phi=\ee^{r\tau}\Phi$ and
$\widehat\Psi=\ee^{r\tau}\Psi$,
\[
    (\partial_\tau-a\partial_{\zeta\zeta}
      -b_0\partial_\zeta)\widehat\Phi
    =\ee^{r\tau+\zeta},
    \qquad
    (\partial_\tau-a\partial_{\zeta\zeta}
      -b_0\partial_\zeta)\widehat\Psi=0,
\]
and $\widehat\Psi/\widehat\Phi=\Psi/\Phi$.
On each fixed local boundary neighborhood the coefficients and
right-hand sides satisfy the regularity and boundedness assumptions
of the boundary Harnack estimates.

It remains to verify the linear nondegeneracy of $\Phi$ near the
free boundary. By the preceding paragraph, $\X_\tau$ tends to zero as
the free boundary is approached from $\mathcal O$, and so do $\X$ and
$\X_\zeta$. Using
\[
    \X_\tau-a\X_{\zeta\zeta}
    -b_0\X_\zeta+r\X=\ee^\zeta-1
    \quad\text{in }\mathcal O,
\]
we obtain
\[
    \X_{\zeta\zeta}
    \bigl(\tau,\zeta_L(\tau)+\bigr)
    =\frac{1-L(\tau)}{a}>0,
    \qquad \tau\in I.
\]
Moreover, the equation implies that $\X_{\zeta\zeta}$ extends
continuously to the free boundary from inside $\mathcal O$.
Since $I$ is compact and $L<1$, there exist constants
$c_I,\rho_I>0$ such that
$\X_{\zeta\zeta}(\tau,\zeta)\ge c_I$ whenever
$\tau\in I$ and
$0<\zeta-\zeta_L(\tau)<\rho_I$.
Therefore, using
$\Phi(\tau,\zeta_L(\tau))=\X_\zeta(\tau,\zeta_L(\tau))=0$,
\[
    \Phi(\tau,\zeta)
    =
    \int_{\zeta_L(\tau)}^\zeta
    \X_{\zeta\zeta}(\tau,y)\dd y
    \ge
    c_I\bigl(\zeta-\zeta_L(\tau)\bigr).
\]
Thus $\Phi$ satisfies the linear nondegeneracy required for the
boundary Harnack estimates.

We write $C^\beta_{\mathrm{par}}$ for the parabolic H\"older spaces,
in which one time derivative counts as two spatial derivatives. The boundary Harnack
inequalities of \citet[Theorems 1.2 and 1.3]{Kukuljan2022} apply to
$\Psi/\Phi$: this quotient belongs to $C^{\alpha}_{\mathrm{par}}$ up to
$\partial\mathcal O$ for every $\alpha<1$, and if
$\partial\mathcal O$ is a parabolic $C^\beta$ domain,
$\beta>1$ not an integer, the quotient belongs to $C^\beta_{\mathrm{par}}$
up to the boundary.

We use the following convention. In one space dimension the domain $\mathcal O$ is a parabolic
$C^\beta$ domain, $\beta>1$ not an integer, when
$\zeta_L\in C^{\beta/2}$ as a function of time, because parabolic
distance scales like the square root of time.
If $f$ belongs to the parabolic H\"older space
$C^\beta_{\mathrm{par}}(\overline{\mathcal O})$ and
$\zeta_L\in C^{\beta/2}$, then the trace
$\tau\mapsto f(\tau,\zeta_L(\tau))$ belongs to $C^{\beta/2}$:
for $\beta<2$ this follows directly from the parabolic H\"older
estimate, and for larger $\beta$ it follows from the chain rule,
since $\partial_\tau f\in C^{\beta-2}_{\mathrm{par}}$ and
$\partial_\zeta f\in C^{\beta-1}_{\mathrm{par}}$ by the equivalence of the usual
characterizations of parabolic H\"older spaces
\citep[Chapter IV]{Lieberman1996}.

To identify the boundary trace of $\Psi/\Phi$, return to the level
curves $\zeta^\kappa$ introduced in the proof of \cref{lem:cone}.
They satisfy
\[
    (\zeta^\kappa)'(\tau)
    =
    \frac{\Psi}{\Phi}
    \bigl(\tau,\zeta^\kappa(\tau)\bigr).
\]

For each $\tau\in I$,
$\zeta^\kappa(\tau)\downarrow\zeta_L(\tau)$ as
$\kappa\downarrow0$.
Both $\zeta^\kappa$ and $\zeta_L$ are continuous on $I$, and the
convergence is monotone in $\kappa$. By Dini's theorem it is uniform
on $I$.
Since $\Psi/\Phi$ extends continuously to the boundary, it is
uniformly continuous on a smaller compact boundary neighborhood.
Consequently,
\[
    \frac{\Psi}{\Phi}
    \bigl(\tau,\zeta^\kappa(\tau)\bigr)
    \longrightarrow
    \frac{\Psi}{\Phi}
    \bigl(\tau,\zeta_L(\tau)\bigr)
\]
uniformly for $\tau\in I$.
For $\tau_a,\tau_b\in I$, integration of the equation of the level
curves gives
\[
    \zeta^\kappa(\tau_b)-\zeta^\kappa(\tau_a)
    =
    \int_{\tau_a}^{\tau_b}
    \frac{\Psi}{\Phi}
    \bigl(s,\zeta^\kappa(s)\bigr)\dd s .
\]
Passing to the limit $\kappa\downarrow0$ yields
\[
    \zeta_L(\tau_b)-\zeta_L(\tau_a)
    =
    \int_{\tau_a}^{\tau_b}
    \frac{\Psi}{\Phi}
    \bigl(s,\zeta_L(s)\bigr)\dd s .
\]
The boundary trace on the right-hand side is continuous, and hence
\begin{equation}\label{eq:traceidentity}
    \zeta_L\in C^1,\qquad
    \zeta_L'(\tau)
    =
    \frac{\Psi}{\Phi}
    \bigl(\tau,\zeta_L(\tau)\bigr).
\end{equation}

The regularity of $\zeta_L$ is now improved in steps.
Since $\zeta_L$ is Lipschitz,
$\zeta_L\in C^{\beta/2}$ for every $\beta\in(1,2)$, so
$\mathcal O$ is a parabolic $C^\beta$ domain for such $\beta$.

The coefficients of the equations for
$\widehat\Phi$ and $\widehat\Psi$ are constant and their right-hand
sides are smooth, so all higher order hypotheses of the boundary
Harnack estimate are satisfied after localizing and rescaling around
each boundary point. The higher order boundary Harnack inequality gives
$\Psi/\Phi\in C^\beta_{\mathrm{par}}$ up to the boundary, and by the trace
property and \eqref{eq:traceidentity},
$\zeta_L'\in C^{\beta/2}$.
Thus $\zeta_L\in C^{1+\beta/2}$, and in particular
$\zeta_L\in C^{1+\alpha}$ for every
$\alpha\in(0,\frac12)$.

Suppose now that
$\zeta_L\in C^{k+\alpha}$ for an integer $k\ge1$ and some
$\alpha\in(0,\frac12)$.
Then $\mathcal O$ is a parabolic $C^\beta$ domain with
$\beta=2k+2\alpha$, which is not an integer.
The higher order boundary Harnack inequality gives
$\Psi/\Phi\in C^\beta_{\mathrm{par}}$, hence
$\zeta_L'\in C^{\beta/2}=C^{k+\alpha}$ and
$\zeta_L\in C^{k+1+\alpha}$.
By induction $\zeta_L\in C^\infty$.
Since $I\Subset(0,\bar\tau]$ was arbitrary,
$\zeta_L$, and hence $L=\ee^{\zeta_L}$, is smooth on
$(0,\bar\tau]$.

Boundary Schauder estimates for
$\PP\X=\lambda-1$ with zero Dirichlet data on the smooth boundary
give smoothness of $\X$ up to $\partial\mathcal O$ from inside,
locally away from $\tau=0$. This justifies evaluating the equation
up to the free boundary and differentiating the boundary identity
classically. Differentiating $\X(\tau,L(\tau))=0$ gives
\[
    0
    =
    \X_\tau(\tau,L(\tau))
    +L'(\tau)\X_\lambda(\tau,L(\tau)).
\]
Using $\X_\lambda(\tau,L)=0$ yields
$\X_\tau(\tau,L)=0$.
Evaluating $\PP\X=\lambda-1$ at $\lambda=L+$ then gives
$-aL^2\X_{\lambda\lambda}(\tau,L+)=L-1$, which is
\eqref{eq:lowercurv}.

Finally,
$W=\lambda\X_\lambda-\X$ solves $\PP W=1$ on $\mathcal O$, so
$W_\zeta=\lambda^2\X_{\lambda\lambda}$ solves
$\PP W_\zeta=0$.
It is nonnegative by convexity and, by \eqref{eq:lowercurv},
has boundary value $(1-L)/a>0$.
The strong maximum principle gives
$W_\zeta>0$ in $\mathcal O$, and hence
$\X_{\lambda\lambda}>0$ on $\{\X>0\}$.
\end{proof}
\subsubsection{The Euler derivative and its tail}

The two-boundary regime will be constructed by using
$W=\lambda\X_\lambda-\X$ as the spatial variable in a partial
hodograph transformation. At $\tau=\tau_p$, where
$\varphi(\tau_p)=p$, the initial hodograph profile is obtained by
inverting $w=W(\tau_p,\ee^\zeta)$ on $0\le w<p$.
We therefore need the strict monotonicity and range of $W$, together
with sufficient control of its tail to obtain the required bounds on
the derivative of the inverse profile. The following lemma collects
these properties.

\begin{lemma}\label{lem:singleW}
In the frozen-income problem let $W=\lambda\X_\lambda-\X$ for
$\lambda>L$ and $W=0$ for $\lambda\le L$.
Then, for every $\tau>0$ and finite $\lambda>L(\tau)$,
\begin{equation}\label{eq:Wsinglebounds}
 0<W(\tau,\lambda)<\varphi(\tau),\qquad
 W_\lambda>0,\qquad
 \lim_{\lambda\to\infty}W(\tau,\lambda)=\varphi(\tau),
\end{equation}
and $\varphi(\tau)-W(\tau,\ee^\zeta)$ and its $\zeta$-derivatives
decay faster than any exponential as $\zeta\to\infty$, locally
uniformly in $\tau$. For fixed $\tau_0>0$ the inverse
$\zeta=\rho_0(w)$ of $w=W(\tau_0,\ee^\zeta)$,
$0\le w<\varphi(\tau_0)$, is smooth up to $w=0$, tends to
$+\infty$ as $w\uparrow\varphi(\tau_0)$, and satisfies
\begin{equation}\label{eq:initialhodobounds}
 \rho_0'(0)=\frac{a}{1-\ee^{\rho_0(0)}},\qquad
 0<m_0\le\rho_0'(w)\le
 \frac{C_0}{\varphi(\tau_0)-w}
\end{equation}
for constants $m_0,C_0>0$.
\end{lemma}

\begin{proof}
Convexity and $\X\le\lambda/K$ imply
$\X_\lambda\le1/K$, and with
$\X\ge\lambda/K-\varphi$ this gives $W\le\varphi$.
Since
$W_\lambda=\lambda\X_{\lambda\lambda}>0$ and
$W(\tau,L)=0$, we have $W>0$ for $\lambda>L$.
The function $U=\varphi-W$ solves $\PP U=0$ in the positive region
(note $\PP\varphi=1$), vanishes at $\tau=0$ and equals
$\varphi(\tau)>0$ at $\lambda=L$, so the strong maximum principle gives
$U>0$, and hence $W<\varphi$.

Fix $\zeta_0>\log L^0_\infty\ge\log L(\tau)$ and let
$g_0(\tau)=U(\tau,\ee^{\zeta_0})$.
On the half-line $\zeta>\zeta_0$, with
$\eta=\zeta-\zeta_0$, the first passage representation for Brownian
motion with drift $b_0$ and variance $2a$ per unit time gives
\begin{equation}\label{eq:tailkernel}
 U(\tau,\ee^\zeta)
 =
 \int_0^\tau
 g_0(\tau-v)\,\ee^{-rv}\,k(v,\eta)\dd v,\qquad
 k(v,\eta)
 =
 \frac{\eta}{\sqrt{4\pi av^3}}
 \exp\Bigl[-\frac{(\eta+b_0v)^2}{4av}\Bigr].
\end{equation}
The Gaussian factor in $k$ gives the superexponential decay of
$U$ and its $\zeta$ derivatives, locally uniformly in $\tau$.
In particular, $W(\tau,\lambda)\to\varphi(\tau)$ as
$\lambda\to\infty$. Since $W_\lambda>0$, the inverse
$\rho_0$ is therefore well defined on
$[0,\varphi(\tau_0))$ and
$\rho_0(w)\to+\infty$ as
$w\uparrow\varphi(\tau_0)$.

Moreover,
\[
 -\frac{\partial_\eta k}{k}
 =
 \frac{\eta}{2av}+\frac{b_0}{2a}-\frac1\eta
 \ge
 \frac{\eta}{2a\tau_0}+\frac{b_0}{2a}-\frac1\eta
 \qquad(v\le\tau_0),
\]
which is bounded below by a positive constant $c$ for large $\eta$.
Since $g_0\ge0$, \eqref{eq:tailkernel} then gives
$-U_\zeta\ge cU$ at $\tau=\tau_0$ for large $\zeta$.
Since $\rho_0'=1/W_\zeta$ and
$(\varphi(\tau_0)-w)\rho_0'=U/(-U_\zeta)$,
the upper bound in \eqref{eq:initialhodobounds} holds near
$w=\varphi(\tau_0)$. On each remaining compact subinterval,
$W_\zeta$ is continuous and strictly positive, so
$\rho_0'=1/W_\zeta$ is bounded above.
The lower bound follows because $W_\zeta$ is continuous and bounded
on $[\log L(\tau_0),\infty)$, with
$W_\zeta\to0$ as $\zeta\to\infty$.

Finally, \eqref{eq:lowercurv} gives
$W_\zeta(\tau,L+)=(1-L)/a>0$. Thus the inverse is smooth up to
$w=0$, and
\[
 \rho_0'(0)
 =
 \frac{1}{W_\zeta(\tau_0,L+)}
 =
 \frac{a}{1-L(\tau_0)}
 =
 \frac{a}{1-\ee^{\rho_0(0)}}.
\]
\end{proof}

For later comparison with the long-run behavior of the two-boundary problem, we also record
the infinite-horizon limit of the frozen-income solution.

\begin{lemma}\label{lem:fixedlimit}
As $\tau\to\infty$, the frozen-income boundary increases to $L^0_\infty$ and $\X(\tau,\cdot)$ decreases to $\X_\infty$, locally uniformly together with the first derivative.
\end{lemma}

\begin{proof}
By \cref{thm:single}, $\X(\tau,\lambda)$ is nonincreasing in $\tau$ and bounded below by $\X_\infty$, and $L$ is nondecreasing and bounded by $L^0_\infty$. Let $\widehat\X=\lim_{\tau\to\infty}\X(\tau,\cdot)$ and $\widehat L=\lim_{\tau\to\infty}L(\tau)\in(0,L^0_\infty]$. Since $\PP\X=(\lambda-1)\ind_{\{\X>0\}}$ is bounded on compact sets uniformly in $\tau$, the parabolic $L^{\mathfrak q}$ estimates give bounds in $W^{1,2}_{\mathfrak q}$ on the windows $[\tau,\tau+1]\times[\lambda_1,\lambda_2]$ that do not depend on $\tau$. Hence the convergence holds in $C^1_{\mathrm{loc}}$, and $\widehat\X\in C^1$ with $\widehat\X=0$ on $(0,\widehat L]$. For $\lambda>\widehat L$ we have $\lambda>L(\tau)$ for all $\tau$, so $\X$ solves $\PP\X=\lambda-1$ near $\lambda$ for all large $\tau$, and the limit of the solutions shifted in time shows $-\LL\widehat\X=\lambda-1$ on $(\widehat L,\infty)$. Moreover $\widehat\X-\lambda/K$ is bounded by \cref{thm:single}(a). The solutions of $-\LL f=\lambda-1$ are $\lambda/K-1/r+c_+\lambda^{m_+}+c_-\lambda^{m_-}$, where $m_+>0>m_-$ are the roots of $am^2+b_0m-r=0$. Boundedness of $\widehat\X-\lambda/K$ forces $c_+=0$, and the conditions $\widehat\X(\widehat L)=\widehat\X'(\widehat L)=0$ then give $\widehat L=(K/r)\,m_-/(m_--1)=L^0_\infty$ and $\widehat\X=\X_\infty$.
\end{proof}

\subsection{The two-boundary regime}\label{sec:hodo}

Assume $rp<1$ and consider $\tau>\tau_p$.
At $\tau=\tau_p$ the frozen-income solution provides, through
\cref{lem:singleW} with $\varphi(\tau_p)=p$, an inverse profile
$\rho_0$ on $[0,p)$ satisfying
$\rho_0(w)\to+\infty$ as $w\uparrow p$.

The hodograph transformation below fixes the moving no-action region
by sending its lower and upper boundaries to $w=0$ and $w=p$,
respectively. In exchange, the fact that the upper boundary starts
at $H(\tau_p)=+\infty$ becomes a singular initial value,
$\rho_0(p-)=+\infty$. The main difficulty in the two-boundary regime
is therefore the initial corner $(h,w)=(0,p)$.

\subsubsection{The partial hodograph transformation}

Suppose for the moment that on the no-action region
$L(\tau)<\lambda<H(\tau)$ the function $W$ is smooth and strictly
increasing in $\zeta=\log\lambda$, with $W=0$ at $L$ and $W=p$ at
$H$. We invert it in the spatial variable as follows:
\begin{equation}\label{eq:hododef}
 w=W(\tau,\ee^{\zeta})\iff \zeta=\rho(h,w),\qquad
 h=\tau-\tau_p>0,\quad 0<w<p .
\end{equation}
Differentiating the identity
$W(\tau,\ee^{\rho(h,w)})=w$ and using $\PP W=1$ gives the
quasilinear equation
\begin{equation}\label{eq:hodoPDE}
 \rho_h
 =\frac{a\rho_{ww}}{\rho_w^2}
 -b_0-(1-rw)\rho_w
 =-a\Bigl(\frac1{\rho_w}\Bigr)_w
 -b_0-(1-rw)\rho_w .
\end{equation}
The free boundaries are the traces
$L=\ee^{\rho(h,0)}$ and $H=\ee^{\rho(h,p)}$.
At the lower boundary, \eqref{eq:lowercurv} gives
$W_\zeta=(1-L)/a$, that is,
$\rho_w=a/(1-\ee^\rho)$.
At the upper boundary the relevant condition is the smooth fit of
the value function, which in the dual variables reads
$Q=p$ and $Q_\lambda=0$ at $H$.
Since $\lambda Q_\lambda=\gamma(Q-W)$ by \eqref{eq:Qmeaning},
$Q_\lambda=0$ forces $W=p$ at $H$, and differentiating once more and
evaluating $\PP Q=1$ at $H$ gives
$W_\zeta(\tau,H-)=(1-rp)/(a\gamma)$.
We therefore impose
\begin{equation}\label{eq:hodoBC}
 \rho_w(h,0)=\frac{a}{1-\ee^{\rho(h,0)}},\qquad
 \rho_w(h,p)=c_H:=\frac{a\gamma}{1-rp},\qquad
 \rho(0,w)=\rho_0(w),\quad 0\le w<p .
\end{equation}
The derivation is only heuristic, and \cref{prop:reconPDE} shows that the
functions reconstructed from a solution of
\eqref{eq:hodoPDE}--\eqref{eq:hodoBC} satisfy all the free boundary
conditions of the original problem.
Note that the upper condition expresses smooth fit of $Q$, not of
$W$: the function $W$ has a kink at $H$ when it is extended by the
constant $p$.

The interval $(0,p)$ is now fixed, the principal coefficient
$a/\rho_w^2$ is positive wherever $\rho_w$ is positive and finite,
and the boundary conditions are oblique.

The remaining singularity is concentrated at $(0,p)$.
Indeed, by \cref{lem:singleW}, the frozen-income profile satisfies
$\rho_0(w)\to+\infty$ as $w\uparrow p$, since
$W(\tau_p,\lambda)\to p$ only as $\lambda\to\infty$.
Moreover,
$\rho_0'(w)=1/W_\zeta(\tau_p,\ee^{\rho_0(w)})\to+\infty$
because $W_\zeta\to0$ as $\zeta\to\infty$.
Thus the initial value and its slope are singular at $w=p$, while
for every $h>0$ the upper boundary condition prescribes the finite
slope $\rho_w(h,p)=c_H$.
To control this corner we introduce the time dependent scale
\begin{equation}\label{eq:G}
 G(h)=(1-rp)\varphi(h),\qquad
 G'(h)=1-rp-rG(h),\qquad
 G(0)=0 .
\end{equation}
It is adapted to the equation satisfied by $\rho_w$: for every
constant $C$, the function $C/(p-w+G(h))$ is an exact solution of
that equation (see \eqref{eq:qPDE} below).

\begin{theorem}\label{thm:hodo}
Let $0<rp<1$ and let $\rho_0\in C^\infty([0,p))$ satisfy
\begin{equation}\label{eq:rhodata}
\begin{aligned}
 &\rho_0(0)<0,\qquad
 \rho_0'(0)=\frac{a}{1-\ee^{\rho_0(0)}},\qquad
 \rho_0(w),\,\rho_0'(w)\to+\infty\quad
 \text{as }\ w\uparrow p,\\
 &0<m_0\le\rho_0'(w)\le\frac{C_0}{p-w}.
\end{aligned}
\end{equation}
For every $\bar h>0$ the problem \eqref{eq:hodoPDE}--\eqref{eq:hodoBC}
has a solution $\rho\in C^\infty((0,\bar h]\times[0,p])$ whose initial
trace is $\rho_0$, in $L^2(0,p)$ and locally uniformly on $[0,p)$.
There are constants $M,C,N,\underline q>0$ and
$0<\underline L\le\overline L<1$ such that on
$(0,\bar h]\times[0,p]$
\begin{align}
 &\underline L\le\ee^{\rho(h,0)}\le\overline L,\qquad
 \underline q\le\rho_w(h,w)\le\frac{C}{p-w+G(h)},\label{eq:hodobounds}\\
 &-M\le\rho(h,w)\le\log\overline L+
 N\log\frac{p+G(h)}{p-w+G(h)}.\label{eq:rhovaluebound}
\end{align}
Only $M$ and $\underline L$ depend on $\bar h$, while the constants
$C$, $N$, $\underline q$ and $\overline L$ depend only on the data.
The solution is unique among classical solutions on
$(0,\bar h]\times[0,p]$ with $L^2$ initial trace $\rho_0$ that satisfy
$\sup_h\rho(h,0)<0$ and bounds of the form \eqref{eq:hodobounds},
with possibly different constants. If, in addition, $\rho_0$
satisfies the compatibility conditions of every order for
\eqref{eq:hodoPDE}--\eqref{eq:hodoBC} at the corner $(0,0)$, then
$\rho\in C^\infty([0,\bar h]\times[0,p'])$ for every $p'<p$.
\end{theorem}

For the proof, we approximate the singular initial datum by truncating
it near $w=p$.

For all sufficiently large $n$, let
$\rho_{0,n}\in C^\infty([0,p])$ coincide with $\rho_0$ on
$[0,p-1/n]$ and satisfy $\rho_{0,n}'(p)=c_H$.
Since $\rho_0'(w)\to\infty$ as $w\uparrow p$, the modification on
$[p-1/n,p]$ may be chosen so that
$\rho_{0,n}'\le\rho_0'$, and hence $\rho_{0,n}\le\rho_0$.
Moreover,
\begin{equation}\label{eq:approxdata}
 \rho_{0,n}(0)=\rho_0(0),\qquad
 m_*\le\rho_{0,n}'(w)\le\frac{C_*}{p-w}
 \quad(0\le w<p)
\end{equation}
with $m_*=\min\{m_0,c_H\}$ and
$C_*=\max\{C_0,pc_H\}$, independent of $n$.
Each $\rho_{0,n}$ is bounded and satisfies the compatibility
conditions of order zero at both corners.

All comparison arguments for the hodograph problem in this paper rest
on the following energy estimate. It needs no strict inequalities at
the boundary and applies near the singular corner as long as the two
functions are ordered in a neighborhood of the corner.

\begin{lemma}\label{lem:comparison}
Let $g$ be a positive nondecreasing function on an interval $I$. Let $0\le h_a<h_b$, and let $\rho^1,\rho^2$ be of class $C^{2,1}$ on $(h_a,h_b]\times[0,p]$, with positive $w$-derivatives that are bounded above on $(h_a,h_b]\times[0,p/2]$, and with $\rho^1(h,0),\rho^2(h,0)\in I$. Suppose that
\[
 \rho^1_h\le\frac{a\rho^1_{ww}}{(\rho^1_w)^2}-b_0-(1-rw)\rho^1_w,\qquad \rho^1_w(h,0)\ge g\bigl(\rho^1(h,0)\bigr),\qquad \rho^1_w(h,p)\le c_H,
\]
that $\rho^2$ satisfies the reverse inequalities, and that $\|(\rho^1-\rho^2)_+(h,\cdot)\|_{L^2(0,p)}\to0$ as $h\downarrow h_a$. Then $\rho^1\le\rho^2$ on $(h_a,h_b]\times[0,p]$. In the hodograph problem $g(\rho)=a/(1-\ee^\rho)$ on $I=(-\infty,0)$.
\end{lemma}

\begin{proof}
Let $\omega=\rho^1-\rho^2$ and $A=a/(\rho^1_w\rho^2_w)$. Since $a\rho_{ww}/\rho_w^2=-a(1/\rho_w)_w$, subtracting the two inequalities gives $\omega_h\le(A\omega_w)_w-(1-rw)\omega_w$. Multiply by $\omega_+$ and integrate over $(0,p)$:
\[
 \frac12\frac{\mathrm d}{\mathrm dh}\|\omega_+\|^2\le\bigl[A\omega_w\omega_+\bigr]_{w=0}^{w=p}-\int_0^pA(\partial_w\omega_+)^2\dd w-\frac{1-rp}2\omega_+(h,p)^2+\frac12\omega_+(h,0)^2-\frac r2\|\omega_+\|^2 .
\]
At $w=p$, $\omega_w\le c_H-c_H=0$. At $w=0$, wherever $\omega>0$ we have $\rho^1_w\ge g(\rho^1)\ge g(\rho^2)\ge\rho^2_w$, so $\omega_w\ge0$. Both boundary terms are therefore nonpositive. On $[0,p/2]$, $A\ge A_0>0$ by the upper bounds on the slopes, and the trace inequality $\frac12\omega_+(h,0)^2\le\frac{A_0}{2}\int_0^{p/2}(\partial_w\omega_+)^2\dd w+C\|\omega_+\|^2$ leads to $\frac{\mathrm d}{\mathrm dh}\|\omega_+\|^2\le2C\|\omega_+\|^2$. Gronwall's lemma and the assumption at $h_a$ give $\omega_+\equiv0$.
\end{proof}

\begin{lemma}\label{lem:truncated}
For every $n$ and $\bar h>0$ the problem
\eqref{eq:hodoPDE}--\eqref{eq:hodoBC} with initial datum
$\rho_{0,n}$ has a unique solution
$\rho_n\in C^{2+\alpha,1+\alpha/2}([0,\bar h]\times[0,p])$
for some $\alpha\in(0,1)$, and
$\rho_n\in C^\infty((0,\bar h]\times[0,p])$.
\end{lemma}

\begin{proof}
Uniqueness follows from \cref{lem:comparison} applied in both directions.

For existence, fix $n$, let $\underline q=\min\{m_*,a,c_H\}$ and put
\[
 \overline q_n=\ee^{r\bar h}
 \max\Bigl\{
     \sup_{[0,p]}\rho_{0,n}',\
     \frac{a}{1-\overline L},\
     c_H
 \Bigr\}.
\]
We first regularize the coefficients outside the range of gradients
relevant to the solution. Choose a smooth positive function $A_n$,
bounded above and below by positive constants, such that
\[
    A_n(s)=\frac{a}{s^2}
    \qquad\text{for }\underline q/2\le s\le2\overline q_n,
\]
and replace $a\rho_{ww}/\rho_w^2$ by $A_n(\rho_w)\rho_{ww}$. We also
replace the boundary nonlinearity $g(\rho)=a/(1-\ee^\rho)$ by a smooth
nondecreasing function $\tilde g\ge a$ with bounded derivatives that
coincides with $g$ on $(-\infty,\log\overline L_1]$, where
$\overline L_1=(1+\overline L)/2$ and $\overline L<1$ is the
constant of the upper barrier below. The resulting equation is
uniformly parabolic, and the boundary operators
\[
    \rho_w-\tilde g(\rho)=0\quad\text{at }w=0,
    \qquad
    c_H-\rho_w=0\quad\text{at }w=p
\]
are uniformly oblique. Since $\rho_{0,n}$ is smooth and satisfies
the compatibility conditions of order zero,
\[
    \rho_{0,n}'(0)=\tilde g\bigl(\rho_{0,n}(0)\bigr),
    \qquad
    \rho_{0,n}'(p)=c_H,
\]
the local solvability argument for quasilinear oblique problems
\citep[Theorem 5.18 and the proof of Theorem 13.19]{Lieberman1996}
gives a solution in $C^{2+\alpha,1+\alpha/2}$ on a short time
interval.

Let $h^*$ be the supremum of the times $h'\le\bar h$ such that the
solution exists on $[0,h']$ and satisfies
$\underline q/2<\rho_w<2\overline q_n$ there. By continuity $h^*>0$.
On $[0,h^*)$ the regularized principal coefficient equals
$a/\rho_w^2$, so $\rho$ solves \eqref{eq:hodoPDE} with $\tilde g$ in
place of $g$ at $w=0$. The barriers $\overline\rho$ and
$\underline\rho$ constructed below take values at $w=0$ in
$(-\infty,\log\overline L]$, where $\tilde g=g$, so
\cref{lem:comparison}, with $g$ replaced by $\tilde g$ and $I=\R$,
gives
\[
    -M_0-C_1\bar h
    \le\rho(h,0)\le\log\overline L,\qquad 0\le h<h^*.
\]
In particular $\tilde g(\rho(h,0))=g(\rho(h,0))$, so the modified
boundary condition is inactive. Moreover $q=\rho_w$ satisfies
\eqref{eq:qPDE} on $[0,h^*)$, with boundary values
$q(h,0)=g(\rho(h,0))\in[a,a/(1-\overline L)]$ and $q(h,p)=c_H$.
At a first interior contact of $q$ with a level below $\underline q$
one would have $q_w=0$, $q_{ww}\ge0$ and $q_h\le0$, whereas
\eqref{eq:qPDE} gives $q_h>0$. The function $\ee^{-rh}q$ solves an
equation without zero order term, and the weak maximum principle
\citep[Theorem 2.4]{Lieberman1996} bounds it by its initial and
boundary values. Together,
\[
    \underline q\le\rho_w(h,w)\le\overline q_n,\qquad 0\le h<h^*.
\]
Hence $\rho_w$ stays away from the endpoints of
$(\underline q/2,2\overline q_n)$, and $h^*$ coincides with the
maximal existence time in $[0,\bar h]$. On the whole interval of
existence both regularizations are inactive, and the solution solves
the original problem \eqref{eq:hodoPDE}--\eqref{eq:hodoBC}.

The bounds on $\rho(h,0)$ and $\rho_w$ also bound $\rho$ itself and
keep the equation uniformly parabolic for fixed $n$.
The one-dimensional H\"older gradient estimate for quasilinear
oblique problems described in
\citet[Chapter XIII, Section 7, pp.~356--357]{Lieberman1996},
followed by the nonlinear oblique regularity estimate
\citep[Theorem 8.9]{Lieberman1996}, gives a
$C^{2+\alpha,1+\alpha/2}$ bound up to the lateral boundaries
$w=0$ and $w=p$. These estimates prevent loss of regularity at a
finite continuation time, and the continuation argument of
\citet[Theorem 8.3]{Lieberman1996}, as used in the proof of
\citet[Theorem 13.19]{Lieberman1996}, extends the solution to
$[0,\bar h]\times[0,p]$.

Finally, fix $\delta>0$. On
$[\delta/2,\bar h]\times[0,p]$ the equation is uniformly parabolic
with smooth coefficients and both boundary conditions are smooth and
uniformly oblique. After localization away from $h=0$, successive
differentiation of the equation and the boundary conditions, together
with the linear oblique Schauder estimate
\citep[Theorem 5.18]{Lieberman1996}, gives bounds of every order up
to both lateral boundaries on $[\delta,\bar h]\times[0,p]$.
Since $\delta>0$ is arbitrary,
$\rho_n\in C^\infty((0,\bar h]\times[0,p])$.
\end{proof}

\subsubsection{Uniform estimates using barriers}
We next derive bounds for the truncated solutions $\rho_n$ that do not
depend on $n$, beginning with the bounds on $\rho_n$ itself obtained
from explicit barriers.

\begin{proof}[Proof of \cref{thm:hodo}: bounds on $\rho_n$]
Let $\rho_n$ be the solutions of \cref{lem:truncated}. The bounds below hold for any classical solution with initial datum $\rho_{0,n}$ on any subinterval $[0,h']$, as used in that lemma. Write $G_\infty=\sup_hG(h)=(1-rp)/r$.

\emph{Upper barrier.} Choose $N\ge C_*$ so large that $rN+b_0-a/N\ge0$ and $N\ge c_HG_\infty$, and then $\overline L\in(\ee^{\rho_0(0)},1)$ so close to one that $N/p\le a/(1-\overline L)$. Put
\begin{equation}\label{eq:rhobar}
 \overline\rho(h,w)=\log\overline L+N\log\frac{p+G(h)}{p-w+G(h)} .
\end{equation}
Using $G'=1-rp-rG$ one finds
\begin{equation}\label{eq:barresidual}
 \overline\rho_h-\frac{a\overline\rho_{ww}}{\overline\rho_w^2}+b_0+(1-rw)\overline\rho_w=\frac{NG'}{p+G}+rN+b_0-\frac aN\ge0 ,
\end{equation}
so $\overline\rho$ is a supersolution in the interior. At $w=0$, $\overline\rho=\log\overline L$ and $\overline\rho_w=N/(p+G)\le a/(1-\overline L)$, which is the supersolution inequality at $w=0$ required in \cref{lem:comparison}. At $w=p$, $\overline\rho_w=N/G(h)\ge c_H$. At $h=0$, integration of \eqref{eq:approxdata} gives $\rho_{0,n}\le\overline\rho(0,\cdot)$. The barrier is infinite at the corner $(0,p)$, whereas $\rho_n$ is bounded, so $(\rho_n-\overline\rho)_+$ vanishes near the corner and tends to zero in $L^2$ as $h\downarrow0$. \Cref{lem:comparison} gives $\rho_n\le\overline\rho$.

\emph{Lower barrier.} Let $f_0$ be smooth with $f_0(0)=0$, $f_0'>0$, $f_0'(0)=2a$ and $f_0'(p)=c_H/2$, for example with $f_0'$ affine. Put $\underline\rho(h,w)=-M_0-C_1h+f_0(w)$ with $M_0\ge\log2$ so large that $\underline\rho(0,\cdot)\le\rho_{0,n}$ for all $n$, and $C_1$ so large that
\[
 -C_1-\frac{af_0''}{(f_0')^2}+b_0+(1-rw)f_0'\le0\quad\text{on }[0,p].
\]
At $w=0$, $2a\ge a/(1-\ee^{-M_0-C_1h})$ because $\ee^{-M_0}\le1/2$, and at $w=p$, $\underline\rho_w=c_H/2\le c_H$. Hence $\underline\rho$ is a subsolution and \cref{lem:comparison} gives $\underline\rho\le\rho_n$. In particular,
\[
 \rho_n(h,0)\in[-M_0-C_1\bar h,\log\overline L].
\]
Thus the corresponding bounds on $\ee^{\rho_n(h,0)}$ and $\rho_n$,
including the right-hand side of \eqref{eq:rhovaluebound}, hold
uniformly in $n$. Passing to the limit $n\to\infty$ in the convergence
step below yields the stated bounds for $\rho$.
\end{proof}

\subsubsection{Uniform parabolicity away from the corner}
\begin{proof}[Proof of \cref{thm:hodo}: gradient bounds]
We next derive upper and lower bounds for $\rho_{n,w}$ that do not
depend on $n$. Since the principal coefficient in \eqref{eq:hodoPDE}
is $a/\rho_{n,w}^2$, they yield uniform parabolicity away from the
singular initial corner $(0,p)$.

Let $q_n=\partial_w\rho_n$. Differentiating \eqref{eq:hodoPDE},
\begin{equation}\label{eq:qPDE}
 q_h=a\Bigl(\frac{q_{ww}}{q^2}-\frac{2q_w^2}{q^3}\Bigr)-(1-rw)q_w+rq ,
\end{equation}
with Dirichlet data
$q_n(h,0)=a/(1-\ee^{\rho_n(h,0)})\in[a,a/(1-\overline L)]$
and $q_n(h,p)=c_H$. The first contact argument in the proof of
\cref{lem:truncated} gives $q_n\ge\underline q=\min\{m_*,a,c_H\}$,
which does not depend on $n$. The upper bound obtained there does
depend on $n$, through $\sup\rho_{0,n}'$.

The decisive estimate is independent of $n$. Let
\[
 C\ge\max\Bigl\{C_*,\ \frac{a\,(p+G_\infty)}{1-\overline L},\ c_HG_\infty\Bigr\},\qquad
 \overline q(h,w)=\frac{C}{p-w+G(h)} .
\]
For this function the diffusion term and the quadratic gradient term in \eqref{eq:qPDE} cancel,
$a(\overline q_{ww}/\overline q^2-2\overline q_w^2/\overline q^3)=0$,
and the relation $G'=1-rp-rG$ shows that $\overline q$ solves
\eqref{eq:qPDE} exactly. It dominates the initial values on $[0,p)$
and both boundary values, and it is infinite at the corner, whereas
$q_n$ is bounded for fixed $n$. Hence $q_n<\overline q$ in a
neighborhood of $(0,p)$. On the complement of such a neighborhood the
difference $q_n-\overline q$ satisfies a linear parabolic equation
whose coefficients, obtained from \eqref{eq:qPDE} by the mean value
theorem, are bounded. The classical maximum principle for linear
equations with Dirichlet data, applied to
$\ee^{-c_0h}(q_n-\overline q)$ with $c_0$ larger than the zero order
coefficient, gives
\begin{equation}\label{eq:qbound}
 \underline q\le q_n(h,w)\le\frac{C}{p-w+G(h)} .
\end{equation}
For every $\delta>0$ the equation for $\rho_n$ is thus uniformly
parabolic on $[\delta,\bar h]\times[0,p]$, uniformly in $n$, and the
lower boundary nonlinearity and its derivatives are bounded because
$\rho_n(h,0)\le\log\overline L<0$. H\"older gradient estimates for
quasilinear oblique problems, Schauder estimates and differentiation
of the equation give bounds on all derivatives of $\rho_n$ on
$[\delta,\bar h]\times[0,p]$, uniform in $n$. On
$[0,\bar h]\times[0,p']$ with $p'<p$ the same holds down to $h=0$
under the additional compatibility assumption of \cref{thm:hodo}.
Fix $p''\in(p',p)$. On $[0,p'']$ we have
$q_n\le C/(p-p'')$, and the initial data are smooth and independent
of $n$ for $1/n<p-p''$. If they satisfy the compatibility conditions
of every order at the corner $(0,0)$, estimates that are local in $w$
on $[0,p'']$ bound all derivatives on $[0,\bar h]\times[0,p']$. The
profile of \cref{lem:singleW} with $\tau_0=\tau_p$ satisfies these
conditions, because it is the trace at $\tau_p$ of the smooth
frozen-income hodograph, which solves the same boundary value problem
near $w=0$ (\cref{thm:single}).
\end{proof}

\subsubsection{Convergence and uniqueness}
\begin{proof}[Proof of \cref{thm:hodo}: completion]
Let $q_i=\partial_w\rho_i$ for $i=n,j$ and
$A=a/(q_nq_j)$. The divergence form of \eqref{eq:hodoPDE}
gives
\begin{equation}\label{eq:diffPDE}
 \partial_h(\rho_n-\rho_j)
 =(A\partial_w(\rho_n-\rho_j))_w
 -(1-rw)\partial_w(\rho_n-\rho_j).
\end{equation}
At $w=p$, $\partial_w(\rho_n-\rho_j)=0$.
At $w=0$ the boundary condition yields the exact identity
\begin{equation}\label{eq:diffBC}
 A\partial_w(\rho_n-\rho_j)
 =\frac{a}{q_j}-\frac{a}{q_n}
 =\ee^{\rho_n(h,0)}-\ee^{\rho_j(h,0)}
 =\vartheta(h)(\rho_n-\rho_j)(h,0),
 \qquad 0<\vartheta(h)<1,
\end{equation}
where $\vartheta$ is a difference quotient of the exponential
on $(-\infty,0)$.

Multiplying \eqref{eq:diffPDE} by $\rho_n-\rho_j$
and integrating by parts, we obtain
\begin{equation}\label{eq:energy}
\begin{aligned}
 &
 \frac12\frac{\mathrm d}{\mathrm dh}
 \|\rho_n-\rho_j\|_{L^2(0,p)}^2
 +\int_0^p A
 |\partial_w(\rho_n-\rho_j)|^2\dd w\\
 &\quad
 +
 \frac{1-rp}{2}(\rho_n-\rho_j)^2(h,p)
 +\frac r2\|\rho_n-\rho_j\|_{L^2(0,p)}^2\\
 &=
 \Bigl(\frac12-\vartheta(h)\Bigr)
 (\rho_n-\rho_j)^2(h,0).
\end{aligned}
\end{equation}
On $[0,p/2]$, \eqref{eq:qbound} gives
$A\ge A_0:=ap^2/(4C^2)$, uniformly in $h$ and in $n,j$.

The trace inequality
\[
 \frac12(\rho_n-\rho_j)^2(h,0)
 \le\frac{A_0}{2}
 \int_0^{p/2}|\partial_w(\rho_n-\rho_j)|^2\dd w
 +C_2\|\rho_n-\rho_j\|_{L^2(0,p)}^2
\]
and Gronwall's lemma give
\begin{equation}\label{eq:energyconvergence}
 \sup_{0\le h\le\bar h}
 \|\rho_n(h,\cdot)-\rho_j(h,\cdot)\|_{L^2(0,p)}^2
 \le
 \ee^{2C_2\bar h}
 \|\rho_{0,n}-\rho_{0,j}\|_{L^2(0,p)}^2 .
\end{equation}
The bounds
$\rho_0(0)+m_*w\le\rho_{0,n}(w)
\le\rho_0(0)+C_*\log\frac{p}{p-w}$
and pointwise convergence imply
$\rho_{0,n}\to\rho_0$ in $L^2(0,p)$.
Hence $(\rho_n)$ is Cauchy in
$C([0,\bar h];L^2(0,p))$.
By the uniform estimates of the previous step,
a subsequence converges smoothly on every
$[\delta,\bar h]\times[0,p]$, and the limit $\rho$
solves \eqref{eq:hodoPDE}--\eqref{eq:hodoBC} for $h>0$,
satisfies \eqref{eq:hodobounds}--\eqref{eq:rhovaluebound},
and has initial trace $\rho_0$ in $L^2$ and
locally uniformly on $[0,p)$.

For uniqueness, apply the same computation to two solutions
in the stated class on $[\delta,\bar h]$ and let
$\delta\downarrow0$. Only the lower bound on $A$ near $w=0$
enters, and it is uniform down to $h=0$. The degeneration
of the coefficient near the corner plays no role.
\end{proof}

We now recover the free boundaries from the hodograph solution.
The lower free boundary $L^0$ of the frozen-income problem is
already smooth for $\tau\le\tau_p$, while the hodograph
construction gives the boundaries $L$ and $H$ for
$\tau>\tau_p$. To show that $L$ and $L^0$ join smoothly at $\tau_p$,
we combine the regularity of $L^0$ and of the frozen-income solution
up to its free boundary with the regularity of the hodograph solution
up to the initial time near $w=0$.

\begin{corollary}\label{cor:bothboundaries}
Let $L^0$ denote the free boundary of the frozen-income problem
in \cref{thm:single}. Let $\rho_0$ be the profile of
\cref{lem:singleW} with $\tau_0=\tau_p$, let $\rho$ be the
solution of \cref{thm:hodo}, and for $\tau>\tau_p$ put
\begin{equation}\label{eq:LHhodo}
 L(\tau)=\exp\rho(\tau-\tau_p,0),\qquad
 H(\tau)=\exp\rho(\tau-\tau_p,p).
\end{equation}
Then $L,H\in C^\infty((\tau_p,\infty))$, $0<L<1$,
and the two
boundaries are uniformly separated in logarithmic coordinates:
\[
 \log H-\log L=\int_0^p\rho_w\dd w
 \ge\underline q\,p>0 .
\]
Moreover, $H(\tau)\to\infty$ as $\tau\downarrow\tau_p$ and,
for each $\bar h>0$,
$H(\tau)\le C_{\bar h}(\tau-\tau_p)^{-N}$ on
$(\tau_p,\tau_p+\bar h]$, with $N$ from \cref{thm:hodo}
and a constant $C_{\bar h}$.
Finally, let $\rho^0(h,\cdot)$ be the inverse of
$\zeta\mapsto W^0(\tau_p+h,\ee^\zeta)$ on $[0,\varphi(\tau_p+h))$, so
that $\rho^0(0,\cdot)=\rho_0$. For every $p'<p$ there is $h_0>0$ such
that the function equal to $\rho^0$ for $-h_0<h\le0$ and to $\rho$
for $0<h\le\bar h$ is $C^\infty$ on $(-h_0,\bar h]\times[0,p']$. In
particular, $L$ extends to $(0,\infty)$ by setting
$L(\tau)=L^0(\tau)$ for $\tau\le\tau_p$, and this extension belongs to
$C^\infty((0,\infty))$.
\end{corollary}

\begin{proof}
For every $\delta>0$, \cref{thm:hodo} gives smoothness of $\rho$
up to both lateral boundaries on $[\delta,\bar h]\times[0,p]$.
Since
$L(\tau)=\exp\rho(\tau-\tau_p,0)$ and
$H(\tau)=\exp\rho(\tau-\tau_p,p)$,
this gives $L,H\in C^\infty((\tau_p,\infty))$.
The only singular point associated with the upper boundary is the
initial right corner $(h,w)=(0,p)$.
Moreover, by the definitions of $L$ and $H$,
\[
 \log H(\tau)-\log L(\tau)
 =\rho(h,p)-\rho(h,0)
 =\int_0^p\rho_w(h,w)\,\dd w
 \ge\underline q\,p>0,
\]
where $h=\tau-\tau_p$. Thus the two boundaries are separated.
The solutions for different $\bar h$ coincide by uniqueness,
so $L$ and $H$ are defined on $(\tau_p,\infty)$.
For fixed $w<p$,
$\rho(h,p)\ge\rho(h,w)\to\rho_0(w)$ as $h\downarrow0$,
and $\rho_0(w)\to\infty$ as $w\uparrow p$.
Hence $H\to\infty$.
Evaluating \eqref{eq:rhovaluebound} at $w=p$ and using
$G(h)\ge(1-rp)\ee^{-r\bar h}h$ gives the polynomial bound.

It remains to verify the smooth joining at $\tau_p$. Fix $p'<p$ and
choose $h_0>0$ with $\varphi(\tau_p-h_0)>p'$. For $h\le0$, the
regularity of the free boundary of the frozen-income problem, together
with boundary Schauder estimates, gives smoothness of $W^0$ up to its
lower free boundary. Since $W^0_\zeta>0$ on the no-action region and
\[
 W^0_\zeta\bigl(\tau,L^0(\tau)+\bigr)
 =\frac{1-L^0(\tau)}{a}>0,
\]
the inverse function theorem shows that $\rho^0$ is smooth on
$(-h_0,0]\times[0,p']$, with
\[
 \rho^0(h,0)=\log L^0(\tau_p+h),
 \qquad
 \rho^0(0,w)=\rho_0(w).
\]
For $h\ge0$, \cref{thm:hodo} applies with the compatibility conditions
of every order, as shown at the end of its proof, and gives
$\rho\in C^\infty([0,\bar h]\times[0,p'])$. Both $\rho$ and $\rho^0$
satisfy \eqref{eq:hodoPDE} on $[0,p']$ and the lower boundary
condition in \eqref{eq:hodoBC}, and they coincide at $h=0$. Each time
derivative at a point $(0,w)$, $0\le w\le p'$, is therefore determined
through \eqref{eq:hodoPDE} by the $w$ derivatives of $\rho_0$ at $w$,
and repeated differentiation gives
\[
 \partial_h^k\rho^0(0-,w)
 =
 \partial_h^k\rho(0+,w),
 \qquad k\ge0,\ 0\le w\le p'.
\]
Differentiating these identities in $w$ shows that all mixed
derivatives match as well, so the glued function is $C^\infty$ on
$(-h_0,\bar h]\times[0,p']$. At $w=0$ this says that all one-sided
derivatives of $\log L$ agree at $\tau_p$, and the extension of $L$
belongs to $C^\infty((0,\infty))$.
\end{proof}

In the numerical example of \cref{sec:numerics} the purchase boundary grows much more slowly than the polynomial bound allows. The reason is the Gaussian tail \eqref{eq:tailkernel} of the frozen-income solution: $p-W(\tau_p,\ee^\zeta)$ decays roughly like $\exp(-\zeta^2/(4a\tau_p))$, so the initial profile grows only like $\rho_0(w)=O\bigl(\sqrt{\log(1/(p-w))}\bigr)$ near the corner, and the region where purchases start is reached at moderate values of $\lambda$ even for small $\tau-\tau_p$.

\subsection{Reconstruction and regularity}\label{sec:reconstruction}

Recall that $W=\lambda\X_\lambda-\X$ is the Euler derivative
introduced in \eqref{eq:XQF}, that $Q=V_\xi/V_x$, and that $F$ is the
scaled dual value function \eqref{eq:Fdef}. Having constructed the
hodograph $\rho$, we now reconstruct $W$ by inverting it and then
recover $\X$, $Q$ and $F$ from $W$.

\subsubsection{One reconstruction formula for both regimes}
For $0<\tau\le\tau_p$ let $W$ be the function of \cref{lem:singleW} and $H(\tau)=\infty$. For $\tau>\tau_p$ define $W(\tau,\cdot)$ on $(L(\tau),H(\tau))$ by $W\bigl(\tau,\ee^{\rho(\tau-\tau_p,w)}\bigr)=w$, $0<w<p$, and extend it by
\begin{equation}\label{eq:Wextension}
 W(\tau,\lambda)=0\quad(\lambda\le L(\tau)),\qquad W(\tau,\lambda)=p\quad(\lambda\ge H(\tau)).
\end{equation}
In both regimes $W(\tau,\cdot)$ is continuous, nondecreasing and bounded by $p$, and strictly increasing on $(L,H)$. Define, for $\lambda>0$,
\begin{align}
 \X(\tau,\lambda)&=\begin{cases}0,&\lambda\le L(\tau),\\[2pt]
 \displaystyle\lambda\int_{L(\tau)}^\lambda\frac{W(\tau,\eta)}{\eta^2}\dd\eta,&\lambda>L(\tau),\end{cases}\label{eq:Xrecon}\\
 Q(\tau,\lambda)&=\gamma\lambda^\gamma\int_\lambda^\infty W(\tau,\eta)\,\eta^{-\gamma-1}\dd\eta,\qquad
 F(\tau,\lambda)=\frac{Q(\tau,\lambda)+\gamma\X(\tau,\lambda)}{1-\gamma}.\label{eq:QFrecon}
\end{align}
The integral defining $Q$ converges for every $\gamma>0$ because $W$ is bounded, so no divergent integral of the linearly growing $\X$ appears when $\gamma<1$. By construction $(D-1)\X=W$ and $(D-\gamma)Q=-\gamma W$. For $\tau\le\tau_p$, \eqref{eq:Xrecon} returns the solution of \cref{thm:single}, since both functions satisfy $(D-1)\X=W$ and vanish at $L$. All identities in \eqref{eq:XQF} and \eqref{eq:Qmeaning} follow by differentiation. Integration by parts gives
\begin{equation}\label{eq:QminusW}
 Q-W=\lambda^\gamma\int_\lambda^\infty W_\eta(\tau,\eta)\,\eta^{-\gamma}\dd\eta>0\quad\text{on }\Nr,
\end{equation}
with $W_\eta$ understood piecewise. Hence $Q_\lambda=\gamma(Q-W)/\lambda>0$ on $\Nr$. If $H<\infty$, then $Q=p$ on $\Ar$ and $Q<p$ below $H$. If $H=\infty$, the strict bound $W<\varphi(\tau)$ gives $Q<\varphi(\tau)\le p$ at every finite $\lambda$. On $\Zr$, $Q(\tau,\lambda)=Q(\tau,L)(\lambda/L)^\gamma$. In particular, at $\tau=\tau_p$ the gradient constraint is still strictly inactive at every finite state, consistent with the purchase boundary entering from infinity.

\begin{proposition}\label{prop:reconPDE}
On $\Nr\cap\{\tau\ne\tau_p\}$ the functions
\eqref{eq:Xrecon}--\eqref{eq:QFrecon} satisfy
\eqref{eq:linearPDEs}, and
\begin{align}
 &\X(\tau,L)=\X_\lambda(\tau,L)=0,\qquad
 \X_{\lambda\lambda}(\tau,L+)=\frac{1-L}{aL^2},
 \label{eq:lowerall}\\
 &Q(\tau,H)=p,\qquad Q_\lambda(\tau,H)=0,\qquad
 H\X_\lambda(\tau,H)-\X(\tau,H)=p,
 \label{eq:upperall}\\
 &Q_{\lambda\lambda}(\tau,H-)
 =-\frac{1-rp}{aH^2},\qquad
 \X_{\lambda\lambda}(\tau,H-)
 =\frac{1-rp}{a\gamma H^2}.
 \label{eq:uppercurv}
\end{align}
The conditions at $H$ refer to $\tau>\tau_p$.
\end{proposition}

\begin{proof}
For $\tau\le\tau_p$ the equation for $\X$ is part of
\cref{thm:single}. Writing
$Q=\gamma\int_0^\infty\ee^{-\gamma u}
W(\tau,\lambda\ee^u)\dd u$, every integration point lies in
$\Nr$ when $\lambda>L$, and the tail estimates of
\cref{lem:singleW} allow $\PP$ to be taken inside the integral.
Hence
$\PP Q=\gamma\int_0^\infty\ee^{-\gamma u}\dd u=1$.

For $\tau>\tau_p$, inverting \eqref{eq:hodoPDE} gives
$\PP W=1$ on $\Nr$, and \eqref{eq:hodoBC} gives
\begin{equation}\label{eq:Wsboundaries}
 W_\zeta(\tau,L+)=\frac{1-L}{a},\qquad
 W_\zeta(\tau,H-)=\frac{1-rp}{a\gamma}.
\end{equation}
Since $\X_{\lambda\lambda}=W_\zeta/\lambda^2$, this gives the
curvatures of $\X$ in \eqref{eq:lowerall} and
\eqref{eq:uppercurv}. The values $\X=\X_\lambda=0$ at $L$ hold
by construction, and differentiating $\X(\tau,L(\tau))=0$ gives
$\X_\tau(\tau,L)=0$.

Since $W=\lambda\X_\lambda-\X$ and $\PP$ commutes with
$\lambda\partial_\lambda$, the identity $\PP W=1$ gives
\[
 \lambda\partial_\lambda(\PP\X-\lambda+1)
 -(\PP\X-\lambda+1)=0.
\]
Hence $\PP\X-\lambda+1=c(\tau)\lambda$ on $\Nr$. At
$\lambda=L$ we have
$\PP\X=-aL^2\X_{\lambda\lambda}(\tau,L+)=L-1$, and therefore
$c=0$.

Similarly, since
$\lambda Q_\lambda-\gamma Q=-\gamma W$,
\[
 \lambda\partial_\lambda(\PP Q-1)
 -\gamma(\PP Q-1)=0,
\]
and hence $\PP Q-1=c_1(\tau)\lambda^\gamma$. Since
$W=p$ on $[H,\infty)$, \eqref{eq:QFrecon} gives
$Q(\tau,H)=p$, and
$\lambda Q_\lambda=\gamma(Q-W)$ gives
$Q_\lambda(\tau,H)=0$. Differentiating
$Q(\tau,H(\tau))=p$ yields $Q_\tau(\tau,H)=0$.
Differentiating
$\lambda Q_\lambda=\gamma(Q-W)$ with respect to $\lambda$ and
evaluating at $H$ gives
\[
 H^2Q_{\lambda\lambda}(\tau,H-)
 =-\gamma W_\zeta(\tau,H-)
 =-\frac{1-rp}{a}.
\]
Therefore $\PP Q(\tau,H)=1$, and hence $c_1=0$.

Finally, from
$F=(Q+\gamma\X)/(1-\gamma)$ we obtain
\[
 \PP F
 =\frac{1+\gamma(\lambda-1)}{1-\gamma}
 =\frac{\gamma}{1-\gamma}\lambda+1.
\]
The remaining identity follows from
$H\X_\lambda-\X=W(\tau,H)=p$.
\end{proof}

\subsubsection{Verification of the global dual variational inequality}

We now verify that the reconstructed functions satisfy the dual
problem outside the no-action region. Since $W=p$ on the purchase
region $\Ar$, the relation $W=\lambda\X_\lambda-\X$ implies that $\X$
is affine in $\lambda$ on $\Ar$. We denote its slope by $d(\tau)$.
This coefficient also determines the purchase threshold and the
explicit form of $F$ on the purchase region. The evolution equation
for $d$ then gives the sign of the residual of the dual equation on
$\Ar$ and completes the verification of the global variational
inequality.

For $\tau>\tau_p$ let
\begin{equation}\label{eq:dB}
 d(\tau)=\int_{L(\tau)}^\infty\frac{W(\tau,\eta)}{\eta^2}\dd\eta>0,\qquad \zb(\tau)=\X(\tau,H(\tau))=d(\tau)H(\tau)-p>0,
\end{equation}
and let $d(\tau)=1/K$ for $0<\tau\le\tau_p$, which is the limit of $\X(\tau,\lambda)/\lambda$ as $\lambda\to\infty$ by \cref{thm:single}(a). On the purchase region \eqref{eq:Xrecon}--\eqref{eq:QFrecon} reduce to
\begin{equation}\label{eq:actionF}
 \X=d\lambda-p,\qquad Q=p,\qquad F=p+\frac{\gamma d}{1-\gamma}\lambda\qquad(\lambda\ge H).
\end{equation}
Dominated convergence in \eqref{eq:dB}, with $L$ bounded below and $0\le W\le p$, shows that $d$ is continuous on $(0,\infty)$ and smooth on $(\tau_p,\infty)$, and that
\begin{equation}\label{eq:escape}
 d(\tau)\to\frac1K,\qquad \zb(\tau)\to\infty\qquad(\tau\downarrow\tau_p).
\end{equation}
Together with \cref{cor:bothboundaries} this gives $\zb(\tau)\le d(\tau)H(\tau)\le C(\tau-\tau_p)^{-N}$ near $\tau_p$. Differentiating $\X(\tau,H(\tau))=d(\tau)H(\tau)-p$ and using $\X_\lambda(\tau,H)=d$ (a consequence of $H\X_\lambda-\X=p$) gives $\X_\tau(\tau,H-)=d'H$. Evaluating $\PP\X=\lambda-1$ at $H-$ with \eqref{eq:uppercurv} then yields
\begin{equation}\label{eq:dODE}
 d'(\tau)=1-Kd(\tau)+\frac{1-\gamma}{\gamma}\,\frac{1-rp}{H(\tau)}\qquad(\tau>\tau_p).
\end{equation}

The verification in \cref{sec:verification} uses only the equations on
$\Nr$ and the boundary conditions. For completeness we check that $F$
also solves the dual variational inequality associated with
\eqref{eq:primalHJB}. Let
$f(\lambda)=\gamma\lambda/(1-\gamma)+1$ and
$\mathcal R=\PP F-f$. On $\Nr$, $\mathcal R=0$. On $\Zr$,
$F=F(\tau,L)(\lambda/L)^\gamma$, and the relation
\[
 \lambda\mathcal R_\lambda-\gamma\mathcal R
 =\gamma(1-\lambda),
\]
together with $\mathcal R(\tau,L)=0$ by the continuity of
$\PP F$ from \cref{prop:regularity}, gives, with
$\hat\lambda=\lambda/L\in(0,1)$,
\begin{equation}\label{eq:lowerresidual}
 \mathcal R=\hat\lambda^\gamma-1
 +\frac{\gamma L}{1-\gamma}
 \bigl(\hat\lambda^\gamma-\hat\lambda\bigr).
\end{equation}
The right-hand side is increasing in $L$ for both ranges of $\gamma$, and at $L=1$ it equals $(\hat\lambda^\gamma-\gamma\hat\lambda+\gamma-1)/(1-\gamma)\le0$ by the tangent line inequality for $\hat\lambda\mapsto\hat\lambda^\gamma$. Since $L<1$, $\mathcal R\le0$ on $\Zr$. On $\Ar$, \eqref{eq:actionF} and \eqref{eq:dODE} give
\begin{equation}\label{eq:upperresidual}
 \mathcal R=(1-rp)\Bigl(\frac{\lambda}{H}-1\Bigr)\ge0 .
\end{equation}
Since $\gamma F-\lambda F_\lambda=-\gamma\X\le0$ everywhere, with equality exactly on $\Zr$, and $p-F+\lambda F_\lambda=p-Q\ge0$, with equality exactly on $\Ar$, we obtain the global variational inequality
\begin{equation}\label{eq:FVI}
 \min\Bigl\{\max\bigl\{\PP F-\tfrac{\gamma\lambda}{1-\gamma}-1,\ \gamma F-\lambda F_\lambda\bigr\},\ p-F+\lambda F_\lambda\Bigr\}=0,\qquad F(0,\lambda)=\frac{\gamma\lambda}{(1-\gamma)K}.
\end{equation}
The first constraint is the dual form of the borrowing constraint and the second that of the purchase option.

\subsubsection{Regularity of the reconstructed functions and the primal value}

We now establish the regularity of the reconstructed functions.
We first verify the matching of the dual quantities across the free
boundaries $L$ and $H$, and then across the transition time
$\tau_p$, where the construction changes from the frozen-income
solution to the hodograph reconstruction. These results are then
transferred to the primal variables through
$z=\X(\tau,\lambda)$. For $z>0$ the corresponding point satisfies
$\lambda>L(\tau)$, where $\X_\lambda>0$, so the change of variables
is regular. In particular, the matching across $H$ gives the
regularity of $v$ across the purchase threshold. We finally verify
the primal HJB equation.

\begin{proposition}\label{prop:regularity}
For $\tau>0$ the function $F$ is $C^{1,2}$ in $(\tau,\lambda)$
across the lower free boundary $L$ and, for
$\tau>\tau_p$, across the upper free boundary $H$.
The functions $\X$ and $Q$ are $C^1$ in $\lambda$ with continuous
time derivatives and belong to $W^{1,2}_{\mathfrak q,\mathrm{loc}}$ for every
$\mathfrak q<\infty$. Their second derivatives in $\lambda$ jump at $L$
and, for $\tau>\tau_p$, at $H$. On each side of a
free boundary all functions are smooth up to the boundary. The
primal function
\begin{equation}\label{eq:vrecon}
 v(\tau,z)=\lambda^{-\gamma}\bigl(F(\tau,\lambda)+z\bigr)
 =\frac{\lambda^{-\gamma}\bigl(z+Q(\tau,\lambda)\bigr)}
 {1-\gamma},\qquad z=\X(\tau,\lambda),
\end{equation}
extended by $\lambda=(z+p)/d(\tau)$ for $z\ge\zb(\tau)$, belongs to
$C^{1,2}((0,\infty)\times(0,\infty))$ and solves
\eqref{eq:primalHJB}. It satisfies
\begin{equation}\label{eq:vderivs}
 v_z=\lambda^{-\gamma}>0,\qquad
 v_{zz}=-\frac{\gamma\lambda^{-\gamma-1}}{\X_\lambda}<0,\qquad
 v(\tau,z)=\frac{d(\tau)^\gamma(z+p)^{1-\gamma}}{1-\gamma}
 \quad(z\ge\zb(\tau)).
\end{equation}
\end{proposition}

\begin{proof}
We first consider the free boundaries in the dual variables.
For $\tau\le\tau_p$, the identities at the lower boundary are contained
in \cref{thm:single}(c). For $\tau>\tau_p$,
\eqref{eq:Xrecon} and $W(\tau,L)=0$ give
$\X(\tau,L)=\X_\lambda(\tau,L)=0$. Since $L$ is smooth,
differentiating $\X(\tau,L(\tau))=0$ gives
$\X_\tau+\X_\lambda L'=0$, and hence
$\X_\tau(\tau,L)=0$. Thus $\X$ and its first derivatives match
the zero extension on $\{\lambda<L\}$. By \eqref{eq:lowerall},
$\X_{\lambda\lambda}(\tau,L+)=(1-L)/(aL^2)>0$, while
$\X_{\lambda\lambda}=0$ on the zero region, so the second spatial
derivative jumps at $L$.

At $H$, \eqref{eq:actionF} gives $\X=d\lambda-p$ on the purchase
region. Using \eqref{eq:upperall} and
$\X(\tau,H)=dH-p$ gives
$\X_\lambda(\tau,H-)=d=\X_\lambda(\tau,H+)$. Differentiating
$\X(\tau,H(\tau))=d(\tau)H(\tau)-p$ then gives
$\X_\tau(\tau,H-)=d'H=\X_\tau(\tau,H+)$. The second derivative
jumps by \eqref{eq:uppercurv}, since the extension to the purchase
region is affine.

For $Q$, on the zero region
$Q(\tau,\lambda)=Q(\tau,L)(\lambda/L)^\gamma$, while on the
no-action region
$\lambda Q_\lambda=\gamma(Q-W)$. Since $W(\tau,L)=0$, the
spatial derivatives match at $L$. Differentiating the value
matching along $\lambda=L(\tau)$ then gives the matching of
$Q_\tau$. At $H$, we have $Q=W=p$, so the same identity gives
$Q_\lambda(\tau,H-)=0$, matching the constant extension $Q=p$.
Differentiating $Q(\tau,H(\tau))=p$ gives the matching of
$Q_\tau$.

The second derivatives of $Q$ jump as well. On the no-action
region,
$\lambda^2Q_{\lambda\lambda}
=(\gamma-1)\gamma(Q-W)-\gamma W_\zeta$.
At $L$, the additional term
$-\gamma W_\zeta$ is nonzero because
$W_\zeta(\tau,L+)=(1-L)/a>0$, while on the zero region
$\lambda^2Q_{\lambda\lambda}=\gamma(\gamma-1)Q$.
At $H$, \eqref{eq:uppercurv} gives
$Q_{\lambda\lambda}(\tau,H-)=-(1-rp)/(aH^2)$, whereas
$Q_{\lambda\lambda}=0$ on the purchase region.

Since $F=(Q+\gamma\X)/(1-\gamma)$, the matching of $\X$, $Q$ and
their time derivatives gives the matching of $F$ and $F_\tau$.
Moreover,
$F_\lambda=\gamma(\X+F)/\lambda$ and
$F_{\lambda\lambda}=-Q_\lambda/\lambda$, so
$F_\lambda$ and $F_{\lambda\lambda}$ also match across the free
boundaries. Hence $F$ is $C^{1,2}$ across $L$ and $H$.
For $\X$ and $Q$, the first derivatives are continuous and the
one-sided second derivatives are locally bounded, which gives
$\X,Q\in W^{1,2}_{\mathfrak q,\mathrm{loc}}$ for every $\mathfrak q<\infty$.

We next consider the transition time $\tau_p$. Fix a finite
$\lambda>L(\tau_p)$ and set
$w_0=W(\tau_p,\lambda)<p$. Choose $p'$ with $w_0<p'<p$.
By \cref{cor:bothboundaries}, the frozen-income hodograph
$\rho^0$ and the hodograph $\rho$ join smoothly across $h=0$ on
$[0,p']$. Since $\rho_w>0$, the corresponding inverse functions
also join smoothly. In particular, $W$ and
$W_\zeta=1/\rho_w$ are continuous across $\tau=\tau_p$ at every
fixed finite $\lambda>L(\tau_p)$.

The reconstruction formula \eqref{eq:Xrecon} involves $W$ only
between $L$ and the fixed point $\lambda$, so the smooth joining
of $L$ and $W$ gives the corresponding smooth joining of $\X$.
The formula for $Q$ involves $W$ on $[\lambda,\infty)$. For every
fixed finite $\eta$,
$W(\tau,\eta)\to W(\tau_p,\eta)$ as
$\tau\downarrow\tau_p$, and $0\le W\le p$. Hence dominated
convergence in \eqref{eq:QFrecon} gives continuity of $Q$.
The identity $Q_\lambda=\gamma(Q-W)/\lambda$ then gives
continuity of $Q_\lambda$, while
$\lambda^2Q_{\lambda\lambda}
=(\gamma-1)\gamma(Q-W)-\gamma W_\zeta$
gives continuity of $Q_{\lambda\lambda}$. Finally,
$Q_\tau=\LL Q+1$ on both sides gives continuity of $Q_\tau$.
Thus $\X$, $Q$ and $F$ are continuous across $\tau=\tau_p$ together with the derivatives listed above, and $F$ is $C^{1,2}$ there.

We now pass to the primal variables. For every $z>0$, the relation
$z=\X(\tau,\lambda)$ determines a unique
$\lambda=\lambda(\tau,z)>L(\tau)$. On the no-action region side,
$\X_\lambda=(\X+W)/\lambda>0$, while on the purchase region
$\X_\lambda=d>0$. Differentiating
$z=\X(\tau,\lambda(\tau,z))$ gives
\[
    \lambda_z=\frac{1}{\X_\lambda},
    \qquad
    \lambda_\tau=-\frac{\X_\tau}{\X_\lambda}.
\]
Thus the strict positivity of $\X_\lambda$ allows the regularity
of $\X$ in $(\tau,\lambda)$ to pass to its inverse
$\lambda(\tau,z)$. At the purchase threshold, the matching
established above gives
$\X_\lambda(\tau,H-)=\X_\lambda(\tau,H+)=d(\tau)>0$, so the same
argument applies across $z=\zb(\tau)$. It also applies across
$\tau=\tau_p$ at every fixed $z>0$.

The formulas \eqref{eq:vderivs} follow from the conjugacy relations.
Continuity of $\X_\lambda$ at $H$ gives continuity of $v_{zz}$ at
$z=\zb(\tau)$, and $v_\tau=\lambda^{-\gamma}F_\tau$ is continuous.
Together with the matching of $F$ established above, this yields
$v\in C^{1,2}((0,\infty)\times(0,\infty))$, including across the
purchase threshold and across $\tau=\tau_p$. Homogeneity transfers this regularity to $V$ for $x,\xi>0$.
Reversing the linearization of \cref{sec:reduction}, the equations
\eqref{eq:linearPDEs} on $\Nr$ and the relations
\eqref{eq:actionF} and \eqref{eq:upperresidual} on $\Ar$ turn into
\eqref{eq:primalHJB} at every $z>0$.
\end{proof}

The three objects have different regularity properties. The dual value
$F$ is $C^{1,2}$ across both free boundaries. The inverse consumption map
$\X$ is only $C^1$ in the spatial variable across the thresholds, with
explicit jumps in its curvature. The primal value $v$ is $C^{1,2}$ at
every positive wealth level, including the purchase threshold, but its
curvature becomes singular at zero wealth.

Here and below, $L=L(\tau)$ denotes the limiting value
\[
L(\tau)=\lambda(\tau,0).
\]

\begin{lemma}\label{lem:zerowealth}
Locally uniformly for $\tau$ in compact subsets of $(0,\infty)$, as
$z\downarrow0$,
\begin{align}
 \lambda(\tau,z)&=L+L\sqrt{\frac{2a}{1-L}}\sqrt z+O(z),\label{eq:lambdaedge}\\
 v_z(\tau,z)&=L^{-\gamma}
 \Bigl[1-\gamma\sqrt{\frac{2a}{1-L}}\sqrt z+O(z)\Bigr],\label{eq:vzedge}\\
 v_{zz}(\tau,z)&=-\gamma L^{-\gamma}
 \sqrt{\frac{a}{2(1-L)}}\,z^{-1/2}+O(1).\label{eq:vzzedge}
\end{align}
\end{lemma}

\begin{proof}
By \eqref{eq:lowerall},
\[
\X(\tau,\lambda)
=
\frac{1-L}{2aL^2}(\lambda-L)^2
+O((\lambda-L)^3),
\]
where the coefficients depend smoothly on $\tau$. Inverting this
expansion gives \eqref{eq:lambdaedge}. The expansion
\eqref{eq:vzedge} follows from
$v_z=\lambda^{-\gamma}$. Finally, differentiating the inverse relation
yields
\[
v_{zz}
=
-\frac{\gamma\lambda^{-\gamma-1}}{\X_\lambda},
\]
and using
\[
\X_\lambda
=
\frac{1-L}{aL^2}(\lambda-L)
+O((\lambda-L)^2)
\]
gives \eqref{eq:vzzedge}.
\end{proof}

Hence $v_z(\tau,\cdot)$ is H\"older continuous of order $1/2$ at zero
wealth, and since the coefficient of $\sqrt z$ in \eqref{eq:vzedge}
is nonzero, this exponent is sharp. As $\tau\downarrow0$,
$L(\tau)^{-\gamma}\sim(K\tau)^{-\gamma}$, so the marginal utility at
zero wealth diverges. Estimates of the behavior near zero wealth
therefore cannot hold uniformly up to the retirement date.

\subsection{Integral equations for the free boundaries}\label{sec:integral}
The construction of \cref{sec:single,sec:hodo} determines $L$ and $H$, but it passes through two different formulations, the obstacle problem before $\tau_p$ and the hodograph equation after it. The free boundaries also satisfy integral equations that hold on the whole horizon. Put, for $u>0$ and $y\in\R$,
\begin{equation}\label{eq:kernels}
 \mathfrak g(u,y)=\frac{\ee^{-ru}}{\sqrt{4\pi au}}\exp\Bigl(-\frac{(y+b_0u)^2}{4au}\Bigr),\qquad
 \mathfrak G(u,y)=\int_{-\infty}^y\mathfrak g(u,y')\dd y' .
\end{equation}
Thus $\ee^{ru}\mathfrak G(u,y)$ is the standard normal distribution function evaluated at $(y+b_0u)/\sqrt{2au}$. The function $\mathfrak g$ is the fundamental solution of $\PP$ in the variables $(\tau,\zeta)$: it solves $\partial_u\mathfrak g=a\mathfrak g_{yy}+b_0\mathfrak g_y-r\mathfrak g$ and tends to the Dirac mass as $u\downarrow0$. Moreover $\mathfrak G(u,\zeta-c)=\int_c^\infty\mathfrak g(u,\zeta-\eta)\dd\eta$. We write $\zeta_L=\log L$ and $\zeta_H=\log H$, and we use the convention that terms containing $\zeta_H(s)$ are absent when $s\le\tau_p$.

\begin{theorem}\label{thm:integral}
Let \cref{ass:standing} hold and $rp<1$. For $\tau>0$ and $\zeta\in\R$,
\begin{multline}\label{eq:IErep}
 W(\tau,\ee^{\zeta})=\int_0^\tau\Bigl[\mathfrak G\bigl(\tau-s,\zeta-\zeta_L(s)\bigr)-(1-rp)\,\mathfrak G\bigl(\tau-s,\zeta-\zeta_H(s)\bigr)\\
 -\bigl(1-L(s)\bigr)\,\mathfrak g\bigl(\tau-s,\zeta-\zeta_L(s)\bigr)+\frac{1-rp}{\gamma}\,\mathfrak g\bigl(\tau-s,\zeta-\zeta_H(s)\bigr)\Bigr]\dd s .
\end{multline}
Denote the right-hand side by $\mathcal J(\tau,\zeta;L,H)$. The free boundaries solve the system of nonlinear Volterra equations
\begin{equation}\label{eq:IEsys}
 \mathcal J\bigl(\tau,\zeta_L(\tau);L,H\bigr)=0\quad(\tau>0),\qquad
 \mathcal J\bigl(\tau,\zeta_H(\tau);L,H\bigr)=p\quad(\tau>\tau_p).
\end{equation}
Moreover, for $\tau>\tau_p$,
\begin{equation}\label{eq:dformula}
 d(\tau)=\frac1K+\frac{(1-\gamma)(1-rp)}{\gamma}\int_{\tau_p}^{\tau}\frac{\ee^{-K(\tau-s)}}{H(s)}\dd s,\qquad \zb(\tau)=d(\tau)H(\tau)-p,
\end{equation}
and $\X$, $Q$, $F$ and $v$ are obtained from $W$ by \eqref{eq:Xrecon}, \eqref{eq:QFrecon} and \eqref{eq:vrecon}. For the frozen-income problem of \cref{sec:single} the same statements hold at every $\tau>0$ without the terms that contain $\zeta_H$: its Euler derivative is given by \eqref{eq:IErep} with these terms removed, and its free boundary solves the first equation in \eqref{eq:IEsys}.
\end{theorem}

\begin{proof}
\emph{The equation for $W$ on the whole line.} Work in the variables $(\tau,\zeta)$. By \cref{thm:single,prop:reconPDE,prop:regularity} and \cref{cor:bothboundaries}, $W$ is continuous on $(0,\infty)\times\R$, it satisfies $0\le W\le\varphi(\tau)$ for $\tau\le\tau_p$ and $0\le W\le p$ for $\tau>\tau_p$, and it is smooth up to the free boundaries from each side, and the boundaries are smooth curves. On $\Zr$, $\Nr$ and $\Ar$ we have $\PP W=0$, $1$ and $rp$, respectively, because $W$ equals $0$ on $\Zr$ and $p$ on $\Ar$. By \eqref{eq:lowercurv} and \eqref{eq:Wsboundaries}, $W_\zeta$ jumps by $(1-L)/a$ across $\zeta_L$ and by $-(1-rp)/(a\gamma)$ across $\zeta_H$. Let $\chi\in C^\infty_c((0,\infty)\times\R)$ and $\PP^*=-\partial_\tau-a\partial_{\zeta\zeta}+b_0\partial_\zeta+r$. On each region,
\[
 \chi\,\PP W-W\,\PP^*\chi=\partial_\tau(\chi W)-\partial_\zeta\bigl(a\chi W_\zeta-aW\chi_\zeta+b_0\chi W\bigr).
\]
Integrate over the three regions and apply the divergence theorem. Since $W$ is continuous across the boundaries, the flux terms that contain $W$ cancel, and only the jumps of $a\chi W_\zeta$ remain. As $H(\tau)\to\infty$ when $\tau\downarrow\tau_p$, the support of $\chi$ meets the upper boundary only over a compact subinterval of $(\tau_p,\infty)$. The result is
\[
 \iint W\,\PP^*\chi\dd\zeta\dd\tau=\iint\chi\,\bigl(\ind_{\Nr}+rp\ind_{\Ar}\bigr)\dd\zeta\dd\tau-\int_0^\infty(1-L)\,\chi(\tau,\zeta_L)\dd\tau+\frac{1-rp}{\gamma}\int_{\tau_p}^\infty\chi(\tau,\zeta_H)\dd\tau .
\]
Thus $\PP W=\mathfrak m$ in the sense of distributions, where $\mathfrak m$ is the measure on the right: the indicator functions of the regions and two single layers on the free boundaries.

\emph{Duhamel's formula.} Let $\widetilde W$ be the right-hand side of \eqref{eq:IErep}. By the properties of $\mathfrak G$ and $\mathfrak g$ listed before the theorem, $\widetilde W(\tau,\zeta)=\int_0^\tau\int\mathfrak g(\tau-s,\zeta-\eta)\,\mathfrak m(\mathrm ds,\mathrm d\eta)$ is the potential of $\mathfrak m$. Fix $\bar\tau>0$. The integrand in \eqref{eq:IErep} is bounded by $1+C(\tau-s)^{-1/2}$, so $\widetilde W$ is continuous and bounded on $(0,\bar\tau]\times\R$ with $|\widetilde W(\tau,\zeta)|\le C(\tau+\sqrt\tau)$. In particular, $\widetilde W(\tau,\cdot)\to0$ uniformly as $\tau\downarrow0$. Moreover $\PP\widetilde W=\mathfrak m$ in the sense of distributions. The difference $\Omega=W-\widetilde W$ therefore satisfies $\PP\Omega=0$ in the sense of distributions and is a classical solution. It is bounded on $(0,\bar\tau]\times\R$, and $\Omega(\tau,\cdot)\to0$ uniformly as $\tau\downarrow0$ because $0\le W\le\varphi(\tau)$ near $\tau=0$. The function $\ee^{r\tau}\Omega$ solves a heat equation with constant drift, and bounded solutions of the Cauchy problem with zero initial values vanish. Hence $\Omega\equiv0$ on $(0,\bar\tau]\times\R$ for every $\bar\tau$, which is \eqref{eq:IErep}. Evaluating \eqref{eq:IErep} at $\zeta_L(\tau)$, where $W=0$, and at $\zeta_H(\tau)$, where $W=p$, gives \eqref{eq:IEsys}. For the frozen-income problem the same argument applies with $\Ar=\emptyset$, using \cref{thm:single} and \cref{lem:singleW} for every $\tau>0$.

\emph{The coefficient $d$.} Formula \eqref{eq:dformula} is the solution of the linear equation \eqref{eq:dODE} with the initial value $d(\tau_p+)=1/K$ from \eqref{eq:escape}, and $\zb=dH-p$ is \eqref{eq:dB}.
\end{proof}

\begin{remark}\label{rem:integral}
(i) The two layers in \eqref{eq:IErep} carry the second order free boundary conditions: their densities are minus the jumps of $aW_\zeta$, that is, of $a\lambda^2\X_{\lambda\lambda}$, across the boundaries, which are given by \eqref{eq:lowerall} and \eqref{eq:uppercurv}. The system \eqref{eq:IEsys} then imposes the values of $W$ on the boundaries. This is the counterpart, for two boundaries, of the early exercise premium representation of the American put \citep{Kim1990,CarrJarrowMyneni1992}. The representation is written for the Euler derivative rather than for the dual value function. This has two advantages: the source term is piecewise constant, so that only the elementary kernels \eqref{eq:kernels} appear, and each boundary condition is imposed on a function that leaves its boundary value with nonzero slope rather than tangentially, which suits time stepping.

(ii) Formula \eqref{eq:dformula} gives the sign of $d-1/K$ directly: $d<1/K$ for $\gamma>1$ and $d>1/K$ for $\gamma<1$, in agreement with \cref{thm:horizon}(b).

(iii) The system \eqref{eq:IEsys} does not involve the retirement payoff directly. If the initial datum $\lambda/K$ of $\X$ is replaced by $c\lambda$ with another constant $c>0$, while the operator $\PP$ is kept, then still $W(0,\cdot)=0$ and the form of the equations is unchanged. The payoff enters through the normalization $L(\tau)/(K\tau)\to1$ as $\tau\downarrow0$ from \cref{thm:single}, and through the identity $\int_{L(\tau)}^\infty W(\tau,\eta)\,\eta^{-2}\dd\eta=1/K$ for $\tau\le\tau_p$, which follows from \eqref{eq:Xrecon} and \cref{thm:single}(a). For $\tau>\tau_p$ the same integral equals $d(\tau)$. In the computations the normalization is imposed through the values of $L$ near $\tau=0$. We do not study uniqueness of solutions of \eqref{eq:IEsys} under this normalization. For the single boundary of the American put, uniqueness of the corresponding equation was proved by \citet{Peskir2005}. \Cref{sec:single,sec:hodo} already determine $L$ and $H$, and in \cref{sec:numerics} the system serves as a numerical method.
\end{remark}

\section{Optimal strategy and horizon effects}\label{sec:verification}

This section constructs the optimal strategy, proves its optimality and the results of \cref{sec:results}, and then studies how the solution depends on the time to retirement.

\subsection{The candidate strategy}
Fix $(t,x,\xi)$ and put $\tau=T-t$. If $\tau>\tau_p$ and $x/\xi>\zb(\tau)$, the investor makes the initial purchase
\begin{equation}\label{eq:initialpurchase}
 \Delta\Xi_t=\frac{x-\zb(\tau)\xi}{p+\zb(\tau)},\qquad \Xi_t=\frac{x+p\xi}{p+\zb(\tau)},\qquad X_t=\zb(\tau)\,\Xi_t ,
\end{equation}
otherwise $\Xi_t=\xi$ and $X_t=x$. After this, let $\Lambda_t\in[L(\tau),H(\tau)]$ be the unique solution of $X_t/\Xi_t=\X(\tau,\Lambda_t)$, with $\Lambda_t=L(\tau)$ if $X_t=0$. Let
\begin{equation}\label{eq:logreflection}
 \mathrm d\log\Lambda_s=\nu\dd s+\frac\theta\gamma\dd B_s+\mathrm d\Gamma^L_s-\mathrm d\Gamma^H_s,\qquad \nu=\frac{r-\beta}{\gamma}+\gamma a,
\end{equation}
where $\Gamma^L,\Gamma^H$ are continuous, nondecreasing, start at zero, and
\[
 L(T-s)\le\Lambda_s\le H(T-s),\qquad \int\ind_{\{\Lambda_s>L(T-s)\}}\dd\Gamma^L_s=0,\qquad \int\ind_{\{\Lambda_s<H(T-s)\}}\dd\Gamma^H_s=0 .
\]
The upper constraint and $\Gamma^H$ are absent when $H(T-s)=\infty$. The controls are
\begin{equation}\label{eq:controls}
 \mathrm d\log\Xi^*_s=\mathrm d\Gamma^H_s,\qquad c^*_s=\Xi^*_s\Lambda_s,\qquad
 \pi^*_s=\frac{\theta}{\gamma\sigma}\,\Xi^*_s\Lambda_s\X_\lambda(T-s,\Lambda_s),\qquad X^*_s=\Xi^*_s\X(T-s,\Lambda_s).
\end{equation}
Since $\lambda\X_\lambda=\X+W$, the portfolio in \eqref{eq:controls} is the one in \eqref{eq:policies}. The drift of $\Lambda$ is $\nu+a=b+\theta^2/\gamma$.

\begin{remark}\label{rem:normalization}
If the reflection is written with additive regulators, $\mathrm d\Lambda_s=(\nu+a)\Lambda_s\dd s+\frac\theta\gamma\Lambda_s\dd B_s+\mathrm d\ell^L_s-\mathrm d\ell^H_s$, then $\mathrm d\ell^H=\Lambda\dd\Gamma^H$ and the income rule becomes
\[
 \mathrm d\log\Xi^*_s=\frac{\mathrm d\ell^H_s}{H(T-s)}=\frac{\X_\lambda}{\X+p}\bigl(T-s,H(T-s)\bigr)\,\mathrm d\ell^H_s .
\]
Indeed, if $\mathrm d\Xi/\Xi=\kappa\,\mathrm d\ell^H$, the singular part of $\mathrm d(\Xi\X)$ at the upper boundary is $\Xi(\kappa\X-\X_\lambda)\,\mathrm d\ell^H$, and the budget equation requires it to equal $-p\Xi\kappa\,\mathrm d\ell^H$, which forces $\kappa=\X_\lambda/(\X+p)=1/H$.
\end{remark}

\begin{lemma}\label{lem:admissible}
The system \eqref{eq:logreflection}--\eqref{eq:controls} has a unique strong solution on $[t,T]$. The strategy $(c^*,\pi^*,\Xi^*)$ is admissible, $X^*$ solves \eqref{eq:wealth}, and $X^*_T=c^*_T/K>0$. The processes $\Sigma_s=X^*_s+p\Xi^*_s$ and $c^*_s$ have finite moments $\E[\sup_{t\le s\le T}\Sigma_s^k]$ and $\E[\sup_{t\le s\le T}(c^*_s)^k]$ for every $k\in\R$.
\end{lemma}

\begin{proof}
On compact time intervals on which both boundaries are finite, their logarithms are smooth and strictly separated by \cref{cor:bothboundaries}, and the Skorokhod problem in a time dependent interval has a unique pathwise solution \citep{BurdzyKangRamanan2009}. Where $H=\infty$ the one-sided Skorokhod map is used. We first construct the solution on $[t,T-\tau_p-\epsilon]$ and on $[T-\tau_p,T-\epsilon]$ for small $\epsilon>0$, and then remove the localization using the bounds below, which do not depend on $\epsilon$.

Apply It\^o's formula to $\Xi^*_s\X(T-s,\Lambda_s)$. By \cref{prop:regularity}, $\X$ is $C^1$ in $\lambda$ across the free boundaries and smooth up to them from $\Nr$, so no local time correction appears. On $\Nr$, $\PP\X=\lambda-1$ and $\nu+a-b=\theta^2/\gamma$ give the drift $rX^*+\sigma\theta\pi^*-c^*+\Xi^*$ and the diffusion coefficient $\sigma\pi^*$. At the lower boundary $\X_\lambda=0$, so $\mathrm d\Gamma^L$ does not enter. At the upper boundary the singular terms are $\Xi^*(\X-\Lambda\X_\lambda)\,\mathrm d\Gamma^H=-p\,\Xi^*\,\mathrm d\Gamma^H=-p\,\mathrm d\Xi^*$. Hence $X^*$ solves \eqref{eq:wealth}, and $X^*\ge0$ by construction.

For the moments, $0\le W\le p$ and $\lambda\X_\lambda=\X+W$ give
\begin{equation}\label{eq:momentbounds}
 \Xi^*_s\le \Sigma_s/p,\qquad |\pi^*_s|\le\frac{|\theta|}{\gamma\sigma}\Sigma_s .
\end{equation}
Moreover $\partial_\lambda\bigl((\X+p)/\lambda\bigr)=(W-p)/\lambda^2\le0$ and $(\X+p)/\lambda\to d(\tau)$ as $\lambda\to\infty$, so $c^*_s/\Sigma_s\le1/d(T-s)\le C_T:=1/\min_{[0,T]}d<\infty$. Purchases cancel in the dynamics of $\Sigma$:
\begin{equation}\label{eq:Sdynamics}
 \mathrm d\Sigma_s=(rX^*_s+\sigma\theta\pi^*_s+\Xi^*_s-c^*_s)\dd s+\sigma\pi^*_s\dd B_s .
\end{equation}
Drift and volatility of $\Sigma$ divided by $\Sigma$ are bounded by deterministic constants, so $\Sigma$ is a stochastic exponential with bounded coefficients. It has finite positive and negative moments of all orders, uniformly on $[t,T]$, and a finite positive limit at the end of each localization interval.

The upper regulator cancels from consumption:
\begin{equation}\label{eq:logc}
 \mathrm d\log c^*_s=\nu\dd s+\frac\theta\gamma\dd B_s+\mathrm d\Gamma^L_s .
\end{equation}
Thus $c^*$ dominates the geometric Brownian motion obtained by dropping $\Gamma^L$, which gives all negative moments, and $c^*\le C_T\Sigma$ gives all positive moments. Since $\Xi_t\le\Xi^*\le\Sigma/p$, $\Gamma^H$ stays bounded, and then \eqref{eq:logc} keeps $\Gamma^L$ bounded. Consequently $\Lambda=c^*/\Xi^*$ stays, path by path, in a compact subset of $(0,\infty)$ that does not depend on $\epsilon$. Since $H(T-s)\to\infty$ as $s\uparrow T-\tau_p$, the upper reflection is inactive on a neighborhood of $T-\tau_p$, so $\Gamma^H$ is eventually constant there and the two pieces join into one solution. Likewise, $L(T-s)\to0$ as $s\uparrow T$, so the lower reflection is inactive near $T$, and the solution extends continuously to $s=T$. Finally $|\X(\tau,\lambda)-\lambda/K|\le\varphi(\tau)$ on $(0,\tau_p]$ and income is constant there, hence $X^*_T=\lim_{s\uparrow T}\Xi^*_s\X(T-s,\Lambda_s)=\Xi^*_T\Lambda_T/K=c^*_T/K$.
\end{proof}

\subsection{Shadow prices}
Along the candidate strategy define
\begin{equation}\label{eq:prices}
 Y_s=(c^*_s)^{-\gamma},\qquad Z_s=Y_s\,Q(T-s,\Lambda_s) .
\end{equation}
$Y$ is the marginal utility of wealth and $Z$ that of income. We have $Y>0$ and $0\le Z\le pY$, and $Z_s\to0$ as $s\uparrow T$ because $Q\le\varphi(T-s)$ on the final interval.

\begin{lemma}\label{lem:adjoint}
The shadow prices satisfy
\begin{align}
 \mathrm dY_s&=(\beta-r)Y_s\dd s-\theta Y_s\dd B_s-\gamma Y_s\dd\Gamma^L_s,\label{eq:Yadjoint}\\
 \mathrm dZ_s&=(\beta Z_s-Y_s)\dd s-\theta Y_sW(T-s,\Lambda_s)\dd B_s .\label{eq:Zadjoint}
\end{align}
In particular the income price has no singular component.
\end{lemma}

\begin{proof}
Equation \eqref{eq:Yadjoint} follows from \eqref{eq:logc}, since $-\gamma\nu+\theta^2/2=\beta-r$. For $Z$, It\^o's formula and $\PP Q=1$ on $\Nr$ give a drift in which the correction $(\nu+a-b)DQ=\frac{\theta^2}{\gamma}DQ$ cancels against the covariation $-\theta Y\cdot\frac\theta\gamma DQ$, leaving $\beta Z-Y$. The singular terms are $Y(DQ-\gamma Q)\,\mathrm d\Gamma^L-Y\,DQ\,\mathrm d\Gamma^H$. The first vanishes because $DQ=\gamma(Q-W)=\gamma Q$ at $L$, and the second because $DQ=0$ at $H$. The Brownian coefficient is $Y(\frac\theta\gamma DQ-\theta Q)=Y(\theta(Q-W)-\theta Q)=-\theta YW$.
\end{proof}

\subsection{Verification}
\begin{theorem}\label{thm:verification}
The strategy of \cref{lem:admissible}, preceded by \eqref{eq:initialpurchase} when applicable, is optimal, and $V(t,x,\xi)=\xi^{1-\gamma}v(T-t,x/\xi)$ with $v$ given by \eqref{eq:vrecon}. The value is finite for every $(t,x,\xi)$, and the conclusion holds for $0<\gamma<1$ and $\gamma>1$.
\end{theorem}

\begin{proof}
Keep the shadow prices $(Y,Z)$ of the candidate and let $(c,\pi,\Xi,X)$ be an arbitrary admissible strategy. By It\^o's product formula and \cref{lem:adjoint}, the process
\begin{equation}\label{eq:budgetprocess}
 \mathcal M_s=\ee^{-\beta(s-t)}\bigl(Y_sX_s+Z_s\Xi_s\bigr)+\int_t^s\ee^{-\beta(u-t)}Y_uc_u\dd u
\end{equation}
has the decomposition $\mathrm d\mathcal M_s=\ee^{-\beta(s-t)}\bigl[(Z_s-pY_s)\,\mathrm d\Xi_s-\gamma Y_sX_s\,\mathrm d\Gamma^L_s\bigr]+\mathrm d(\text{local martingale})$: the absolutely continuous drift cancels exactly. The remaining finite variation part is nonpositive because $Z\le pY$, $X\ge0$ and $\Xi$, $\Gamma^L$ are nondecreasing. All terms of $\mathcal M$ are nonnegative, so $\mathcal M$ is a nonnegative local supermartingale and hence a supermartingale. Including a possible purchase at $t$, which changes $Y_tx+Z_t\xi$ by $(Z_t-pY_t)\Delta\Xi_t\le0$, and using $Z_T=0$, we obtain the budget inequality
\begin{equation}\label{eq:budget}
 \E\Bigl[\int_t^T\ee^{-\beta(s-t)}Y_sc_s\dd s+\ee^{-\beta(T-t)}Y_TX_T\Bigr]\le Y_tx+Z_t\xi .
\end{equation}
No integrability of the competing strategy is needed beyond admissibility.

For $y>0$ let $\widetilde u(y)=\sup_{c\ge0}\{u(c)-yc\}=\frac{\gamma}{1-\gamma}y^{(\gamma-1)/\gamma}$ and $\widetilde R(y)=\sup_{x\ge0}\{R(x)-yx\}=\frac{\gamma}{(1-\gamma)K}y^{(\gamma-1)/\gamma}$. Then $u(c_s)\le\widetilde u(Y_s)+Y_sc_s$ and $R(X_T)\le\widetilde R(Y_T)+Y_TX_T$ pathwise, and $\widetilde u(Y)$ and $\widetilde R(Y_T)$ are integrable by \cref{lem:admissible}, because $Y^{(\gamma-1)/\gamma}=(c^*)^{1-\gamma}$. Combining with \eqref{eq:budget},
\begin{multline}\label{eq:dualbound}
 \E\Bigl[\int_t^T\ee^{-\beta(s-t)}u(c_s)\dd s+\ee^{-\beta(T-t)}R(X_T)\Bigr]\\
 \le\E\Bigl[\int_t^T\ee^{-\beta(s-t)}\widetilde u(Y_s)\dd s+\ee^{-\beta(T-t)}\widetilde R(Y_T)\Bigr]+Y_tx+Z_t\xi<\infty .
\end{multline}
This argument uses only pathwise inequalities and the supermartingale property, so it is valid for either sign of $1-\gamma$.

For the candidate both singular terms in $\mathrm d\mathcal M$ vanish: $X^*=0$ on the support of $\mathrm d\Gamma^L$ and $Z=pY$ on the support of $\mathrm d\Xi^*$. The initial purchase, when it occurs, is made at a state where $Z_t=pY_t$ and therefore has no cost in \eqref{eq:budget}. The candidate's $\mathcal M$ is a true martingale, since $\sup_s\mathcal M_s\le\sup_s(c^*_s)^{-\gamma}\Sigma_s+\int_t^T(c^*_u)^{1-\gamma}\dd u$ is integrable by \cref{lem:admissible}. Hence \eqref{eq:budget} holds with equality. The Fenchel inequalities are also equalities, because $c^*=Y^{-1/\gamma}$ and $X^*_T=c^*_T/K=Y_T^{-1/\gamma}/K$. The candidate therefore attains the upper bound \eqref{eq:dualbound} and is optimal.

It remains to compute the optimal value. Equality in the budget and in the Fenchel inequalities gives $u(c^*)=Y c^*/(1-\gamma)$ and $R(X^*_T)=Y_TX^*_T/(1-\gamma)$, so the optimal objective equals $(Y_tx+Z_t\xi)/(1-\gamma)$. Without an initial purchase, $Y_t=(\xi\Lambda_t)^{-\gamma}$ and $Z_t=Y_tQ(\tau,\Lambda_t)$, and this is $\xi^{1-\gamma}v(\tau,x/\xi)$ by \eqref{eq:vrecon}. With an initial purchase, the extended inverse is $\lambda=(z+p)/d(\tau)$. After the purchase $\Lambda_t=H(\tau)$ and $\Xi_t=(x+p\xi)/(d(\tau)H(\tau))$, so $c^*_t=(x+p\xi)/d(\tau)=\xi\lambda$ and $Y_t$ is unchanged, while $Z_t=pY_t$. The same identity follows. At zero wealth $\Lambda_t=L(\tau)$.
\end{proof}

\subsection{Proofs of the main results}\label{sec:proofsmain}
\begin{proof}[Proof of \cref{thm:main,thm:optimal}]
The functions $L$, $H$, $\zb$ and $d$ are constructed in \cref{thm:single}, \cref{cor:bothboundaries} and \eqref{eq:dB}, and parts (i) and (ii) of \cref{thm:main} are contained in these results and in \eqref{eq:escape}, except for the bound \eqref{eq:blowup}, which is proved in \cref{sec:comparison}. \Cref{thm:verification} gives $V=\xi^{1-\gamma}v$ with $v$ from \eqref{eq:vrecon}, \cref{prop:regularity} gives part (iii), and \cref{lem:zerowealth} gives part (iv). The purchase region is empty on $(0,\tau_p]$ by \cref{thm:window}, and nonempty for $\tau>\tau_p$ because $\zb(\tau)<\infty$ there. If $rp\ge1$, then $\tau_p=\infty$, \cref{thm:single} applies for every $\tau$, and the reconstruction of \cref{sec:reconstruction} with $H\equiv\infty$ satisfies $Q<\varphi(\tau)<1/r\le p$, and the arguments of \cref{sec:reconstruction,sec:verification} then apply verbatim without the upper boundary. \Cref{thm:optimal} is \cref{thm:verification} together with \cref{lem:admissible}.
\end{proof}

\begin{proof}[Proof of \cref{prop:properties}]
(a) On the no-action region $\lambda(\tau,\cdot)$ is the inverse of the strictly increasing, strictly convex function $\X(\tau,\cdot)$, hence strictly increasing and strictly concave with $\lambda_z=1/\X_\lambda$. On the purchase region $\lambda=(z+p)/d$ is affine, and $\X_\lambda(\tau,H)=d$ makes $\lambda_z$ continuous at $\zb$. Thus $\lambda_z$ is nonincreasing on $(0,\infty)$, and $\lambda(\tau,\cdot)$, being continuous at $z=0$, is concave on $[0,\infty)$. The expansion near $z=0$ follows from \cref{lem:zerowealth}.
(b) The formula is \eqref{eq:controls}. The bounds on $W$ come from \cref{lem:singleW} and \eqref{eq:Wextension}. Near zero wealth, $W=\lambda\X_\lambda-\X=\sqrt{2(1-L)z/a}+O(z)$ by \eqref{eq:lowerall} and \eqref{eq:lambdaedge}, and $\frac{\theta}{\gamma\sigma}\sqrt{2/a}=\frac{2\sgn\theta}{\sigma}$.
(c) This was shown after \eqref{eq:QminusW}, and on $\Zr$, $W=0\le Q$.
(d) This is $X^*_T=c^*_T/K$ from \cref{lem:admissible}, combined with $c^R=KX$.
\end{proof}

\subsection{Monotonicity in the time to retirement}\label{sec:horizon}
\begin{proof}[Proof of \cref{thm:horizon}]
(a) Fix $z\ge0$ and $0<\tau_1<\tau_2$. By homogeneity and time invariance we may start at time $0$ with wealth $z$ and income $1$ and compare the problems with retirement dates $\tau_1$ and $\tau_2$. Let $(c,\pi,\Xi)$ be admissible for the first problem, with wealth $X$ and a finite objective $J_1$. When $\gamma>1$ this forces $X_{\tau_1}>0$ almost surely. Let $X^M$ be the Merton wealth on $[\tau_1,\infty)$ started from $X_{\tau_1}$ under the feedbacks \eqref{eq:retirement}. It is a geometric Brownian motion, and the process $\ee^{-\beta s}R(X^M_s)+\int_{\tau_1}^s\ee^{-\beta q}u(KX^M_q)\dd q$ is a martingale on $[\tau_1,\infty)$, so
\begin{equation}\label{eq:mertonsplit}
 \E\bigl[\ee^{-\beta\tau_1}R(X_{\tau_1})\bigr]=\E\Bigl[\int_{\tau_1}^{\tau_2}\ee^{-\beta s}u(KX^M_s)\dd s+\ee^{-\beta\tau_2}R(X^M_{\tau_2})\Bigr].
\end{equation}
In the second problem, follow $(c,\pi,\Xi)$ on $[0,\tau_1]$, and on $(\tau_1,\tau_2]$ keep income at $\Xi_{\tau_1}$, consume $KX^M_s$, hold $\theta X^M_s/(\gamma\sigma)$ in the stock, and invest the income in the bond. Wealth is then $X^M_s+\Xi_{\tau_1}(\ee^{r(s-\tau_1)}-1)/r\ge X^M_s\ge0$, so the strategy is admissible. It consumes the same amounts as the first strategy followed by the Merton policy, and its terminal wealth is at least $X^M_{\tau_2}$. Since $R$ is nondecreasing, \eqref{eq:mertonsplit} shows that its objective is at least $J_1$. Taking suprema gives $v(\tau_2,z)\ge v(\tau_1,z)$.

(b) For $\tau_p<\tau_1<\tau_2$ and $z\ge\max\{\zb(\tau_1),\zb(\tau_2)\}$, part (a) and \cref{thm:main}(iii) give $d(\tau_1)^\gamma/(1-\gamma)\le d(\tau_2)^\gamma/(1-\gamma)$. Multiplying by $1-\gamma$ shows that $d$ is nonincreasing when $\gamma>1$ and nondecreasing when $\gamma<1$. Letting $\tau_1\downarrow\tau_p$ and using \eqref{eq:escape} gives the bounds $d\le1/K$ and $d\ge1/K$, respectively. The statement on consumption follows from $\lambda=(z+p)/d(\tau)$ on the purchase region.

(c) Let $\gamma>1$. By (b) and the smoothness of $d$, $d'\le0$ on $(\tau_p,\infty)$. Fix $\tau_1>\tau_p$ and let $\mathcal D=\{(\tau,\lambda):\tau_p<\tau\le\tau_1,\ L(\tau)<\lambda<H(\tau)\}$. The function $\Psi=-\X_\tau$ is smooth in $\mathcal D$, continuous up to its lateral boundaries (\cref{prop:regularity}), and $\PP\Psi=0$ in $\mathcal D$. On $\lambda=L(\tau)$, $\Psi=0$ by \cref{prop:reconPDE}. On $\lambda=H(\tau)$, differentiating $\X(\tau,H(\tau))=d(\tau)H(\tau)-p$ and using $\X_\lambda(\tau,H)=d$ gives $\Psi=-d'H\ge0$. As $\tau\downarrow\tau_p$, $\Psi(\tau,\lambda)\to-\X_\tau(\tau_p,\lambda)\ge0$ for each fixed $\lambda$, by \cref{thm:single}(a) and the continuity of $\X_\tau$ across $\tau_p$ (\cref{prop:regularity}). More precisely, $\Psi$ extends continuously to $\overline{\mathcal D}\cap\{\lambda\le M\}$ for every $M$, including the corner $(\tau_p,L(\tau_p))$, because in hodograph variables $\X_\tau=\lambda\int_L^\lambda W_\tau\eta^{-2}\dd\eta$ with $W_\tau=-\rho_h/\rho_w$, and the glued hodograph is smooth across $\tau_p$ up to $w=0$ (\cref{cor:bothboundaries}). Finally, since $\lambda\X_\lambda=\X+W$ and $\lambda^2\X_{\lambda\lambda}=W_\zeta$,
\[
 \Psi=-aW_\zeta-b(\X+W)+r\X-\lambda+1\qquad\text{on }\mathcal D,
\]
and $W_\zeta=1/\rho_w\le1/\underline q$, $0\le W\le p$ and $0\le\X\le p\lambda/\min_{[\tau_p,\tau_1]}L$ give $|\Psi|\le C_1(1+\lambda)$ on $\mathcal D$.

For $\epsilon>0$ put $\Psi_\epsilon=\Psi+\epsilon\lambda^2\ee^{C_2(\tau-\tau_p)}$ with $C_2=|2a+2b-r|+1$. Then $\PP\Psi_\epsilon=\epsilon(C_2-2a-2b+r)\lambda^2\ee^{C_2(\tau-\tau_p)}>0$ in $\mathcal D$. Choose $M_\epsilon$ with $\epsilon M_\epsilon^2>C_1(1+M_\epsilon)$ and apply the minimum principle on $\mathcal D\cap\{\lambda<M_\epsilon\}$: its parabolic boundary consists of parts of the curves $\lambda=L(\tau)$ and $\lambda=H(\tau)$, where $\Psi_\epsilon>0$, of the segment $\lambda=M_\epsilon$, where $\Psi_\epsilon>0$ by the growth bound, and of the initial segment $\tau=\tau_p$, where $\liminf\Psi_\epsilon\ge\epsilon\lambda^2>0$. Hence $\Psi_\epsilon\ge0$ on $\mathcal D$, and letting $\epsilon\downarrow0$ gives $\Psi\ge0$ on $\mathcal D$.

Outside $\bigcup_{\tau_1}\mathcal D$ the sign of $\X_\tau$ is known: $\X_\tau\le0$ for $\tau\le\tau_p$ by \cref{thm:single}(a), $\X_\tau=0$ at points with $\lambda<L(\tau)$, where $\X$ vanishes in a neighborhood, and $\X_\tau=d'\lambda\le0$ for $\lambda\ge H(\tau)$ by \eqref{eq:actionF}. Thus $\X_\tau\le0$ everywhere. Since $\X$ is $C^1$ in $\tau$, $\X(\tau_2,\lambda)\le\X(\tau_1,\lambda)$ for $\tau_1<\tau_2$, and taking $\lambda\le L(\tau_1)$ gives $\X(\tau_2,\lambda)=0$, that is, $L(\tau_2)\ge L(\tau_1)$. Finally $\lambda(\tau,\cdot)$ is the inverse of the increasing function $\X(\tau,\cdot)$ on $[L(\tau),\infty)$, extended by $L(\tau)$ at $z=0$, so $\X_\tau\le0$ implies that $\lambda(\tau,z)$ is nondecreasing in $\tau$.
\end{proof}

\subsection{Comparison with frozen income}\label{sec:comparison}
\begin{proof}[Proof of \cref{thm:compare}]
For $\tau>\tau_p$ let $\rho^0(\tau-\tau_p,\cdot)$ be the hodograph of the frozen-income solution, $W^0(\tau,\ee^{\rho^0(\tau-\tau_p,w)})=w$ for $0\le w\le p$. This is well defined because $W^0(\tau,\cdot)$ increases from $0$ to $\varphi(\tau)>p$. By \cref{thm:single} and \cref{lem:singleW}, $\rho^0$ is smooth on $(0,\infty)\times[0,p]$, solves \eqref{eq:hodoPDE} with the lower boundary condition in \eqref{eq:hodoBC}, and $\rho^0(h,\cdot)\to\rho_0$ locally uniformly on $[0,p)$ as $h\downarrow0$. Its slope at $w=p$ is $\rho^0_w(h,p)=1/(\lambda W^0_\lambda)(\tau,\Lambda^0_p(\tau))$.

We first show $h_1>0$. Since $W^0(\tau_p,\lambda)<p$ for every $\lambda$ and $W^0(\tau,\cdot)\to W^0(\tau_p,\cdot)$ locally uniformly, $\Lambda^0_p(\tau)\to\infty$ as $\tau\downarrow\tau_p$. By \cref{lem:singleW}, $\lambda W^0_\lambda(\tau,\lambda)\to0$ as $\lambda\to\infty$, uniformly for $\tau$ near $\tau_p$. Hence $\lambda W^0_\lambda(\tau,\Lambda^0_p(\tau))\to0<(1-rp)/(a\gamma)$, and $h_1>0$.

For $0<h\le h_1$ the slope condition means $\rho^0_w(h,p)\ge c_H$, so $\rho^0$ is a supersolution of the two-boundary problem on $(0,h_1]\times[0,p]$. The truncated data satisfy $\rho_{0,n}\le\rho_0=\rho^0(0,\cdot)$, the solutions $\rho_n$ of \cref{lem:truncated} are bounded, and $\rho^0(h,w)\to\infty$ as $(h,w)\to(0,p)$, because $\rho^0(h,w)\ge\rho^0(h,p-\delta)\to\rho_0(p-\delta)$ for $w\ge p-\delta$ and $\rho_0(p-\delta)\to\infty$ as $\delta\downarrow0$. Hence $(\rho_n-\rho^0)_+$ vanishes near the corner and tends to zero in $L^2$ as $h\downarrow0$. The slopes of both functions are bounded on $[0,p/2]$. For $\rho^0$ this holds because $\lambda W^0_\lambda$ is bounded below on $\{W^0\le p/2\}$, uniformly for $\tau$ near $\tau_p$. \Cref{lem:comparison}, applied on $(0,h]$ for every finite $h\le h_1$, gives $\rho_n\le\rho^0$, and letting $n\to\infty$ gives $\rho\le\rho^0$ on $(0,h_1]\times[0,p]$. At $w=0$ this is $L\le L^0$, and at $w=p$ it is $H\le\Lambda^0_p$. For $L(\tau)\le\lambda\le H(\tau)$ write $\lambda=\ee^{\rho(h,w)}$ with $w=W(\tau,\lambda)$. Then $\lambda\le\ee^{\rho^0(h,w)}$, and since $W^0(\tau,\cdot)$ is nondecreasing, $W^0(\tau,\lambda)\le w=W(\tau,\lambda)$.

For the bound on $\Lambda^0_p$ take $\zeta_0=0$ in the first passage representation \eqref{eq:tailkernel}, which is admissible because $L^0\le L^0_\infty<1$. Since $g_0\le\varphi(\tau)$ and $\ee^{-rv}\le1$,
\[
 \varphi(\tau)-W^0(\tau,\ee^\zeta)\le\varphi(\tau)\,\mathbb P\Bigl(\min_{v\le\tau}\bigl(\zeta+b_0v+\sqrt{2a}\,\widetilde B_v\bigr)\le0\Bigr)\le\varphi(\tau)\exp\Bigl(-\frac{(\zeta-|b_0|\tau)^2}{4a\tau}\Bigr)
\]
for $\zeta\ge|b_0|\tau$, by the reflection principle and the Gaussian tail bound (a standard normal variable exceeds $x\ge0$ in absolute value with probability at most $\ee^{-x^2/2}$). At $\zeta=\log\Lambda^0_p(\tau)$ the left side equals $\varphi(\tau)-p=\int_{\tau_p}^\tau\ee^{-rs}\dd s\ge\ee^{-r\tau}(\tau-\tau_p)$. Solving for $\zeta$ gives the stated bound.
\end{proof}

\begin{proof}[Proof of the bound \eqref{eq:blowup} in \cref{thm:main}]
Take $h_*=\min\{h_1,1\}$. On $(\tau_p,\tau_p+h_*]$, \cref{thm:compare} gives $\log H\le\log\Lambda^0_p$, and $\varphi(\tau)\ee^{r\tau}$ is bounded there. Moreover $\zb=dH-p\le dH$, and $d$ is bounded near $\tau_p$ because it is continuous with $d(\tau_p+)=1/K$.
\end{proof}

\subsection{The stationary problem}
The stationary free boundary problem is the time independent version of \eqref{eq:hodoPDE}--\eqref{eq:hodoBC}: find $\rho_s\in C^2([0,p])$ with $\rho_s'>0$ and $\rho_s(0)<0$ such that
\begin{equation}\label{eq:stationaryproblem}
 \frac{a\rho_s''}{(\rho_s')^2}-b_0-(1-rw)\rho_s'=0\quad\text{on }[0,p],\qquad \rho_s'(0)=\frac{a}{1-\ee^{\rho_s(0)}},\qquad \rho_s'(p)=c_H .
\end{equation}
Undoing the hodograph transformation, it is equivalent to finding $0<L^\infty<H^\infty$ and an increasing $W^\infty$ with $-\LL W^\infty=1$ on $(L^\infty,H^\infty)$, $W^\infty=0$ and $\lambda W^\infty_\lambda=(1-L^\infty)/a$ at $L^\infty$, and $W^\infty=p$ and $\lambda W^\infty_\lambda=(1-rp)/(a\gamma)$ at $H^\infty$. These are the free boundary conditions of the perpetual income problem, written in the variables of this paper.

\begin{lemma}\label{lem:stationary}
\textup{(i)} The terminal value problem \eqref{eq:ustar} has a unique solution $u_*$ on $[0,p]$, and $u_*>0$.
\textup{(ii)} Problem \eqref{eq:stationaryproblem} has a solution if and only if $u_*(0)<1/a$. The solution is then unique and given by $\rho_s'=1/u_*$ and $\ee^{\rho_s(0)}=1-a\,u_*(0)$.
\textup{(iii)} The function $u^0(w)=\varsigma(1/r-w)$ solves the differential equation in \eqref{eq:ustar} on $[0,p]$, and $u^0(0)=(1-L^0_\infty)/a$. Moreover $a\gamma\varsigma<r$, so that $u_*>u^0$ on $[0,p]$.
\textup{(iv)} If $\gamma\ge1$, then $u_*(w)<(1-rw)/a$ for $0\le w<p$. In particular, $u_*(0)<1/a$.
\end{lemma}

\begin{proof}
(i) The right side $f_*(w,u)=-(b_0+(1-rw)/u)/a$ is smooth on $[0,p]\times(0,\infty)$, and we solve backward from $w=p$. If $b_0\ge0$, then $f_*<0$, so $u$ increases as $w$ decreases and $u\ge u(p)$. If $b_0<0$, then $f_*<0$ whenever $u<(1-rp)/|b_0|$, so $u$ cannot fall below $\min\{u(p),(1-rp)/|b_0|\}$ as $w$ decreases. In both cases $u$ stays bounded away from zero, $|f_*|$ stays bounded, and the solution extends to all of $[0,p]$.
(ii) If $\rho_s$ solves \eqref{eq:stationaryproblem}, then $u=1/\rho_s'$ satisfies $a\rho_s''/(\rho_s')^2=-au'$, hence $u'=f_*(w,u)$ and $u(p)=1/c_H=(1-rp)/(a\gamma)$, so $u=u_*$. The lower boundary condition reads $\ee^{\rho_s(0)}=1-au_*(0)$, which has a solution with $\rho_s(0)<0$ if and only if $0<au_*(0)<1$, and then $\rho_s(w)=\log(1-au_*(0))+\int_0^w\mathrm dv/u_*(v)$. Conversely this function solves \eqref{eq:stationaryproblem}.
(iii) Substituting $u^0$ gives $-\varsigma=f_*(w,u^0)$ if and only if $a\varsigma^2-b_0\varsigma-r=0$, which is the equation for $\varsigma=-m_-$, and $u^0(0)=\varsigma/r=(1-L^0_\infty)/a$ by \eqref{eq:Linfty}. Since $a\gamma^2=\theta^2/2$ and $b_0\gamma=-(\beta-r+\theta^2/2)$, the number $y=a\gamma\varsigma$ is the positive root of
\[
 y^2+\Bigl(\beta-r+\frac{\theta^2}{2}\Bigr)y-\frac{r\theta^2}{2}=0 .
\]
The left side is negative at $y=0$ and equals $\beta r>0$ at $y=r$, so $y<r$. Hence $u_*(p)-u^0(p)=(1-rp)\bigl(1/(a\gamma)-\varsigma/r\bigr)>0$. Two solutions of the same first order equation, whose right side is locally Lipschitz in $u>0$, cannot meet, so $u_*>u^0$ on $[0,p]$.
(iv) Let $\bar u(w)=(1-rw)/a$. Then $\bar u'=-r/a$ and $f_*(w,\bar u(w))=-(b_0+a)/a=(K-r)/a$, so at any point where $u_*=\bar u$ we have $(u_*-\bar u)'=K/a>0$. If $\gamma>1$, $u_*(p)=(1-rp)/(a\gamma)<\bar u(p)$, and if $\gamma=1$, $u_*(p)=\bar u(p)$ and $u_*<\bar u$ just to the left of $p$. In both cases, if $u_*$ met $\bar u$ at a largest point $w_0<p$, then $u_*-\bar u$ would vanish at $w_0$ and be negative to its right, forcing $(u_*-\bar u)'(w_0)\le0$, a contradiction. Hence $u_*<\bar u$ on $[0,p)$ and $u_*(0)<1/a$.
\end{proof}

The linear solution $u^0$ is the reciprocal slope $1/\rho_w$ of the stationary hodograph of the frozen-income problem. Indeed, \eqref{eq:stationarysingle} gives $W^0_\infty(\lambda)=r^{-1}\bigl(1-(\lambda/L^0_\infty)^{m_-}\bigr)$ for $\lambda\ge L^0_\infty$, so that $\lambda\partial_\lambda W^0_\infty=\varsigma(1/r-W^0_\infty)$.

\subsection{Convergence for long times to retirement}
\begin{lemma}\label{lem:decay}
Suppose that there are $h_0>0$ and $\varepsilon\in(0,1)$ such that $\rho(h,0)\ge\log\varepsilon$ for all $h\ge h_0$. Then there are $C,\kappa>0$ such that $|\rho_h(h,w)|\le C\ee^{-\kappa h}$ for $h\ge h_0+1$, and $\rho(h,\cdot)$ converges in $C^2([0,p])$, exponentially fast, to a solution of \eqref{eq:stationaryproblem}.
\end{lemma}

\begin{proof}
\emph{Uniform bounds.} By \cref{thm:hodo}, the constants $N$, $\overline L$, $C$ and $\underline q$ do not depend on the horizon. Hence, for $h\ge h_0$, $\log\varepsilon\le\rho(h,0)\le\rho\le\log\overline L+N\log((p+G_\infty)/G(h_0))$ and $\underline q\le\rho_w\le C/G(h_0)$. The equation is therefore uniformly parabolic for $h\ge h_0$, and the boundary nonlinearity $g(\rho)=a/(1-\ee^\rho)$ is evaluated on the fixed interval $[\log\varepsilon,\log\overline L]$, where $g'(\rho)=a\ee^\rho/(1-\ee^\rho)^2\ge a\varepsilon$. The H\"older gradient and Schauder estimates for oblique problems \citep[Chapters XII and XIII]{Lieberman1996}, applied on the time windows $[h-2,h]\times[0,p]$, give bounds on $\rho$ in $C^{2+\alpha,1+\alpha/2}([h-1,h]\times[0,p])$ that do not depend on $h\ge h_0+2$.

\emph{Decay of $\rho_h$.} The function $\vartheta=\rho_h$ solves the linear equation
\[
 \vartheta_h=A\vartheta_{ww}+\mathcal B\vartheta_w,\qquad A=\frac{a}{\rho_w^2},\quad \mathcal B=-\frac{2a\rho_{ww}}{\rho_w^3}-(1-rw),
\]
with $\vartheta_w(h,p)=0$ and $\vartheta_w(h,0)=g'(\rho(h,0))\,\vartheta(h,0)$, obtained by differentiating the boundary conditions in $h$. By the previous step, $A\ge A_0>0$ and $|\mathcal B|\le \mathcal B_1$ for $h\ge h_0+1$. Choose $\alpha_1>0$ with $A_0\alpha_1^2-\mathcal B_1\alpha_1\ge1$, then $M_\psi\ge1+\alpha_1/(a\varepsilon)$, and put $\psi(w)=M_\psi-\ee^{-\alpha_1 w}$ and $\kappa=\ee^{-\alpha_1 p}/(2M_\psi)$. Then $A\psi''+\mathcal B\psi'\le-\ee^{-\alpha_1 w}(A_0\alpha_1^2-\mathcal B_1\alpha_1)\le-\ee^{-\alpha_1 p}$, $\psi'(p)>0$ and $\psi'(0)=\alpha_1\le g'(\rho(h,0))\psi(0)$. The function $\Theta=K_0\ee^{-\kappa(h-h_0-1)}\psi$ with $K_0=\sup|\vartheta(h_0+1,\cdot)|/(M_\psi-1)$ therefore satisfies $\Theta_h-A\Theta_{ww}-\mathcal B\Theta_w>0$, $\Theta_w(h,p)>0$, $\Theta_w(h,0)\le g'\Theta(h,0)$ and $\Theta\ge|\vartheta|$ at $h=h_0+1$. If $\vartheta-\Theta$ had a positive maximum on $[h_0+1,h]\times[0,p]$, it could not be interior because of the strict differential inequality, not at $w=p$ because there $(\vartheta-\Theta)_w=-\Theta_w<0$, and not at $w=0$ because there $(\vartheta-\Theta)_w=g'(\vartheta-\Theta)+(g'\Theta-\Theta_w)>0$. Hence $\vartheta\le\Theta$, and the same argument applied to $-\vartheta$ gives $|\rho_h|\le\Theta\le C\ee^{-\kappa h}$.

\emph{Convergence.} It follows that $\rho(h,\cdot)$ converges uniformly to a limit $\rho_s$ at the rate $\ee^{-\kappa h}$. By the uniform $C^{2+\alpha}$ bounds and interpolation, the convergence holds in $C^2([0,p])$, still at an exponential rate. Passing to the limit in the equation and the boundary conditions, and using $\rho_h\to0$, shows that $\rho_s$ solves \eqref{eq:stationaryproblem}. Moreover, $\rho_s'\ge\underline q$ and $\log\varepsilon\le\rho_s(0)\le\log\overline L$.
\end{proof}

\begin{proof}[Proof of \cref{thm:longrun}]
Part (a) is \cref{lem:stationary}(ii), and the formula for $H^\infty$ follows from $\rho_s(p)=\rho_s(0)+\int_0^p\mathrm dw/u_*$.

(b) Let $\rho_s$ be the stationary solution. For every $M\ge0$ the function $\rho_s-M$ is a subsolution of the two-boundary problem: the equation is invariant under constant shifts, the condition at $w=p$ holds with equality, and at $w=0$, $(\rho_s-M)'(0)=g(\rho_s(0))\ge g(\rho_s(0)-M)$ because $g$ is increasing. Choose $M$ with $\rho_s-M\le\rho(1,\cdot)$. \Cref{lem:comparison} on $(1,h]\times[0,p]$ for every $h>1$, where no corner singularity is present, gives $\rho(h,0)\ge\rho_s(0)-M$ for $h\ge1$. \Cref{lem:decay} applies, and by uniqueness the limit is $\rho_s$. Hence $L=\ee^{\rho(h,0)}$ and $H=\ee^{\rho(h,p)}$ converge exponentially fast. In hodograph variables \eqref{eq:dB} reads, after an integration by parts,
\[
 d(\tau)=\int_0^pw\,\ee^{-\rho}\rho_w\dd w+p\,\ee^{-\rho(h,p)}=\int_0^p\ee^{-\rho(\tau-\tau_p,w)}\dd w,
\]
so $d$ and $\zb=dH-p$ converge exponentially fast as well. Since $W(\tau,\cdot)$ is the inverse of $\rho(\tau-\tau_p,\cdot)$, extended by $0$ and $p$, and $\rho_w\ge\underline q$, the functions $W(\tau,\cdot)$ converge locally uniformly, and \eqref{eq:Xrecon} and \eqref{eq:QFrecon} then give the locally uniform convergence of $\X$ and $Q$.

(c) This is \cref{lem:stationary}(iv). (It can also be derived dynamically: for $\gamma>1$, $L$ is nondecreasing by \cref{thm:horizon}(c), so \cref{lem:decay} applies and produces a stationary solution.)

(d) If $\liminf_{\tau\to\infty}L(\tau)>0$, there are $\varepsilon>0$ and $h_0$ with $\rho(h,0)\ge\log\varepsilon$ for $h\ge h_0$, and \cref{lem:decay} produces a stationary solution, which contradicts $u_*(0)\ge1/a$. By (c), this case requires $\gamma<1$.

(e) By \cref{lem:stationary}(iii), $u_*(0)>u^0(0)$, and $L^\infty-L^0_\infty=-a\bigl(u_*(0)-u^0(0)\bigr)<0$.
\end{proof}

\section{Numerical results and implications}\label{sec:numerics}

\subsection{Baseline parameters and numerical method}
The baseline parameters are listed in \cref{tab:params}. A working life of forty years corresponds to ages 25 to 65. The price $p=12.5$ means that an increase of the income rate by one unit, until retirement, costs 12.5 units of wealth. Forty years before retirement this is an investment with an internal rate of return of about $7.6\%$, well above the interest rate. With these values $K=0.03625$, $a=0.01125$ and $\tau_p=14.38$: no income is bought after age $50.6$.

\begin{table}[t]
\centering\small
\caption{Baseline parameters.}\label{tab:params}
\begin{tabular}{@{}llll@{}}
\toprule
Parameter & Symbol & Value & Implied quantity\\
\midrule
Interest rate & $r$ & $0.02$ & \\
Expected stock return, volatility & $\mu,\ \sigma$ & $0.08,\ 0.20$ & $\theta=0.30$\\
Discount rate & $\beta$ & $0.03$ & \\
Relative risk aversion & $\gamma$ & $2$ & $K=0.03625$, $a=0.01125$\\
Price of a unit of income flow & $p$ & $12.5$ & $\tau_p=14.38$\\
Time to retirement at age 25 & $T$ & $40$ & \\
\bottomrule
\end{tabular}
\end{table}

The free boundaries, value functions and policies below are computed from the integral equations of \cref{thm:integral}, for the problem with income adjustment and for the frozen-income benchmark. The stationary quantities of \cref{thm:longrun} come from the scalar problem \eqref{eq:ustar}. The system \eqref{eq:IEsys} is solved by marching in $\tau$ on a grid that is refined geometrically near $\tau=0$ and near $\tau=\tau_p$, with $L(\tau)=K\varphi(\tau)$ imposed below $\tau=10^{-8}$ as the normalization of \cref{rem:integral}(iii). The function $W$ is then evaluated from \eqref{eq:IErep}, the functions $\X$, $Q$ and $v$ from \eqref{eq:Xrecon}--\eqref{eq:vrecon}, and $d$ and $\zb$ from \eqref{eq:dformula}. For $\tau\ge20$ the computed boundaries agree with those of a finite difference solution of the hodograph problem of \cref{sec:hodo} to at least five significant digits, and at $\tau=300$ they agree with the stationary values $L^\infty$ and $H^\infty$ of \cref{thm:longrun} to six digits.

\begin{figure}[t]
\centering
\includegraphics[width=\textwidth]{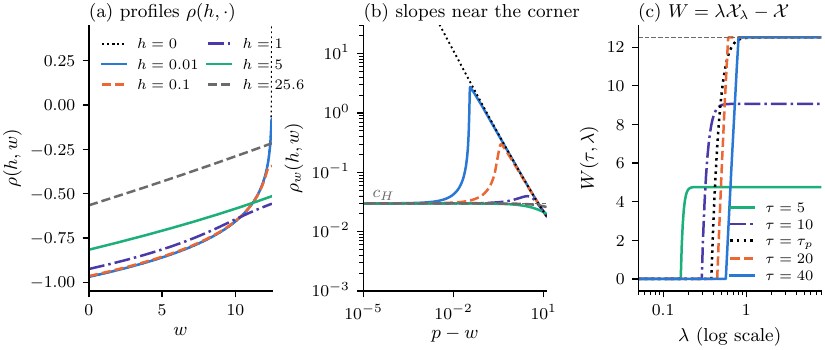}
\caption{The hodograph construction in the baseline. (a) Profiles $\rho(h,\cdot)$ on the fixed interval $[0,p]$ for several $h=\tau-\tau_p$, obtained by inverting $W(\tau,\cdot)$. The dotted curve is the initial profile $\rho_0$. (b) The slope $\rho_w=1/W_\zeta$ against the distance $p-w$ to the corner, on logarithmic scales. The dashed line is the boundary value $c_H$. (c) The function $W(\tau,\lambda)=\lambda\X_\lambda-\X$ given by \eqref{eq:IErep} for several $\tau$.}\label{fig:hodograph}
\end{figure}

\Cref{fig:hodograph} illustrates the hodograph construction. Panel (a) shows that away from the corner the profiles $\rho(h,\cdot)$ move only slowly, while near $w=p$ the solution detaches immediately from the diverging initial profile. Panel (b) shows the mechanism: for small $h$ the slope $\rho_w$ is large inside the interval and drops to $c_H$ in a layer of width comparable to $G(h)$ at $w=p$, which is the structure encoded in the barrier $C/(p-w+G(h))$. Panel (c) shows $W(\tau,\cdot)$. Because $a$ is small in the baseline, $W$ rises steeply from $0$ to its upper limit, which is $\varphi(\tau)$ for $\tau\le\tau_p$ and $p$ afterwards.

\subsection{Free boundaries and the timing of income purchases}
\begin{figure}[t]
\centering
\includegraphics[width=\textwidth]{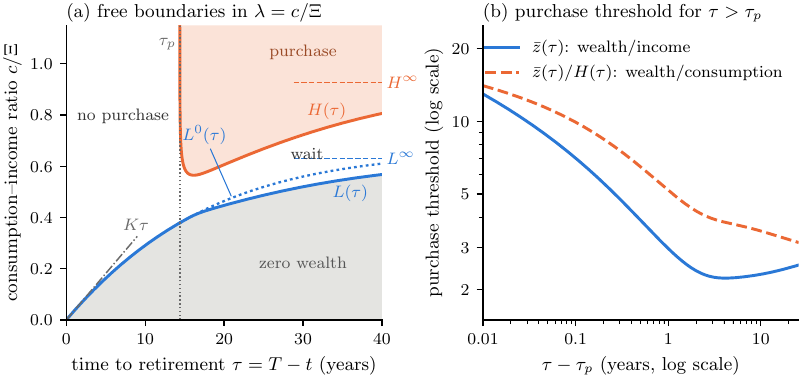}
\caption{(a) Free boundaries in the consumption--income variable $\lambda=c/\Xi$: zero-wealth boundary $L$, purchase boundary $H$, the zero-wealth boundary $L^0$ without income adjustment, the stationary levels $L^\infty$ and $H^\infty$ of the perpetual income model, and the asymptote $K\tau$. (b) The purchase threshold for $\tau>\tau_p$, in years of income ($\zb$) and in years of consumption ($\zb/H$).}\label{fig:boundaries}
\end{figure}

\Cref{fig:boundaries}(a) displays the free boundaries in the consumption--income variable. In the last $\tau_p$ years only the zero-wealth boundary $L$ is present. It starts from zero at retirement with slope $K$, as \cref{thm:single} predicts, so that an investor without savings shortly before retirement consumes a small fraction of her income and saves the rest for retirement. Beyond $\tau_p$ the purchase boundary $H$ enters from above. Between the two boundaries lies the no-action region. Below $L$ wealth is zero, above $H$ income is bought. For long horizons both boundaries approach the boundaries $L^\infty$ and $H^\infty$ of the perpetual income problem, but the approach is slow and at $\tau=40$ they are still well below these values. Since $\gamma=2>1$, $L$ is nondecreasing on the whole horizon, as \cref{thm:horizon}(c) asserts. The dotted curve is the zero-wealth boundary $L^0$ of an investor who cannot buy income. It lies above $L$: the purchase option makes a constrained investor consume less out of current income, because savings now have an additional use. \Cref{thm:compare} proves this inequality in the first years after $\tau_p$, and in the baseline the slope condition of that theorem holds for all $\tau\le40$, and at age 25, $L^0(40)=0.610$ against $L(40)=0.568$. The inequality is not universal: for $\gamma=5$ we find $L(40)=0.725>L^0(40)=0.707$, because a less substitutable investor responds more to the income effect of the option. It is restored in the long run whenever the stationary problem is solvable, in particular for every $\gamma>1$, since \cref{thm:longrun}(c),(e) then give $L^\infty<L^0_\infty$. For $\gamma=5$, $L^\infty=0.852$ and $L^0_\infty=0.868$.

Panel (b) shows the wealth threshold $\zb(\tau)$ at which income is bought, on logarithmic axes. It diverges at $\tau_p$, falls steeply over the next few years, reaches a minimum of about $2.24$ years of income around $\tau\approx18.5$, and then rises slowly to $2.53$ at $\tau=40$. The threshold measured in years of consumption, $\zb/H$, decreases monotonically. Two effects produce the U-shape. A longer time to retirement raises the value of an additional unit of income and lowers the wealth needed before buying, and this dominates near $\tau_p$. At the same time a young investor anticipates future income growth and consumes more relative to her current income ($H$ increases with $\tau$ once $\tau\gtrsim\tau_p+2$), so the same liquidity buffer, measured in consumption, corresponds to a larger multiple of current income. The U-shape is a property of the baseline parameters rather than a theorem: the general results only guarantee that $\zb$ diverges at $\tau_p$ (\cref{thm:main}) and converges to $\zb^\infty=d^\infty H^\infty-p$ as $\tau\to\infty$ (\cref{thm:longrun}).

\subsection{Consumption, portfolio and the value of income}
\begin{figure}[t]
\centering
\includegraphics[width=\textwidth]{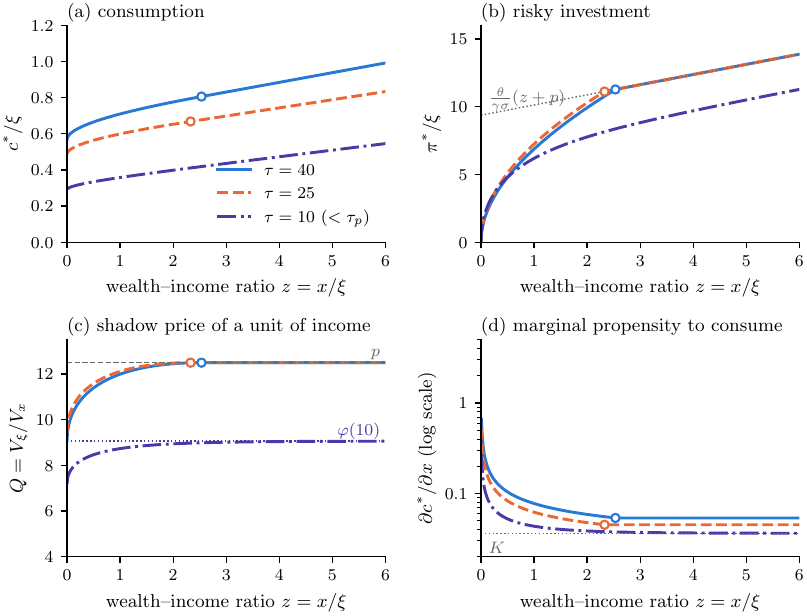}
\caption{Optimal policies against the wealth--income ratio $z=x/\xi$ at three horizons ($\tau=10$ lies in the no-purchase window). Circles mark the purchase threshold $\zb(\tau)$. (a) Consumption--income ratio. (b) Risky investment per unit of income. The dotted line is $\frac{\theta}{\gamma\sigma}(z+p)$. (c) Shadow price $Q=V_\xi/V_x$. (d) Marginal propensity to consume out of wealth, logarithmic scale.}\label{fig:policies}
\end{figure}

\Cref{fig:policies} shows the optimal policies at $\tau=40$, $25$ and $10$ as functions of the wealth--income ratio. Circles mark the purchase threshold, beyond which the investor immediately buys income. Consumption (panel a) is concave in wealth, as stated in \cref{prop:properties}, and the marginal propensity to consume (panel d) is infinite at zero wealth and falls quickly: at $\tau=40$ it is $0.22$ at $z=0.1$ and $0.077$ at $z=1$. In the purchase region it equals $1/d(\tau)$, which is $0.054$ at $\tau=40$ and exceeds the Merton value $K=0.036$: a marginal unit of wealth above the threshold is converted partly into income, which supports a permanently higher level of consumption.

Risky investment (panel b) rises like $\sqrt z$ near zero wealth and reaches the line $\frac{\theta}{\gamma\sigma}(z+p)$ at the purchase threshold. The investor leverages her liquid wealth heavily when it is small: at $\tau=40$ and $z=1$ the risky position is $6.9$ times liquid wealth. The effective human wealth is then $W=8.2$ units of income, against a market value of $\varphi(40)=27.5$ and a price of $12.5$. The liquidity constraint thus cuts the value of human capital that matters for the portfolio far below its market value. The shadow price $Q$ of a unit of income (panel c) increases with wealth and reaches $p$ at the threshold. For $\tau=10<\tau_p$ it stays below the present value $\varphi(10)=9.06$ of the remaining income, reflecting the illiquidity discount of \cref{prop:properties}(c).

\subsection{Life cycle profiles}
\begin{figure}[t]
\centering
\includegraphics[width=\textwidth]{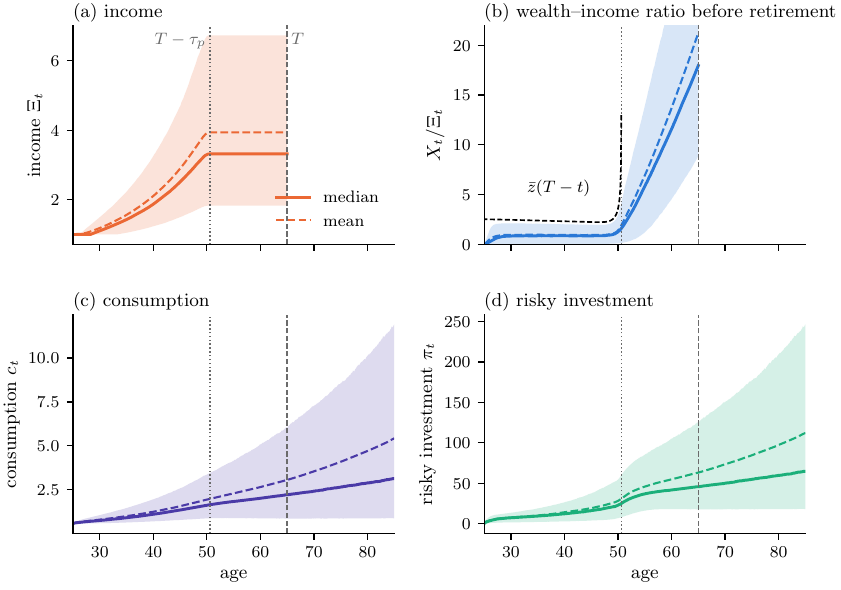}
\caption{Simulated life cycles ($2\times10^4$ paths, initial wealth $0$, initial income $1$). Solid lines are medians, dashed lines means, shaded bands the 10\%--90\% range. Vertical lines mark the last purchase date $T-\tau_p$ (age $50.6$) and retirement $T$ (age $65$). (a) Income. (b) Wealth--income ratio before retirement, with the purchase threshold $\zb(T-t)$. (c) Consumption. (d) Risky investment.}\label{fig:lifecycle}
\end{figure}

\Cref{fig:lifecycle} summarizes $2\times10^4$ simulated paths of an investor who starts at age 25 with no wealth and income $\xi=1$. The paths follow the strategy of \cref{thm:optimal}, with time step $0.01$ and the reflections at $L$ and $H$ computed from the Brownian bridge extremes of each increment. The first purchase occurs at a median age of $28.3$, and almost every path buys income at some point. Income then grows steadily, flattens, and is frozen from age $50.6$ on. At age $65$ the median income is $3.32$ and the mean $3.94$. The resulting income profile is concave and front loaded, the shape that \citet{BenPorath1967} derived from optimal human capital accumulation. Here it emerges together with a portfolio and a liquidity constraint. During the investment phase liquid wealth is kept low, between zero and the threshold $\zb$, because every unit of wealth above the threshold is converted into income (panel b). Once purchases stop, the investor accumulates financial wealth for retirement: the median wealth--income ratio rises from $1.6$ at age $50.6$ to $17.8$ at age $65$. Consumption grows smoothly over the whole life and does not jump at retirement (panel c), in line with \cref{prop:properties}(d). Risky investment (panel d) is small in absolute terms early in working life, when liquid wealth is small, and increases sharply after age $50.6$.

These profiles suggest a sequencing that is familiar from the life cycle literature but arises here from a single mechanism. Early in life the investor invests in human capital and holds little financial wealth, and later she invests in financial assets. Purchases fade out before $T-\tau_p$, because the threshold $\zb$ rises without bound as $\tau\downarrow\tau_p$, and $T-\tau_p$ is the last date at which a purchase can occur. This date is determined by the price of income and the interest rate alone.

\subsection{Long horizons and the elasticity of substitution}
\begin{figure}[t]
\centering
\includegraphics[width=\textwidth]{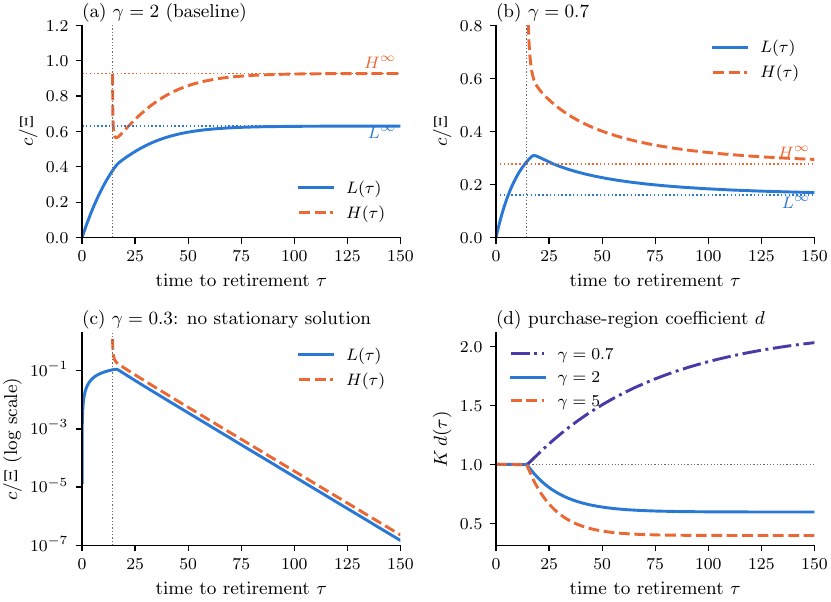}
\caption{Long times to retirement. (a) Baseline boundaries and their stationary limits (dotted). (b) $\gamma=0.7$, $\mu=0.06$, $\beta=0.04$: the zero-wealth boundary is not monotone. (c) $\gamma=0.3$, $\mu=0.04$, $\beta=0.03$: no stationary solution exists and both boundaries tend to zero (logarithmic scale). (d) The coefficient $d(\tau)$ of the purchase region, normalized by $1/K$, for $\gamma=0.7$ (parameters of panel b), $\gamma=2$ and $\gamma=5$ (other parameters as in the baseline). Vertical dotted lines mark $\tau_p$.}\label{fig:longrun}
\end{figure}

\Cref{fig:longrun} illustrates \cref{thm:horizon,thm:longrun} over times to retirement of up to 150 years, a range chosen to display the limits, not for realism. In the baseline (panel a) the boundaries approach $L^\infty$ and $H^\infty$ over several decades beyond the baseline horizon, and $L$ increases throughout. Panel (b) uses $\gamma=0.7$, $\mu=0.06$ and $\beta=0.04$. The stationary problem still has a solution ($u_*(0)=20.58<1/a=24.5$), but $L$ is no longer monotone: it rises to $0.31$ around $\tau=18$ and then declines toward $L^\infty=0.160$. Far from retirement an investor with a high elasticity of intertemporal substitution cuts consumption at zero wealth in order to exploit the purchase option. Panel (c) takes $\gamma=0.3$, $\mu=0.04$ and $\beta=0.03$, for which $u_*(0)=38.5\ge1/a=18$. There is no stationary solution, and both free boundaries decay exponentially in $\tau$, the case described in \cref{thm:longrun}(d). The purchase threshold $\zb$ nevertheless stays bounded, at about $2.65$ years of income. Panel (d) displays the income--substitution dichotomy of \cref{thm:horizon}(b): the normalized coefficient $Kd(\tau)$ equals one on the final interval, falls below one when $\gamma>1$ and rises above one when $\gamma<1$.

\subsection{Comparative statics and the value of the purchase option}
\begin{figure}[t]
\centering
\includegraphics[width=\textwidth]{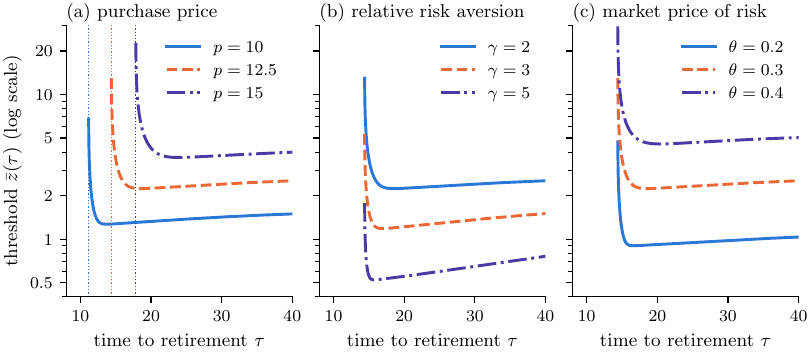}
\caption{Purchase threshold $\zb(\tau)$ (logarithmic scale) for (a) different prices $p$, with the corresponding windows $\tau_p$ as dotted lines, (b) different coefficients of relative risk aversion, and (c) different market prices of risk $\theta$ (varying $\mu$). Other parameters as in \cref{tab:params}.}\label{fig:comparative}
\end{figure}

\Cref{fig:comparative} shows how the purchase threshold depends on the price, on risk aversion and on the market price of risk. A higher price lengthens the no-purchase window and raises the threshold. At $\tau=40$, $\zb$ equals $1.50$, $2.53$ and $3.99$ for $p=10$, $12.5$ and $15$. Higher risk aversion lowers the threshold sharply ($2.53$, $1.50$ and $0.76$ for $\gamma=2,3,5$). In the model income is a riskless asset. A more risk averse investor holds less stock, so converting liquid wealth into illiquid, riskless income costs her less of the equity premium. A larger market price of risk has the opposite effect ($\zb(40)=1.04$, $2.53$ and $5.04$ for $\theta=0.2,0.3,0.4$): liquid wealth gives access to the equity premium, which raises the opportunity cost of locking wealth into income. The window $\tau_p$ does not respond to $\gamma$ or $\theta$ at all.

\begin{figure}[t]
\centering
\includegraphics[width=\textwidth]{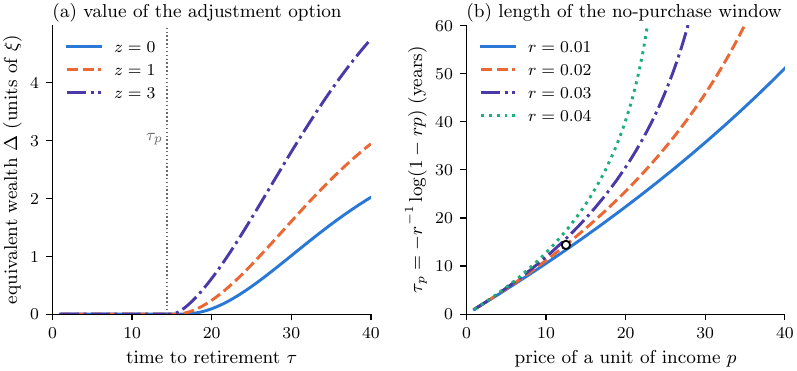}
\caption{(a) Wealth equivalent value $\Delta$ of the option to adjust income, in units of current income, for three wealth--income ratios. (b) Length of the no-purchase window $\tau_p$ as a function of the price, for several interest rates. The circle marks the baseline.}\label{fig:welfare}
\end{figure}

\Cref{fig:welfare}(a) measures the value of the option to buy income by the amount $\Delta$ of additional wealth, in units of current income, that would make an investor without the option as well off as one with it: $V^0(t,x+\Delta\xi,\xi)=V(t,x,\xi)$. The option is worthless in the last $\tau_p$ years by \cref{thm:window}, and its value rises smoothly with the remaining working life. At age 25 it is worth $2.0$ years of income to an investor without savings and $4.7$ years to an investor with $z=3$. Wealthier investors value the option more because they can exercise it sooner. Panel (b) plots $\tau_p$ as a function of the price for several interest rates. The window is roughly proportional to $p$ when $rp$ is small, since $\tau_p=p+rp^2/2+O(r^2p^3)$, and it grows without bound as $rp\uparrow1$.

The model has a simple implication for retirement policy. Raising the mandatory retirement age by one year shifts the no-purchase window by one year, so the last date at which investors find it optimal to raise their income moves with the retirement age one for one, independently of preferences and market parameters. A lower price of income shortens the window. Prices that vary with age, such as training subsidies targeted at older workers, require the extension discussed in \cref{sec:conclusion}.

\FloatBarrier
\section{Concluding remarks}\label{sec:conclusion}

We have solved a finite-horizon consumption and investment problem in which the investor can raise her income irreversibly before a mandatory retirement date and cannot borrow against future income. A pathwise dominance argument yields a no-purchase window that does not depend on preferences. The free boundary problem is then solved in three steps. During the window it is a single obstacle problem, treated by established methods. Before the window a partial hodograph transformation maps the region between the two free boundaries onto a fixed interval and turns the problem into a quasilinear parabolic equation with oblique boundary conditions, in which the entry of the purchase boundary from infinity becomes an infinite initial corner that explicit barriers can handle. The inverse transformation, glued to the obstacle problem along the opening time of the window, gives the regularity of the value function and of both free boundaries on the whole horizon. A dual budget process verifies optimality for every $\gamma\ne1$, and the logarithmic case is covered by \cref{rem:log}. The formulation on a fixed interval also yields the comparison with the frozen-income problem near $\tau_p$ and, through the decay of the time derivative of the hodograph, the description of the long-run limit in terms of the scalar function $u_*$. Finally, the Euler derivative satisfies a linear equation on the whole line with single layers on the free boundaries, which leads to Volterra integral equations for both boundaries on the whole horizon. These equations are the basis of our computations.

Several extensions seem worth pursuing. Income could be risky, correlated with the stock, or subject to depreciation. The structure based on the Euler derivative relies on the linearity of the dual equation, which can be expected to persist when income risk is spanned by the stock but not when it is unspanned. The price of income could depend on age. The pathwise comparison of \cref{thm:window} then still rules out purchases at dates where the present value of the remaining income is below the current price, but this set need not be an interval, and the structure of the purchase region outside it is an open question. Finally, making retirement a choice rather than a date, as in \citet{FarhiPanageas2007} or \citet{ChoiShim2006}, would couple the purchase boundary with an optimal stopping boundary. More generally, the combination of an obstacle problem before a critical time with a hodograph construction after it may be useful whenever a free boundary appears from infinity at a positive time.

\appendix
\section{Two explicit computations}\label{app:computations}

\paragraph{The supersolution $\Upsilon_n$.}
Write $\varpi=(\ell_n/\lambda)^{n-1}$. For $\lambda>\ell_n$,
\[
 \partial_\ell\Upsilon_n=\frac{n}{K(n-1)}(\varpi-1),\qquad \lambda\partial_\lambda\Upsilon_n=\frac1K(\lambda-\ell_n\varpi),\qquad \lambda^2\partial_{\lambda\lambda}\Upsilon_n=\frac nK\,\ell_n\varpi .
\]
Since $\ell_n'=K(n-1)/n-r\ell_n$, we have $\frac{n}{K(n-1)}\ell_n'=1-\frac{rn\ell_n}{K(n-1)}$, and collecting terms,
\[
 \PP\Upsilon_n-\lambda+1=(\varpi-1)\Bigl(1-\frac{rn\ell_n}{K(n-1)}\Bigr)-\frac{an}{K}\ell_n\varpi-\frac bK(\lambda-\ell_n\varpi)+\frac rK\Bigl(\lambda-\frac{n\ell_n}{n-1}+\frac{\ell_n\varpi}{n-1}\Bigr)-\lambda+1 .
\]
The terms linear in $\lambda$ cancel because $r-b=K$, the constant terms cancel, the term $\varpi$ survives with coefficient $1$, and the coefficient of $\ell_n\varpi$ is $(-r-an+b)/K=-(an+K)/K$. This gives the formula used in the proof of \eqref{eq:singleL}.

\paragraph{The hodograph barriers.}
For $\overline\rho$ in \eqref{eq:rhobar}, $\overline\rho_w=N/(p-w+G)$ and $\overline\rho_{ww}=N/(p-w+G)^2$, so $a\overline\rho_{ww}/\overline\rho_w^2=a/N$, and
\[
 \overline\rho_h=\frac{NG'}{p+G}-\frac{NG'}{p-w+G},\qquad
 \frac{(1-rw)-G'}{p-w+G}=\frac{r(p-w+G)}{p-w+G}=r ,
\]
which gives \eqref{eq:barresidual}. For $\overline q=C/(p-w+G)$, $\overline q_{ww}/\overline q^2=2\overline q_w^2/\overline q^3=2/(C(p-w+G))$, and
\[
 -(1-rw)\overline q_w+r\overline q=\frac{C(-1+rp+rG)}{(p-w+G)^2}=-\frac{CG'}{(p-w+G)^2}=\overline q_h .
\]
The Python script \texttt{check\_identities.py}, supplied with the paper, verifies these and the other algebraic identities used above by high order numerical differentiation at random parameter values.

\end{document}